\PassOptionsToPackage{table,xcdraw}{xcolor}

\documentclass[10pt]{amsart}

\makeatletter
\def\ifdraft{\ifdim\overfullrule>\z@
  \expandafter\@firstoftwo\else\expandafter\@secondoftwo\fi}
\makeatother

\usepackage{amsmath, amsthm, marginnote, ltxcmds, adjustbox, graphicx, stmaryrd}

\usepackage{mathtools, stackengine, scalerel}

\usepackage[nice]{nicefrac}

\usepackage[normalem]{ulem}

\usepackage{xr}
\newcommand{\crefnolink}[1]{\cref*{#1}}

\usepackage{todonotes}

\makeatletter
\newcommand*{\saved@uline}{}
\let\saved@uline\uline

\newcommand*{\mathuline}{%
  \mathpalette{\math@uline\saved@uline}%
}
\newcommand*{\math@uline}[3]{%
  \mbox{#1{$#2#3\m@th$}}%
}

\renewcommand*{\uline}{%
  \relax  
  \ifmmode
    \expandafter\mathuline%
  \else
    \expandafter\saved@uline%
  \fi
}
\makeatother

\usepackage[style=alphabetic, backref=true, backend=biber, hyperref=true, giveninits=true]{biblatex}

\AtEveryBibitem{
    \clearfield{urlyear}
    \clearfield{urlmonth}
}

\AtBeginBibliography{\small}

\renewbibmacro*{pageref}{%
	\iflistundef{pageref}
	{}
	{%
		\printtext{\addperiod\space $\uparrow${\printlist[pageref][-\value{listtotal}]{pageref}}}}}

\DeclareLabelalphaTemplate{
	\labelelement{
		\field[final]{shorthand}
		\field{label}
		\field[strwidth=1,strside=left,ifnames=1]{labelname}
		\field[strwidth=1,strside=left]{labelname}
	}
}

\definecolor{cite}{HTML}{11871E}
\definecolor{url}{HTML}{698996}
\definecolor{link}{HTML}{912F1B}

\usepackage[colorlinks=true, linkcolor=link, citecolor=cite, urlcolor=url, linktocpage]{hyperref}
\usepackage[hyphenbreaks]{breakurl}
\usepackage{xurl}
\hypersetup{breaklinks=true}

\hypersetup{final}

\usepackage{tikz}
\usetikzlibrary{cd}
\usetikzlibrary{arrows, arrows.meta, positioning, calc}
\tikzcdset{arrow style=tikz, diagrams={>={Straight Barb[scale=0.8]}}}

\tikzstyle{arrow} = [-{Straight Barb[scale=0.8]}, line width=0.2mm]

\usepackage{enumitem}
\setlist[itemize,1]{label=--}
\setlist[itemize,2]{label=+}
\setlist[itemize,3]{label=$\bullet$}
\setlist[itemize,4]{label=$\circ$}

\setlist[enumerate,1]{label=(\roman*)}

\setlist{nosep}

\usepackage[
	letterpaper,
	twoside=true,
	textheight=22cm,
	textwidth=15cm,
	marginparsep=0.75cm,
	marginparwidth=2.5cm,
	heightrounded,
	centering
]{geometry}
\usepackage[protrusion=true]{microtype}

  \usepackage[bitstream-charter]{mathdesign}
  \usepackage[T1]{fontenc}

\usepackage{biolinum}

\renewcommand{\mathsf}[1]{\text{\normalfont\sffamily#1}}

\DeclareMathAlphabet{\eur}{U}{zeus}{m}{n}
\renewcommand{\mathcal}[1]{\eur{#1}}

\usepackage[capitalise]{cleveref}

\Crefname{prop}{Proposition}{Propositions}
\Crefname{lem}{Lemma}{Lemmas}
\Crefname{cor}{Corollary}{Corollaries}
\Crefname{thm}{Theorem}{Theorems}

\Crefname{defn}{Definition}{Definitions}
\Crefname{notation}{Notation}{Notations}
\Crefname{conj}{Conjecture}{Conjectures}
\Crefname{ass}{Assumption}{Assumptions}
\Crefname{expt}{Expectation}{Expectations}

\Crefname{rmk}{Remark}{Remarks}
\Crefname{question}{Question}{Questions}
\Crefname{expl}{Example}{Examples}

\Crefname{figure}{Figure}{Figures}

\crefformat{equation}{(#2#1#3)}
\crefmultiformat{equation}{(#2#1#3)}{ \&~(#2#1#3)}{, (#2#1#3)}{, \&~(#2#1#3)}
\crefformat{section}{\S#2#1#3}
\crefmultiformat{section}{\S\S#2#1#3}{ \&~#2#1#3}{, #2#1#3}{, \&~#2#1#3}

\theoremstyle{plain}
\newtheorem{prop}[subsubsection]{Proposition}
\newtheorem{lem}[subsubsection]{Lemma}
\newtheorem{cor}[subsubsection]{Corollary}
\newtheorem{thm}[subsubsection]{Theorem}
\newtheorem*{thm*}{Theorem}

\theoremstyle{definition}
\newtheorem{defn}[subsubsection]{Definition}

\theoremstyle{remark}
\newtheorem{rmk}[subsubsection]{Remark}
\newtheorem*{rmk*}{Remark}
\newtheorem{expl}[subsubsection]{Example}

\numberwithin{equation}{subsection}

\newcommand{\teq}{\addtocounter{subsubsection}{1}\tag{\thesubsubsection}}

\DeclareMathOperator{\Ad}{\mathsf{Ad}}
\newcommand{\add}{\mathsf{add}}

\DeclareMathOperator{\Alg}{\mathsf{Alg}}

\DeclareMathOperator{\BiMod}{\mathsf{BiMod}}
\newcommand{\Br}{\mathrm{Br}}

\DeclareMathOperator{\Cat}{\mathsf{Cat}}
\newcommand{\CC}{\mathbb{C}}

\DeclareMathOperator{\cEnd}{\mathcal{E}\mathsf{nd}}
\DeclareMathOperator{\Ch}{\mathsf{Ch}}
\newcommand{\ChT}{\wtilde{\Ch}{}}
\DeclareMathOperator{\cHom}{\mathcal{H}\mathsf{om}}

\DeclareMathOperator{\cMap}{\mathcal{M}\mathsf{ap}}
\DeclareMathOperator{\Co}{\mathsf{C}}

\DeclareMathOperator{\coev}{\mathsf{coev}}
\DeclareMathOperator{\Coh}{\mathsf{Coh}}
\newcommand{\coh}{\mathsf{coh}}
\DeclareMathOperator{\coinv}{\mathsf{coinv}}
\DeclareMathOperator*{\colim}{\mathsf{colim}}
\DeclareMathOperator{\ComAlg}{\mathsf{ComAlg}}

\newcommand{\cont}{\mathsf{cont}}
\DeclareMathOperator{\Corr}{\mathsf{Corr}}
\DeclareMathOperator{\cuEnd}{\smash{\uline{\smash{\mathcal{E}\mathsf{nd}}}}}
\DeclareMathOperator{\cuHom}{\smash{\uline{\smash{\mathcal{H}\mathsf{om}}}}}

\newcommand{\defeq}{\coloneqq}
\newcommand{\DG}{\mathsf{DG}}
\DeclareMathOperator{\DGCat}{\mathsf{DGCat}}
\newcommand{\DGCatidemex}{\DGCat_{\idem,\ex}}
\newcommand{\DGCatprescont}{\DGCat_{\pres,\cont}}
\DeclareMathOperator{\Disk}{\mathsf{Disk}}

\newcommand{\En}{\mathsf{E}}
\DeclareMathOperator{\End}{\mathsf{End}}
\newcommand{\enh}{\mathsf{enh}}
\newcommand{\eqdef}{\eqqcolon}
\newcommand{\etale}{\'etale}
\DeclareMathOperator{\ev}{\mathsf{ev}}
\newcommand{\ex}{\mathsf{ex}}
\DeclareMathOperator{\Ext}{\mathsf{Ext}}

\DeclareMathOperator{\FHilb}{\mathsf{FHilb}}
\newcommand{\fin}{\mathsf{fin}}
\newcommand{\Fq}{\mathbb{F}_q}
\newcommand{\Fqn}{\mathbb{F}_{q^n}}
\newcommand{\Fqbar}{\overline{\mathbb{F}}_q}

\DeclareMathOperator{\Fun}{\mathsf{Fun}}

\DeclareMathOperator{\GL}{GL}

\newcommand{\Gm}{\mathbb{G}_\mathsf{m}}
\newcommand{\gr}{\mathsf{gr}}

\DeclareMathOperator{\Hecke}{\mathsf{H}}
\DeclareMathOperator{\HH}{\mathsf{HH}}
\DeclareMathOperator{\HHH}{\mathsf{HHH}}
\DeclareMathOperator{\Hilb}{\mathsf{Hilb}}
\newcommand{\ho}{\mathop{}\!\mathsf{h}}
\DeclareMathOperator{\Ho}{\mathsf{H}}
\DeclareMathOperator{\Hom}{\mathsf{Hom}}
\newcommand{\homflypt}{\textsc{homfly-pt}}
\newcommand{\horiz}{\mathsf{horiz}}
\newcommand{\hphBiMod}{\smash{\mhyph\BiMod}}
\newcommand{\hphMod}{\smash{\mhyph\Mod}}

\DeclareMathOperator{\id}{\mathsf{id}}
\DeclareMathOperator{\IC}{\mathrm{IC}}
\newcommand{\idem}{\mathsf{idem}}

\DeclareMathOperator{\Ind}{\mathsf{Ind}}
\DeclareMathOperator{\IndCoh}{\mathsf{IndCoh}}

\DeclareMathOperator{\inv}{\mathsf{inv}}
\newcommand{\Iso}{\mathsf{Iso}}

\DeclareMathOperator{\K}{\mathsf{K}}
\newcommand{\Kuenneth}{K\"unneth}

\DeclareMathOperator{\LS}{\mathsf{LS}}

\DeclareMathOperator{\Map}{\mathsf{Map}}
\newcommand{\mixed}{\mathsf{m}}
\DeclareMathOperator{\Mnfd}{\mathsf{Mnfd}}
\DeclareMathOperator{\Mod}{\mathsf{Mod}}

\DeclareMathOperator{\munit}{\mathsf{1}}

\newcommand{\Negut}{Neguț}
\newcommand{\nilp}{\mathsf{nilp}}

\DeclareMathOperator{\oblv}{\mathsf{oblv}}

\newcommand{\opp}{\mathsf{op}}

\newcommand{\per}{\operatorname{\mathsf{2-per}}}

\newcommand{\perf}{\mathsf{perf}}
\DeclareMathOperator{\Perf}{\mathsf{Perf}}

\newcommand{\proper}{\mathsf{prop}}

\newcommand{\pres}{\mathsf{pres}}

\newcommand{\pt}{\mathrm{pt}}

\DeclareMathOperator{\QCoh}{\mathsf{QCoh}}
\newcommand{\Qlbar}{\QQbar_{\ell}}
\def\QlbarA_#1{\QQbar_{\ell,#1}}
\newcommand{\QQbar}{\overline{\mathbb{Q}}}

\newcommand{\ren}{\mathsf{ren}}
\DeclareMathOperator{\Rep}{\mathsf{Rep}}

\DeclareMathOperator{\Res}{\mathsf{Res}}
\newcommand{\rev}{\mathsf{rev}}

\newcommand{\SBim}{\mathsf{SBim}}

\DeclareMathOperator{\sh}{\mathsf{sh}}
\DeclareMathOperator{\Shv}{\mathsf{Shv}}
\newcommand{\sign}{\mathsf{sign}}

\newcommand{\sm}{\mathsf{sm}}
\DeclareMathOperator{\Spc}{\mathsf{Spc}}
\DeclareMathOperator{\Spec}{\mathsf{Spec}}
\DeclareMathOperator{\Spr}{\mathsf{Spr}}

\newcommand{\St}{\mathsf{St}}

\DeclareMathOperator{\Stk}{\mathsf{Stk}}
\DeclareMathOperator{\Sym}{\mathsf{Sym}}
\newcommand{\SymGrp}{\mathrm{S}}

\newcommand{\theart}{\heartsuit\!_t}
\DeclareMathOperator{\Tot}{\mathsf{Tot}}
\DeclareMathOperator{\tr}{\mathsf{tr}}
\DeclareMathOperator{\Tr}{\mathsf{Tr}}
\DeclareMathOperator{\triv}{\mathsf{triv}}

\newcommand{\unip}{\mathsf{u}}

\DeclareMathOperator{\Vect}{\mathsf{Vect}}
\DeclareMathOperator{\VectgrBiMod}{\smash{\Vect^\gr\!\!\hphBiMod}}
\DeclareMathOperator{\VectgrMod}{\smash{\Vect^\gr\!\!\hphMod}}
\DeclareMathOperator{\VectgrcMod}{\smash{\Vect^{\gr,c}\!\!\hphMod}}
\newcommand{\vertc}{\mathsf{vert}}

\newcommand{\weightheart}{\heartsuit\!_w}
\newcommand{\wex}{w\mhyph\,\mathsf{ex}}
\newcommand{\wHH}{\wtilde{\mathsf{HH}}}
\DeclareMathOperator{\wt}{\mathsf{wt}}

\DeclareMathOperator{\z}{\mathsf{z}}
\DeclareMathOperator{\Z}{\mathsf{Z}}
\newcommand{\ZZ}{\mathbb{Z}}

\mathchardef\mhyphensymb="2D
\newcommand{\mhyph}{\mhyphensymb\!}

\ifdraft{
  \let\oldtimes\otimes
  \newcommand{\rotimes}{\otimes}
  
}{
  \let\oldtimes\otimes
  \renewcommand{\otimes}{\mathbin{\mkern1mu\scalerel*{\tikz[line width=.75]{
    \clip(0,0) circle[radius=1mm+.3pt];
    \draw(0,0) circle[radius=.1]; \draw (-135:.1) -- (45:.1); \draw (135:.1) -- (-45:.1);
    }}{\oldtimes}\mkern1mu}}
  
  \newcommand{\rotimes}{\mathbin{\mkern1mu\scalerel*{\tikz[line width=.75]{
    \clip(0,0) circle[radius=1mm+.3pt];
    \draw(-45:.1)--(135:.1) arc (135:-135:.1) -- (45:.1);
    }}{\otimes}\mkern1mu}}
}

\newcommand{\arrdisp}{0.33ex}

\newcommand{\alignsep}{\!\!\!\!\!\!} %

\newcommand{\interior}[1]{\oset{\circ}{#1}}
\newcommand{\lbar}[1]{\overline{#1}}

\newcommand{\lrangle}[1]{\langle#1\rangle}

\newcommand{\oset}[3][0pt]{\mathop{\mathrel{\ensurestackMath{\stackon[#1]{#3}{\scriptstyle #2}}}}}

\newcommand{\wtilde}[1]{\widetilde{#1}}

\title[HOMFLY-PT homology and Hilbert scheme of points on $\CC^2$]{Graded character sheaves, HOMFLY-PT homology, and \\ Hilbert schemes of points on $\CC^2$}
\author{Quoc P. Ho}
\address{Department of Mathematics, The Hong Kong University of Science and Technology (HKUST), Clear Water Bay, Hong Kong}
\email{phuquocvn@gmail.com}
\author{Penghui Li}
\address{YMSC, Tsinghua University, Beijing, China}
\email{lipenghui@mail.tsinghua.edu.cn}
\date{\today}


\begin{document}
\begin{abstract}
  Using a geometric argument building on our new theory of graded sheaves, we compute the categorical trace and Drinfel'd center of the (graded) finite Hecke category $\Hecke_W^\gr = \Ch^b(\SBim_W)$ in terms of the category of (graded) unipotent character sheaves, upgrading results of Ben-Zvi--Nadler and Bezrukavnikov--Finkelberg--Ostrik. In type $A$, we relate the categorical trace to the category of $2$-periodic coherent sheaves on the Hilbert schemes $\Hilb_n(\CC^2)$ of points on $\CC^2$ (equivariant with respect to the natural $\CC^* \times \CC^*$ action), yielding a proof of (a $2$-periodized version of) a conjecture of Gorsky--\Negut--Rasmussen which relates \homflypt{} link homology and the spaces of global sections of certain coherent sheaves on $\Hilb_n(\CC^2)$. As an important computational input, we also establish a conjecture of Gorsky--Hogancamp--Wedrich on the formality of the Hochschild homology of $\Hecke_W^\gr$.
\end{abstract}

\maketitle
\tableofcontents

\section{Introduction}
\subsection{Motivation}
Let $G$ be a (split) reductive group, $T \subseteq B$ a fixed pair of a maximal torus and a Borel subgroup, and $W$ the associated Weyl group. The geometric finite Hecke category associated to $G$, defined to be the category of $B$-bi-equivariant sheaves/$D$-modules on $G$, plays a central role in representation theory. Of particular interest to us are results of Bezrukavnikov--Finkelberg--Ostrik (underived) and Ben-Zvi--Nadler (derived) in~\cites{ben-zvi_character_2009, bezrukavnikov_character_2012} which identify its \emph{categorical trace} and \emph{Drinfel'd center} with the category of \emph{character sheaves} constructed by Lusztig in a series of seminal papers~\cites{lusztig_character_1985-1, lusztig_character_1985-2, lusztig_character_1985, lusztig_character_1986-1, lusztig_character_1986} as a tool for studying the characters of the finite group $G(\Fq)$.

The finite Hecke category has a \emph{graded} cousin $\Ch^b(\SBim_W)$, the monoidal ($\DG$-)category of bounded chain complexes of Soergel bimodules.\footnote{\label{ftn:SBim_Weyl_vs_Coxeter}Note that the notation is somewhat abusive here since $\SBim_W$ depends on the Coxeter system rather than just the Weyl group.} Like its geometric counterpart, (variations of) graded finite Hecke categories play an important role in geometric representation theory, especially in Koszul duality phenomena. See, for example,~\cites{beilinson_koszul_1996, bezrukavnikov_koszul_2013, achar_koszul_2018, lusztig_endoscopy_2020}. Of particular interest to us is the role these categories play in one construction of the \homflypt{} link homology. See~\cites{khovanov_triply-graded_2007, webster_geometric_2017, trinh_hecke_2021, bezrukavnikov_monodromic_2022}.

It is thus natural to study the categorical trace and Drinfel'd center of $\Ch^b(\SBim_W)$ with an eye toward \emph{both} categorified link invariants and character sheaves. Indeed, the categorical trace has been studied by Beliakova--Putyra--Wehrli and Queffelec--Rose--Sartori (underived) and Gorsky--Hogancamp--Wedrich (derived) in~\cites{beliakova_quantum_2019, queffelec_sutured_2017, queffelec_annular_2018, gorsky_derived_2021} with applications to \homflypt{} link homology. Moreover, it was proposed in~\cite{gorsky_derived_2021} that an in-depth study of the categorical trace should yield a proof of a conjecture by Gorsky--\Negut--Rasmussen~\cite{gorsky_flag_2021} which relates \homflypt{} link homology and global sections of coherent sheaves on the Hilbert schemes of points $\Hilb_n(\CC^2)$ on $\CC^2$.

Despite their importance in both categorified link homology and geometric representation theory, many questions are still left unanswered by these previous works:
\begin{enumerate}
  \item What is the Drinfel'd center of $\Ch^b(\SBim_W)$?
  \item How are the categorical trace and Drinfel'd center of $\Ch^b(\SBim_W)$ related to the category of character sheaves?
  \item \label{item:intro_questions_SBim_vs_Hilb} In type $A$, how are the trace and centers related to the category of coherent sheaves on $\Hilb_n(\CC^2)$?
\end{enumerate}
Indeed, various aspects of these questions have been raised by the experts previously but without a precise answer. See~\cite[\S1.1.5]{ben-zvi_character_2009}, \cite[discussion after Theorem A]{webster_geometry_2011}, \cite[\S1.5]{gorsky_derived_2021}, and~\cite[\S7.5]{gorsky_algebra_2022}.

In this paper, we provide complete answers to the questions above (except the center case of part \ref{item:intro_questions_SBim_vs_Hilb}) and as an application, give a proof of a version\footnote{See \cref{rmk:differences_with_GNR} for a precise comparison.} of the conjecture of Gorsky--\Negut--Rasmussen.\footnote{The link between the Drinfel'd centers and Hilbert schemes of points on $\CC^2$ requires an entirely different technique to establish. The details will appear in a forthcoming paper. We will not need this fact in the current paper. See also \cref{rmk:drinfeld_center_vs_hilb}.} Moreover, as the main computational input for answering question~\ref{item:intro_questions_SBim_vs_Hilb} above, we also establish a conjecture of Gorsky--Hogancamp--Wedrich on the formality of the Hochschild homology of $\Ch^b(\SBim_W)$,~\cite[Conjecture 1.7]{gorsky_derived_2021}, which is a consequence of the formality of the Grothendieck--Springer sheaf in our geometric setup. This formality conjecture is the main obstacle for~\cite{gorsky_derived_2021} to make the trace of $\Ch^b(\SBim_W)$ more explicit.

\subsection{The approach}
\subsubsection{Categorical traces and Drinfel'd centers}
Unlike~\cites{beliakova_quantum_2019, gorsky_derived_2021}, our computations of the categorical traces and Drinfel'd centers are completely geometric, inspired by the arguments of Ben-Zvi--Nadler found in~\cites{ben-zvi_character_2009} and our previous work~\cite{ho_eisenstein_2022}. This is possible thanks to our new theory of graded sheaves developed in~\cite{ho_revisiting_2022}, which provides a geometric realization of $\Ch^b(\SBim_W)$ (as the category $\Shv_{\gr, c}(B\backslash G/B)$ of $B$-bi-equivariant constructible \emph{graded sheaves} on $G$) \emph{within a sheaf theory} with weight structures and perverse $t$-structures and sufficient functoriality to \emph{do geometry and categorical algebra}. Such a theory is indispensable for our strategy since other geometric constructions of $\Ch^b(\SBim_W)$ via $B\backslash G/B$, such as those found in~\cites{beilinson_koszul_1996, rider_formality_2013, soergel_equivariant_2018, eberhardt_integral_2023}, have not been extended to include those spaces such as $\frac{G}{B}$ and $\frac{G}{G}$ (with conjugation actions), geometric objects that appear naturally in the study of categorical traces and Drinfel'd centers.

While the general idea for computing the categorical traces and Drinfel'd centers (see \cref{thm:intro:trace_center_vs_character_sheaves}) is similar to the one in~\cite{ben-zvi_character_2009,ho_eisenstein_2022}, the implementation is more involved. This is because compared to quasi-coherent sheaves or $D$-modules, the categorical \Kuenneth{} formula fails much more frequently for constructible sheaves. But now, with the modifications done in this paper, the proof of \cref{thm:intro:trace_center_vs_character_sheaves} works equally well with other sheaf theoretic settings, such as $D$-modules (which recovers the result of Ben-Zvi--Nadler in~\cite{ben-zvi_character_2009}) or sheaves with positive characteristic coefficients (either in the classical analytic topology when the geometric objects involved are defined over $\CC$ or in the \etale{} topology when working over $\Fq$, as long as $\ell\neq p$). We choose to work with $\ell$-adic sheaves throughout due to the availability of the theory of graded sheaves in this setting, crucial for applications to Soergel bimodules.

\subsubsection{Formality of Hochschild homologies of Hecke categories}
One interesting feature of our approach to the formality problem conjectured by Gorsky, Hogancamp, and Wedrich in~\cite{gorsky_derived_2021} is that its solution does not go through the usual ``purity implies formality'' route. This is not possible since the algebra of derived endomorphisms of the Grothendieck--Springer sheaf is not pure, already for the simplest case when $G$ is a torus. Instead, it relies on a transcendental argument carried out by the second-named author in~\cite{li_derived_2018} and a spreading argument. In type $A$, once the categories of graded unipotent character sheaves are computed explicitly using the formality result, their relation to the Hilbert schemes of points on $\CC^2$ can then be deduced as a combination of Koszul duality and results of~\cite{bridgeland_mckay_2001,krug_remarks_2018}.

\subsubsection{Hilbert schemes and link invariants}
In terms of link invariants, the categorical trace of $\Hecke^\gr_{\GL_n}$ is a universal receptacle for (derived) annular link invariants coming out of $\Hecke^\gr_{\GL_n}$. As such, \homflypt{} homology factors through it. Using an argument involving weight structures in the sense of Bondarko and Pauksztello, we prove a co-representability result for the degree-$a$ parts of the \homflypt{} homology. Matching the co-representing objects on the Hilbert scheme side, we obtain a proof of a version\footnote{See \cref{rmk:differences_with_GNR} for a precise comparison.} of a conjecture of Gorsky--\Negut--Rasmussen.

\subsubsection{The use of $\infty$-categories and higher algebra}
We remark that the use of homotopical/categorical algebra and $\infty$-categories as developed in~\cites{lurie_higher_2017-1, lurie_higher_2017, gaitsgory_study_2017} is indispensable to our approach at many levels, from the definitions/constructions of the objects and the formulation of the statements to the actual proofs. While being sophisticated, the theory is packaged in such a clean and convenient way that the arguments presented in this paper can still be followed by readers who are not familiar with it. The readers might also consult~\cite[\crefnolink{mg:sec:generalities_DGCats}]{ho_revisiting_2022} where most of the relevant background is reviewed.

\begin{rmk}
  Employing completely different techniques involving matrix factorizations, Oblomkov and Rozansky proved a 2-periodized version of~\cite{gorsky_flag_2021} in a series of papers~\cites{oblomkov_knot_2018, oblomkov_categorical_2018,oblomkov_soergel_2020}. It is not clear to us how their work fits with known results in geometric representation theory revolving around \emph{finite} (as opposed to \emph{affine}) Hecke categories and character sheaves. The appearance of affine Hecke categories and their interactions with the finite ones in their work are themselves extremely interesting and fascinating and deserve a closer look.

  Our approach, on the other hand, draws a direct connection between finite Hecke categories, character sheaves, and categorified knot invariants, where the passage from the first to the last one is as in the original construction of \homflypt{} homology theory via Soergel bimodules by Khovanov in~\cite{khovanov_triply-graded_2007}. In fact, our main result is geometric representation theoretic: we compute both the trace and the Drinfel'd center of the graded finite Hecke categories and relate them to character sheaves, upgrading previous results of Ben-Zvi--Nadler in~\cite{ben-zvi_character_2009} and Bezrukavnikov--Finkelberg--Ostrik in~\cite{bezrukavnikov_character_2012}. The relation to \homflypt{} homology is then deduced as a consequence. We note that the geometric description of the finite Hecke categories is still conjectural in their framework (see~\cite[Conjecture 7.3.1]{oblomkov_soergel_2020}).

  Their work lives in the \emph{coherent} world whereas our arguments happen in the \emph{constructible} world. Under the Langlands philosophy, there should be a duality between the objects considered in their work and ours. We expect that the two approaches are related via a graded version of Bezrukavnikov's theorem on two geometric realizations of affine Hecke algebras~\cite{bezrukavnikov_two_2016}. We plan revisit this question in a subsequent paper.
\end{rmk}

In the remainder of the introduction section, we will recall the basic objects in~\cref{subsec:notation_convention,subsec:main_players} and then provide the precise statements of the main results in~\cref{subsec:main_results}.

\subsection{Notation and conventions} \label{subsec:notation_convention}
Let us now quickly review the basic notation and conventions used throughout the paper.

\subsubsection{Category theory}
We work within the framework of $(\infty, 1)$-categories (or more briefly, $\infty$-categories) as developed by Lurie in~\cites{lurie_higher_2017-1,lurie_higher_2017}. By default, our categories are all $\infty$-categories. In particular, we use the language of $\DG$-categories as developed by Gaitsgory and Rozenblyum in this framework~\cites{gaitsgory_study_2017}. See also~\cite[\crefnolink{mg:sec:generalities_DGCats}]{ho_revisiting_2022} for a quick recap. All of our functors and categories are ``derived'' by default when it makes sense.

We let $\DGCat_{\pres,\cont}$ denote the category whose objects are presentable $\DG$-categories and whose morphisms are continuous (i.e., \emph{colimit} preserving) functors. Similarly, we let $\DGCat_{\idem,\ex}$ denote the category whose objects are idempotent complete $\DG$-categories and whose morphisms are exact (i.e., finite colimit and limit preserving) functors. Finally, let $\DGCat_{\pres,\cont,c}$ denote the \emph{1-full} subcategory (see~\cite[Volume I, Chapter 1, \S1.2.5]{gaitsgory_study_2017}) of $\DG_{\pres,\cont}$ whose objects are compactly generated and whose morphisms are quasi-proper functors, i.e., compact preserving functors. Each is equipped with a symmetric monoidal structure---the Lurie tensor product. The operation $\Ind$ of taking ind-completion is symmetric monoidal and factors as follows
\[
  \begin{tikzcd}
    \DGCat_{\idem,\ex} \ar[bend right=9]{rr}[swap]{\Ind} \ar{r}{{\Ind'}} & \DGCat_{\pres,\cont,c} \ar{r} & \DGCat_{\pres,\cont}.
  \end{tikzcd}
\]
Moreover, the first functor is an equivalence of categories.

We will informally refer to objects in $\DGCat_{\pres, \cont}$, resp. $\DGCat_{\idem,\ex}$, as ``large''/``big'', resp. ``small'', categories.

\subsubsection{Algebraic geometry} \label{subsubsec:AG_convention}
Unless otherwise specified, our stacks are Artin, with smooth affine stabilizers, and are of finite type over $\Fq$ (or $\Fqbar$, depending on the field of definition), where $q = p^k$ for some prime number $p$. We let $\Stk$ and $\Stk_{\Fq}$ denote the categories of such Artin stacks over $\Fq$ and $\Fqbar$, respectively.

We write $\pt_0 = \Spec \Fq$ and $\pt = \Spec \Fqbar$. Moreover, for any Artin stack $\mathcal{Y}$ over $\Fqbar$, we usually use $\mathcal{Y}_0$ to denote \emph{an} $\Fq$-form of $\mathcal{Y}$. We will make heavy use of the theory of graded sheaves developed in~\cite{ho_revisiting_2022}, obtained by categorically semi-simplifying Frobenius actions in the category of mixed sheaves~\cite{beilinson_faisceaux_2018}. All notation involving the theory of graded sheaves will be the same as in~\cite{ho_revisiting_2022}. We will now recall only the main pieces of notation from there.

As the theory of graded sheaves is built on top of the theory of mixed $\ell$-adic sheaves, we fix a prime number $\ell \neq p$ throughout this paper. For any Artin stack $\mathcal{Y}_0$ over $\Fq$ and $\mathcal{Y}$ its base change to $\Fqbar$, we will use $\Shv_{\mixed, c}(\mathcal{Y}_0)$, $\Shv_\mixed(\mathcal{Y}_0)$, and $\Shv_\mixed(\mathcal{Y}_0)^\ren \defeq \Ind(\Shv_{\mixed,c}(\mathcal{Y}_0))$ (resp. $\Shv_{\gr, c}(\mathcal{Y})$, $\Shv_\gr(\mathcal{Y})$, and $\Shv_\gr(\mathcal{Y})^\ren \defeq \Ind(\Shv_{\gr,c}(\mathcal{Y})^\ren)$) to denote the $\DG$-category of constructible mixed (resp. graded) sheaves, the $\DG$-category of mixed (resp. graded) sheaves, and the renormalized $\DG$-category of mixed (resp. graded) sheaves on $\mathcal{Y}_0$ (resp. $\mathcal{Y}$). The first one is an object in $\Shv_{\mixed,c}(\pt_0)\hphMod$ (resp. $\VectgrcMod$) whereas the last two are objects in $\Shv_\mixed(\pt_0)\hphMod$ (resp. $\VectgrMod$). All the usual operations on sheaves, when defined, are linear over these categories. Note that here, $\Shv_\gr(\pt) \simeq \Vect^\gr$, the ($\infty$-derived) symmetric monoidal category of graded chain complexes over $\Qlbar$ and $\Shv_{\gr,c}(\pt) \simeq \Vect^{\gr,c}$, the full symmetric monoidal subcategory consisting of perfect complexes.

$\Shv_{\gr, c}(\mathcal{Y})$ is equipped with a six-functor formalism, a perverse $t$-structure, and a weight/co-$t$-structure in the sense of Bondarko and Pauksztello. Moreover, it fits into the following diagram
\[
  \begin{tikzcd}
    \Shv_{\mixed, c}(\mathcal{Y}_0) \ar{r}{\gr_{\mathcal{Y}_0}} & \Shv_{\gr, c}(\mathcal{Y}) \ar{r}{\oblv_{\gr}} & \Shv_c(\mathcal{Y})
  \end{tikzcd}
\]
that is compatible with the six-functor formalism, the perverse $t$-structures, and the Frobenius weights. Furthermore, $\oblv_\gr \circ \gr_{\mathcal{Y}_0}$ is simply the pullback functor along $\mathcal{Y} \to \mathcal{Y}_0$. Roughly speaking, $\gr_{\mathcal{Y}_0}$ turns Frobenius weights into an actual grading and $\oblv_\gr$ forgets this grading.

For any $\mathcal{F}, \mathcal{G} \in \Shv_{\gr, c}(\mathcal{Y})$, one can talk about the graded $\cHom$-space $\cHom^\gr_{\Shv_{\gr, c}(\mathcal{Y})}(\mathcal{F}, \mathcal{G}) \in \Vect^\gr$, the category of graded (chain complexes) of vector spaces over $\Qlbar$, such that
\[
  \cHom_{\Shv_c(\mathcal{Y})}(\oblv_\gr \mathcal{F}, \oblv_\gr \mathcal{G}) \simeq \oblv_\gr \cHom^\gr_{\Shv_{\gr, c}(\mathcal{Y})}(\mathcal{F}, \mathcal{G}) \simeq \bigoplus_k \cHom^\gr_{\Shv_{\gr, c}(\mathcal{Y})}(\mathcal{F}, \mathcal{G})_k
\]
where for any $V \in \Vect^\gr$, $V_k \in \Vect$ denotes the $k$-th graded component of $V$. We also write
\[
  \cEnd^\gr_{\Shv_{\gr, c}(\mathcal{Y})}(\mathcal{F}) \defeq \cHom^\gr_{\Shv_{\gr, c}(\mathcal{Y})}(\mathcal{F}, \mathcal{F}).
\]
See~\cite[\crefnolink{mg:subsubsec:enriched_Hom}]{ho_revisiting_2022} for a quick review on enriched $\Hom$-spaces.

Since we always make use of the renormalized categories of sheaves (for example, $\Shv_{\mixed}(\mathcal{Y}_0)^\ren$ and its graded counterpart $\Shv_{\gr}(\mathcal{Y})^\ren$) introduced in~\cite{arinkin_stack_2020} and used extensively in~\cite{ho_revisiting_2022} rather than the usual categories, we will omit $\ren$ from the notation of the various pull and push functors. For example, we will simply write $f^*$, $f_*$, $f^!$, and $f_!$ rather than $f^*_\ren$, $f_{*,\ren}$, $f^!_\ren$, and $f_{!,\ren}$ as in~\cite{ho_revisiting_2022}.

\begin{rmk} \label{rmk:notation_change_vs_mixed_geometry}
  We deviate slightly from the convention used in~\cite{ho_revisiting_2022} in two places:
  \begin{enumerate}
    \item The category $\VectgrMod$ of $\Vect^\gr$-module categories is denoted by $\Mod_{\Vect^\gr}$ there.
    \item An $\Fq$-form of $\mathcal{Y}$ is denoted by $\mathcal{Y}_1$ there rather than $\mathcal{Y}_0$ as we do here. What we use in this paper conforms to the standard convention employed in, for example,~\cite{kiehl_weil_2010}. More generally, in~\cite{ho_revisiting_2022}, $\mathcal{Y}_n$ is used to denote an $\Fqn$-form of $\mathcal{Y}$. This is necessary because we have to deal with different forms of the same stack over different fields of definition. For example, see \cite[\crefnolink{mg:subsec:invariance_extensions_scalars}]{ho_revisiting_2022}.
  \end{enumerate}
\end{rmk}

\subsubsection{Grading conventions}
The different grading conventions appearing in the theory of \homflypt{} link homology can be a source of confusion (at least to the authors of the current paper). We will now collect the various grading conventions used in this paper and where they appear. Serving as a point of orientation, this subsubsection should thus be skipped, to be returned to only when the need arises.

In $\Vect^\gr$, we use $X$ to denote the formal grading and $C$ the cohomological grading. Sometimes, to emphasize the grading convention in the presence of other gradings, we also use $\Vect^{\gr_X}$ in place of $\Vect^\gr$. As most $\DG$-categories in this paper are module categories over $\Vect^\gr$, given two objects $c_1$ and $c_2$ in such a category $\mathcal{C}$, the $\Vect^\gr$-enriched $\Hom$, $\cHom^\gr_\mathcal{C}(c_1, c_2) \in \Vect^\gr$ is a graded chain complex. The gradings $X$ and $C$ therefore also apply to $\cHom^\gr_\mathcal{C}(c_1, c_2)$. We refer to this as the singly graded situation (as the cohomological grading is not a formal grading). This is the \emph{default} grading in the current paper.

The $2$-periodic construction in \cref{subsec:coh_shear_2_period} allows one to exchange an extra grading with $2$-periodization. It is used in \cref{subsec:sheared_and_2-periodic_Hilb} to introduce an extra grading on (the $\cHom^\gr$ in) any $\Vect^\gr$-module category by turning it into a $\Vect^{\gr_X, \gr_Y}$-module category. As the notation suggests, the extra grading is denoted by $Y$. Moreover, all the cohomological shearing and $2$-periodization in this paper will be with respect to this $Y$-grading. In other words, the $Y$-grading is artificially introduced and then gets ``canceled out.''

In \cref{subsubsec:X_Y_gradings_to_X'_Y'}, the $\wtilde{X}, \wtilde{Y}$-grading is introduced, which is related to the $X, Y$-grading via a linear change of coordinates. This is the grading that is used in the statements of the main theorems regarding the conjecture of~\cite{gorsky_flag_2021}. The switch is necessary to match with the usual grading convention on the \homflypt{} homology side.

Finally, on the \homflypt{} homology, the most natural gradings, from our geometric point of view, are given by $Q'$, $A'$, and $T'$, defined in \cref{subsubsec:AQT_aqt_gradings}. The relations between $Q',A',T'$, $Q, A, T$, and $q, a, t$ are also given there, where the last two grading conventions appear in~\cite{gorsky_algebra_2022}. See also \cref{rmk:difference_QAT_Q'A'T'} for the source of the difference between $Q', A', T'$ and $Q, A, T$.

\subsection{The main players}
\label{subsec:main_players}
We will now recall the definitions of the various objects/constructions that appear in the main results of the paper.

\subsubsection{The various versions of finite Hecke categories}
\label{subsubsec:different_versions_finite_Hecke}
Let $G_0$ be a split reductive group over $\Fq$, equipped with a pair $T_0 \subseteq B_0$ of a maximal torus and a Borel subgroup. Let $W$ be the associated Weyl group. By convention, $G$, $T$, and $B$ are the pullbacks of $G_0$, $T_0$, and $B_0$ to $\Fqbar$.

While our primary interest lies in the graded finite Hecke category $\Hecke^\gr_G \defeq \Shv_{\gr, c}(B\backslash G/B)$, which is equivalent to $\Ch^b(\SBim_W)$ by \cite[\crefnolink{mg:thm:geometric_incarnation_Hecke}]{ho_revisiting_2022},\footnote{See also~\cref{ftn:SBim_Weyl_vs_Coxeter}.} it is necessary to consider its ``big'' or renormalized version $\Hecke^{\gr, \ren}_G \defeq \Ind(\Hecke^\gr_G)$ because the ``big'' world possesses more functorial symmetries, such as the adjoint functor theorem~\cite[Cor. 5.5.2.9]{lurie_higher_2017-1} and the fact that compactly generated (presentable) stable $\infty$-categories are dualizable~\cite[Volume I, Proposition 7.3.2]{gaitsgory_study_2017}.\footnote{ Here, ``big'' refers to the fact that we are working with presentable categories. On the other hand, ``renormalized'' refers to the fact that we use the renormalized sheaf theory on stacks.} Consequently, we will most of the time start with the renormalized case, from which we deduce the corresponding result for the small variants.

Due to their geometric nature, our arguments apply equally well to other variants of the finite Hecke categories, such as the mixed and ungraded versions, which are of independent interest. For the reader's convenience, we summarize all the different variants in \cref{tab:Hecke_categories} below (see also~\cref{subsubsec:AG_convention} for our conventions regarding algebraic geometry).
\begin{table}[ht]
  \resizebox{\textwidth}{!}{%
    \begin{tabular}{|
        >{\columncolor[HTML]{EFEFEF}}l |l|l|}
      \hline
                                                                          &
      \multicolumn{1}{c|}{\cellcolor[HTML]{EFEFEF}\textbf{Small}}         &
      \multicolumn{1}{c|}{\cellcolor[HTML]{EFEFEF}\textbf{Big/Renormalized}}                                                                                             \\ \hline
      \textbf{Mixed}                                                      &
      $\Hecke^\mixed_{G_0} \defeq \Shv_{\mixed,c}(B_0\backslash G_0/B_0)$ &
      $\Hecke^{\mixed,\ren}_{G_0} \defeq \Shv_\mixed(B_0\backslash G_0/B_0)^\ren \defeq \Ind(\Shv_{\mixed, c}(B_0 \backslash G_0/B_0)) \simeq \Ind(\Hecke^\mixed_{G_0})$ \\ \hline
      \textbf{Graded}                                                     &
      $\Hecke^\gr_G \defeq \Shv_{\gr, c}(B\backslash G/B)$                &
      $\Hecke^{\gr,\ren}_G \defeq \Shv_\gr(B\backslash G/B)^\ren \defeq \Ind(\Shv_{\gr, c}(B\backslash G/B)) \simeq \Ind(\Hecke^\gr_G)$                                  \\ \hline
      \textbf{Ungraded}                                                   &
      $\Hecke_G \defeq \Shv_c(B\backslash G/B)$                           &
      $\Hecke^{\ren}_G \defeq \Shv(B\backslash G/B)^\ren \defeq \Ind(\Shv_{c}(B\backslash G/B)) \simeq \Ind(\Hecke_G)$                                                   \\ \hline
    \end{tabular}%
  }
  \caption{\label{tab:Hecke_categories}Different versions of finite Hecke categories.}
\end{table}

We also will adopt the notation $\Hecke^{?,\ren}_G$, $\Hecke^?_G$, $\Shv_?(\pt)$, and $\Shv_?(\mathcal{Y})$ etc. in statements where $?$ can be $\gr$, $\mixed$, or nothing/$\emptyset$, that is, the ungraded case, e.g., $\Hecke^?_G = \Hecke_G = \Shv_c(B\backslash G/B)$. In this notation, $\Hecke^{?,\ren}_G \in \Alg(\Shv_?(\pt)\hphMod)$ and $\Hecke^?_G \in \Alg(\Shv_{?,c}(\pt)\hphMod)$. Unless otherwise specified, $?$ can be any of the three.

\begin{rmk}[Abuse of notation] \label{rmk:abuse_of_notation_forget_n_mixed}
  The subscript $0$ in, for example, $G_0$, $B_0$, $T_0$, and $\pt_0$ etc. can be unwieldy and even conflicts with the use of $?$ introduced above. For example, $\Shv_?(\mathcal{Y}_0)$ is $\Shv(\mathcal{Y}_0)$ in the ungraded case, which does not really make sense, as it should be $\Shv(\mathcal{Y})$ instead in this case. Therefore, we will, in most cases, drop it altogether without fear of confusion. For example, $\Shv_\mixed(BG)^\ren$ and $\Shv_\mixed(\pt)$ only make sense when we take $BG$ and $\pt$ to be $BG_0$ and $\pt_0$, respectively. Similarly, we will use $\Hecke^\mixed_G$ and $\Hecke^{\mixed,\ren}_G$ rather than the more precise $\Hecke^\mixed_{G_0}$ and $\Hecke^{\mixed,\ren}_{G_0}$.
\end{rmk}

\begin{rmk}
  The ungraded renormalized Hecke category $\Hecke^\ren_G$ is larger than the one appearing in~\cite{ben-zvi_character_2009}. Ignoring the distinction between $D$-modules and constructible sheaves, their category corresponds to $\Shv(B\backslash G/B)$, the usual (as opposed to \emph{renormalized}) category of sheaves on $B\backslash G/B$.
\end{rmk}

\subsubsection{Categorical trace and Drinfel'd center}
As mentioned earlier, the main goal of the current paper is to study the categorical trace and Drinfel'd center of $\Hecke^\gr_G$. Recall the following definition.

\begin{defn}[{\cite[Definition 5.1]{ben-zvi_integral_2010}}] \label{defn:trace_center}
  Let $A$ be an associative algebra object in a closed symmetric monoidal $\infty$-category $\mathcal{C}$.
  \begin{enumerate}
    \item The (derived) \emph{trace} or \emph{categorical Hochschild homology} of $A$, denoted by $\Tr(A) \in \mathcal{C}$, is the relative tensor $A \otimes_{A\otimes A^\rev} A$. It comes with a natural universal trace morphism $\tr: A \to \Tr(A)$ given by $\tr(a) = a\otimes 1_A$.
    \item The (derived) \emph{Drinfel'd center} or \emph{categorical Hochschild cohomology} of $A$, denoted by $\Z(A) \in \mathcal{C}$, is the endomorphism object $\cuEnd_{A\otimes A^\rev}(A)$. It comes with a natural central morphism $\z: \Z(A) \to A$ given by $\z(F) = F(1_A)$.
  \end{enumerate}
  Here, we view $A$ as a bimodule over itself. Moreover, $A^\rev$ denotes an associative algebra with the same underlying object as $A$ but with reversed multiplications.\footnote{In the literature, $A^\opp$ is also used to denote $A^\rev$. We opt to use $A^\rev$ since for a category $A$, $A^\opp$ is usually already understood as the \emph{opposite} category.}
\end{defn}

\begin{rmk}
  The trace also has a natural $S^1$-action. Moreover, the center is naturally equipped with an $\En_2$-structure (Deligne's conjecture), i.e., in the case of categories, it has a braided monoidal structure. See~\cite[\S\S5.3 \& 5.5]{lurie_higher_2017}.
\end{rmk}

Note that $B\backslash G/B \simeq BB \times_{BG} BB$ and the monoidal structures on the various versions of finite Hecke categories are given by $g_! f^*$ (with the appropriate sheaf theory) in the following diagram
\[
  \begin{tikzcd}
    BB\times_{BG} BB \times BB \times_{BG} BB & \ar{l}{\id\times\Delta\times\id}[swap]{f} BB\times_{BG} BB \times_{BG} BB \ar{r}{g}[swap]{p_{13}} & BB \times_{BG} BB.
  \end{tikzcd} \teq\label{eq:Hecke_convolution_correspondence}
\]
Depending on whether we work with the small or big/renormalized versions, these Hecke categories are naturally algebra objects in $\Shv_{?, c}(\pt)\hphMod$ or $\Shv_?(\pt)\hphMod$, where the symmetric monoidal structures of the latter are given by \emph{relative} tensors. See also~\cite[\crefnolink{mg:subsec:module_categories,mg:subsec:large_vs_small_cats}]{ho_revisiting_2022}. Their traces and Drinfel'd centers are computed as objects in these symmetric monoidal categories. In other words, the ambient categories $\mathcal{C}$ in \cref{defn:trace_center} that are relevant to us are the various categories of module categories $\Shv_{?, c}(\pt)\hphMod$ or $\Shv_?(\pt)\hphMod$.

\subsection{The main results} \label{subsec:main_results}
We are now ready to formulate the main results of the current paper.

\subsubsection{Trace and center of Hecke categories}
We start with the computations of the trace and Drinfel'd center. The result below upgrades the main result of Ben-Zvi and Nadler \cite{ben-zvi_character_2009} to the graded setting.

\begin{thm}[\cref{thm:trace_and_center_summary}] \label{thm:intro:trace_center_vs_character_sheaves}
  The trace $\Tr(\Hecke^{?,\ren}_G)$ and center $\Z(\Hecke^{?,\ren}_G)$ of $\Hecke^{?,\ren}_G$ coincide with the full subcategory $\Ch^{\unip,?,\ren}_G$ of $\Shv_?\left(\frac{G}{G}\right)^\ren$ generated under colimits by the essential image of $\Hecke^{?,\ren}_G$ under $q_!p^*$ in the correspondence
  \[
    \begin{tikzcd}[sep=tiny]
      & \frac{G}{B} \ar{dl}[swap]{p} \ar{dr}{q} \\
      B\backslash G/B && \frac{G}{G}.
    \end{tikzcd}
  \]
  Moreover, under this identification, the canonical trace and center maps are adjoint $\tr \dashv \z$ and are identified with the adjoint pair $q_! p^* \dashv p_*q^!$.

  The trace $\Tr(\Hecke^?_G)$ coincides with the full subcategory $\Ch^{\unip,?}_G$ of $\Ch^{\unip,?,\ren}_G$ spanned by compact objects. Namely,
  \[
    \Ch^{\unip,?}_G \defeq (\Ch^{\unip,?,\ren}_G)^c = \Ch^{\unip,?,\ren}_G \cap \Shv_{?, c}\left( \frac{G}{G} \right).
  \]
  Equivalently, they are generated under finite colimits and idempotent splittings by the essential image of $\Hecke^?_G$ under $q_! p^*$.

  The center $\Z(\Hecke^?_G)$ is the full subcategory $\ChT^{\unip,?}_G$ of $\Ch^{\unip,?,\ren}_G$ spanned by the \emph{pre-image} of $\Hecke^?_G$ under the central functor $\z$.
\end{thm}

\begin{rmk}
  We note that the Drinfel'd center of $\Hecke^?_G$, that is, working within the context of small categories, is larger than its trace, unlike in the large category setting. This is true already in the case $G=T$ is a torus. See \cref{rmk:small_center_is_bigger}. This seems to be new.
\end{rmk}

\subsubsection{Formality of Hochschild homology}
By the description of the trace map above, we see that the image of the unit $\tr(1_{\Hecke^{\gr,\ren}_G}) \in \Ch^{\unip,\gr}_G$ is, by definition, the (graded) Grothendieck--Springer sheaf $\Spr^\gr_G$, a graded refinement of the usual Grothendieck--Springer sheaf, i.e., $\oblv_\gr(\Spr^\gr_G) \simeq \Spr_G \in \Ch^\unip_G \subseteq \Shv_c\left(\frac{G}{G}\right)$. By~\cite[Theorem 3.8.5]{gaitsgory_toy_2022}, the graded ring of endomorphisms of this object is the Hochschild homology $\HH(\Hecke^{\gr,\ren}_G)$  of $\Hecke^{\gr,\ren}_G$, where $\Hecke^{\gr,\ren}_G$ is \emph{viewed as a dualizable object in $\VectgrMod$}.

\begin{thm}[\cref{thm:HH_Hecke_as_End,thm:formality_Spr^gr}, {\cite[Conjecture 1.7]{gorsky_derived_2021}}] \label{thm:intro:formality_Spr^gr}
  The graded algebra $\HH(\Hecke^{\gr,\ren}_G) \simeq \cEnd^\gr_{\Ch^{\unip,\gr,\ren}_G}(\Spr^\gr_G) \in \Alg(\Vect^\gr)$ is formal. Moreover, we have an equivalence of algebras
  \[
    \cEnd^\gr_{\Ch^{\unip,\gr,\ren}_G}(\Spr^\gr_G) \simeq \Ho^*(\cEnd^\gr_{\Ch^{\unip,\gr,\ren}_G}(\Spr^\gr_G)) \simeq (\Ho_\gr^*(BT) \otimes \Ho_\gr^*(T)) \rotimes \Qlbar[W],
  \]
  where $\Ho_\gr^*(-)$ is the functor of taking graded cohomology, i.e., we remember the Frobenius weights on the cohomology (see~\cite[\crefnolink{mg:eq:graded_cohomology}]{ho_revisiting_2022}), and where the action of $\Qlbar[W]$ on the first factor is induced by the natural actions of $W$ on $T$ and $BT$.
\end{thm}

In type $A$, a result of Lusztig states that $\Ch^{\unip,\gr}_G$ is generated by $\Spr^\gr_G$. As a consequence, we obtain the following result.

\begin{thm}[\cref{thm:explicit_Char_type_A}] \label{thm:intro:explicit_Char_type_A}
  When $G$ is of type $A$, we have
  \[
    \Ch^{\unip,\gr}_G \simeq \cEnd^\gr_{\Ch^{\unip,\gr}_G}(\Spr^\gr_G)\hphMod(\Vect^\gr)^c \simeq (\Ho_\gr^*(BT) \otimes \Ho_\gr^*(T)) \rotimes \Qlbar[W] \hphMod(\Vect^\gr)^\perf,
  \]
  where for any graded ring $R \in \Alg(\Vect^\gr)$, $R\hphMod(\Vect^\gr)$ denotes the category of $R$-module objects in $\Vect^\gr$, i.e., the category of graded $R$-modules, and the superscript $\perf$  denotes the full subcategory spanned by perfect complexes (or equivalently, compact objects).
\end{thm}

\subsubsection{Hilbert scheme of points on $\mathbb{C}^2$}
While the results above are of independent interest in representation theory, their relation to the \homflypt{} link homology and the conjecture of Gorsky--\Negut--Rasmussen is one of our main motivations.

The first link is given by the following result which relates the categories of unipotent character sheaves and Hilbert schemes of points on $\CC^2$.\footnote{Strictly speaking, the Hilbert schemes that appear in our paper are over $\Qlbar$. However, we have an abstract isomorphism $\CC \simeq \Qlbar$ as fields. We will elide this inconsequential difference in the introduction section.} It is a consequence of \cref{thm:intro:explicit_Char_type_A}, Koszul duality, and a result of Krug in~\cite{krug_remarks_2018}.

\begin{thm}[\cref{thm:2-periodized_hilb_vs_ch_sheaves}] \label{thm:intro:2-periodized_hilb_vs_ch_sheaves}
  When $G=\GL_n$, we have the following equivalence of $\Vect^{\gr_{\wtilde{X}}, \gr_{\wtilde{Y}}, \per}$-module categories (thus, also as $\DG$-categories)
  \[
    {}^\Rightarrow \Ch^{\unip,\gr, \ren}_G
    \simeq \CC[\tilde{\uline{x}}, \tilde{\uline{y}}] \rotimes \CC[\SymGrp_n] \hphMod^{\gr_{\wtilde{X}},\gr_{\wtilde{Y}},\per}_{\nilp_{\tilde{\uline{y}}}}
    \xrightarrow[\simeq]{\Psi^{\per}}
    \QCoh(\Hilb_n/\Gm^2)_{\Hilb_{n, \tilde{\uline{x}}}}^{\per},
  \]
  and similarly, we have the small variant, which is an equivalence of $\Vect^{\gr_{\wtilde{X}}, \gr_{\wtilde{Y}}, \per, c}$-module categories
  \[
    {}^\Rightarrow \Ch^{\unip,\gr}_G
    \simeq \CC[\tilde{\uline{x}}, \tilde{\uline{y}}] \rotimes \CC[\SymGrp_n]\hphMod^{\gr_{\wtilde{X}}, \gr_{\wtilde{Y}}, \per}_{\nilp_{\tilde{\uline{y}}}}
    \xrightarrow[\simeq]{\Psi^{\per}} \Perf(\Hilb_n / \Gm^2)^{\per}_{\Hilb_{n, \tilde{\uline{x}}}}.
  \]
  Here,
  \begin{enumerate}
    \item $\Vect^{\gr_{\wtilde{X}}, \gr_{\wtilde{Y}}, \per}$ denotes the category of $2$-periodic bi-graded chain complexes (see \cref{defn:2_periodized_bigraded_chain} for the precise definition);
    \item $\tilde{\uline{x}}$ and $\tilde{\uline{y}}$ are multi-variables (that is, $\tilde{\uline{x}}$ denotes $\tilde{x}_1, \dots, \tilde{x}_n$ and similarly for $\tilde{\uline{y}}$), living in cohomological degree $0$ and bi-graded degrees $(1, 0)$ and $(0, 2)$ respectively;
    \item similarly, the action of $\Gm^2$ on $\Hilb_n$ is induced by its action on $\mathbb{A}^2$, scaling the two axes with weights $(1, 0)$ and $(0, 2)$, respectively;
    \item the subscript $\nilp_{\tilde{\uline{y}}}$ indicates the fact that we only consider the full subcategories where the variables $\tilde{\uline{y}}$ act (ind-)nilpotently;
    \item the subscript $\Hilb_{n,\tilde{\uline{x}}}$ indicates the fact that we only consider quasi-coherent sheaves with set-theoretic supports on the closed subscheme of $\Hilb_n$ consisting of points supported on the first axis of $\mathbb{A}^2$ (see \cref{subsubsec:supports_hilb_scheme_defn} for the precise definition); and
    \item the superscripts $\Rightarrow$ and $\per$ denote cohomological shear and $2$-periodization constructions (see \cref{subsec:coh_shear_2_period,subsec:sheared_and_2-periodic_Hilb} for an in-depth discussion on these constructions and \cref{defn:2-periodized_hilb} for the definitions of the $2$-periodized categories of sheaves on the Hilbert schemes).
  \end{enumerate}
\end{thm}

\begin{rmk} \label{rmk:drinfeld_center_vs_hilb}
  \cref{thm:intro:2-periodized_hilb_vs_ch_sheaves} above establishes the relation between the categorical trace of $\Hecke^\gr_{\GL_n}$ and the category of coherent sheaves on $\Hilb_n$ set-theoretically supported ``on the $x$-axis.'' The center, on the other hand, corresponds to the whole category of coherent sheaves on $\Hilb_n$ \emph{without any support condition}. The details will appear in a forthcoming paper.
\end{rmk}

\subsubsection{\homflypt{} link homology}
Using an argument involving weight structures in the sense of Bondarko and Pauksztello, we prove a co-representability result for the degree $a$-parts of the \homflypt{} homology. Matching the co-representing objects on the Hilbert scheme side, we obtain a proof of a version of a conjecture of Gorsky--\Negut--Rasmussen.

\begin{thm}[\cref{thm:GNR_conj_a}, {\cite[Conjecture 7.2.(a)]{gorsky_algebra_2022}}] \label{thm:intro:GNR_conj_a}
  For any $R_\beta \in \Hecke^\gr_n$ associated to a braid $\beta$, there exists a natural $\mathcal{F}_\beta \in \Perf(\Hilb_n/\Gm^2)^{\per}_{\Hilb_{n,\tilde{\uline{x}}}}$
  such that the $a$-degree $\alpha$ component of the \homflypt{} homology of $\beta$ is given by
  \[
    \cHom^{\gr_{\wtilde{X}},\gr_{\wtilde{Y}}}_{\Perf(\Hilb_n/\Gm^2)^{\per}} (\wedge^\alpha \mathcal{T}^{\per}, \mathcal{F}_\beta),
  \]
  where $\mathcal{T}$ is the tautological bundle. The two formal gradings (denoted by $\wtilde{X}$ and $\wtilde{Y}$) coming from the $\Gm^2$ action coincide with $q$ and $\sqrt{t}$. Moreover, the cohomological degree corresponds to $\sqrt{qt}$.

  See \cref{thm:GNR_conj_a} for a more precise formulation.
\end{thm}

\begin{rmk} \label{rmk:differences_with_GNR}
  \cref{thm:intro:GNR_conj_a} differs from the conjecture of Gorsky--\Negut--Rasmussen as formulated in~\cite{gorsky_algebra_2022} in two important aspects. First, we work with the $2$-periodized version of the Hilbert scheme of points on $\mathbb{C}^2$ rather than the usual version; see \cref{rmk:necessity_2-periodization} below for brief discussion of why this is necessary for us. And second, the degrees $\wtilde{X}$ and $\wtilde{Y}$ correspond to $q$ and $\sqrt{t}$ for us rather than $q$ and $t$; we believe that this is the correct torus action to consider. Moreover, compared to the conjecture of Gorsky--\Negut--Rasmussen as formulated in~\cite{gorsky_flag_2021}, our version and the version formulated in~\cite{gorsky_algebra_2022} use the usual rather than the flagged version of the Hilbert scheme; see \cref{rmk:flag_Hilb} for a speculation on how to relate the two.
\end{rmk}

\begin{rmk} \label{rmk:necessity_2-periodization}
  \cref{thm:intro:2-periodized_hilb_vs_ch_sheaves} relates the category of graded unipotent character sheaves for $\GL_n$ with the \emph{2-periodic} category of (quasi-)coherent sheaves on $\Hilb_n$, equivariant with respect to the scaling $\Gm^2$-action. As there is only \emph{one} formal grading on the character sheaf side compared to the two $\Gm$ actions on the Hilbert scheme side, $2$-periodization is necessary to relate the two. Because of this, our theorem regarding link invariants above relates \homflypt{} homology and $2$-periodic quasi-coherent sheaves on Hilbert schemes rather than the precise version in~\cite{gorsky_algebra_2022}. The same phenomenon also appeared in the work of Oblomkov and Rozansky.

  A recent work of Elias~\cite{elias_hecke_2023} introduces an extra formal grading on the Hecke category. Even though it is highly suggestive, we do not know how this helps remove the necessity of $2$-periodization.
\end{rmk}

\begin{rmk} \label{rmk:flag_Hilb}
  The main conjecture of~\cite{gorsky_flag_2021} relates $\Hecke^\gr_n$ and the category of coherent sheaves on the (derived) flag Hilbert scheme $\FHilb_n$ of $\mathbb{C}^2$ via a pair of adjoint functors. More precisely, they relate unbounded versions of these categories, which we expect to correspond to $\Hecke^{\gr, \ren}_n$ and $\IndCoh(\FHilb_n)$, respectively. We also expect that their pair of adjoint functors fits with ours in the following diagram
  \[
    \begin{tikzcd}
      \Hecke^{\gr,\ren}_n \ar[shift left=\arrdisp]{r} & \ar[shift left=\arrdisp]{l} \IndCoh(\FHilb_n)_{\FHilb_{n, \uline{x}}} \ar[shift left=\arrdisp]{r} & \ar[shift left=\arrdisp]{l} \IndCoh(\Hilb_n)_{\Hilb_{n, \uline{x}}} \simeq \QCoh(\Hilb_n)_{\Hilb_{n,\uline{x}}},
    \end{tikzcd}
  \]
  where we suppress $2$-periodization and equivariant parameter(s) altogether. Here, the compositions are our trace and central functors. Moreover, the second adjoint pair should be given by pulling and pushing along the natural map $\FHilb_n \to \Hilb_n$.

  In~\cite{elias_gaitsgorys_2018}, the category of coherent sheaves over $\FHilb_n$ is speculated to be related to a subcategory of the Hecke category generated by the Jucys--Murphy elements, which are themselves certain relative centralizers. The methods developed in (the first part of) this paper could be extended to study these relative centers geometrically. But it is unclear to us at the moment how to put \emph{all} the Jucys--Murphy elements \emph{together} in a geometric way.
\end{rmk}

\subsubsection{Support of $\mathcal{F}_\beta$}

Thanks to the Hilbert--Chow morphism, $\Hilb_n \to \mathbb{A}^{2n}\sslash \SymGrp_n$, the geometric realization of \homflypt{} link homology given in \cref{thm:intro:GNR_conj_a} implies that it admits an action of the algebra of symmetric functions in $n$ variables. More algebraically, the definition of \homflypt{} homology via Soergel bimodules also affords such an action (see, for example,~\cite[\S5.1]{gorsky_algebra_2022}). We show that these two are the same.

\begin{thm}[\cref{thm:GNR_conj_b}, {part of~\cite[Conjecture 7.2.(b)]{gorsky_algebra_2022}}] \label{thm:intro:GNR_conj_b}
  The two actions above coincide.
\end{thm}

\cref{thm:intro:GNR_conj_b} is a statement about the ``global'' property of $\mathcal{F}_\beta$ as it contains information about $\mathcal{F}_\beta$ \emph{after} pushing forward to $\mathbb{A}^{2n}\sslash \SymGrp_n$. More precisely, by~\cite[\S5.1]{gorsky_algebra_2022}, the supports of the various $a$-degree parts of the \homflypt{} homology of a braid $\beta$, as a sheaf over $\mathbb{A}^n \sslash \SymGrp_n$, can be bounded above using the number of connected components of the link associated to $\beta$.

The support of $\mathcal{F}_\beta$ over $\Hilb_n$ itself is, however, a more local (and hence, more refined) property. To prove a similar statement, we have to start by transporting the corresponding support condition on the Rouquier complex $R_\beta$. However, this does not quite make sense since there is no quasi-coherent sheaf involved at this stage yet. To circumvent this difficulty, we formulate this support condition using module category structures and the notion of support developed by Benson--Iyengar--Krause and Arinkin--Gaitsgory in~\cites{benson_local_2008, arinkin_singular_2015}. Transporting this structure around to the Hilbert scheme side, where it coincides with the usual notion of support, we obtain the desired statement.

\begin{thm}[\cref{thm:GNR_conj_b_c}, {\cite[Conjecture 7.2.(b) and (c)]{gorsky_algebra_2022}}] \label{thm:intro:GNR_conj_b_c}
  Let $\beta$ be a braid on $n$ strands. Then an upper bound for the support of $\mathcal{F}_\beta$ (as a sheaf over $\Hilb_n$) given in \cref{thm:intro:GNR_conj_a} is determined by the number of connected components of the link obtained by closing up $\beta$. See \cref{thm:GNR_conj_b_c} for a more precise formulation.
\end{thm}

\subsection{An outline of the paper}
We will now give an outline of the paper. In~\cref{sec:cat_trace_and_drinfeld_center}, we study the categorical traces and Drinfel'd centers of the various variants of the finite Hecke categories and relate them to the theory of character sheaves. In \cref{sec:formality_gspringer_sheaf}, we prove the formality result for the graded Grothendieck--Springer sheaf, of which formality of the Hochschild homologies of the graded finite Hecke categories is a consequence. In type $A$, this result gives an explicit description of the categorical trace of $\Hecke^{\gr,\ren}_n$. This is followed by \cref{sec:char_sheaves_vs_Hilb}, where the relation between unipotent character sheaves and Hilbert schemes of points on $\CC^2$ is established. Finally, the paper concludes with a proof of aversion\footnote{See \cref{rmk:differences_with_GNR} for a precise comparison.} of the conjectures by Gorsky--\Negut--Rasmussen in \cref{sec:GNR_conjecture}.

We note that the different sections are largely independent of each other, both in terms of techniques and results. More precisely, the sections can be grouped as follows: (\cref{sec:cat_trace_and_drinfeld_center}), (\cref{sec:formality_gspringer_sheaf}), (\cref{sec:char_sheaves_vs_Hilb,sec:GNR_conjecture}). The reader can read each group mostly independently, referring to the other parts of the paper mostly to look up the notation only.

\section{Categorical traces and Drinfel'd centers of (graded) finite Hecke categories} \label{sec:cat_trace_and_drinfeld_center}

This section is dedicated to the proof of the first main theorem, \cref{thm:intro:trace_center_vs_character_sheaves} (which appears as \cref{thm:trace_and_center_summary} below), which relates the categorical trace and Drinfel'd center of a finite Hecke category and the category of unipotent character sheaves. We work with mixed, graded, and ungraded versions and in both settings of ``big'' and ``small'' categories. One interesting feature is that the trace and center of a ``big'' (also known as \emph{renormalized}) finite Hecke category coincide whereas the center of the ``small'' version is larger than its trace.

In what follows, \cref{subsec:prelim_for_trace_center} reviews basic definitions and outlines the basic strategy. In \cref{subsec:categorical_Kuenneth_finite_orbit,subsec:Beck_Chevalley}, we study the categorical analog of \Kuenneth{} formula for finite orbit stacks and the Beck--Chevalley condition, respectively. As the trace/center is computed by the colimit/limit of a certain simplicial/co-simplicial category, the categorical \Kuenneth{} formula allows one to interpret each term in the (co-)simplicial diagram algebro-geometrically whereas the Beck--Chevalley condition acts as a descent-type result which allows one to realize the colimit/limit more concretely in terms of sheaves on some space. In \cref{subsec:Hecke_mixed_ren_as_functor}, we formulate the monoidal structure of finite Hecke categories in terms of $1$-manifolds, which is then used in \cref{subsec:augmented_cyclic_bar_Mnfd1} to formulate the cyclic bar simplicial category computing the trace of a finite Hecke category in terms of the geometry of a circle. The geometric formulation developed thus far is applied in \cref{subsec:adjointability} to verify the Beck--Chevalley conditions in a uniform way. The computations of the traces and centers conclude in \cref{subsec:trace_center_Hecke}. In \cref{subsec:char_sheaves}, we introduce the various versions of the category of (graded) unipotent character sheaves and state our main results thus far in terms of these categories. In \cref{subsec:wt_and_t_on_Ch}, we show that $\Tr(\Hecke^\gr_G)$ inherits the weight and perverse $t$-structure from the ambient category $\Shv_{\gr,c}(G/G)$. Finally, \cref{subsec:rigidity_etc}, which is mostly independent of the rest of the paper, deduces several consequences for the trace and center of a finite Hecke category from the rigidity of the latter.

\subsection{Preliminaries}
\label{subsec:prelim_for_trace_center}
We will now review the basic definitions of the objects involved and outline the strategy employed to prove \cref{thm:intro:trace_center_vs_character_sheaves}.

\subsubsection{Generalities regarding big vs. small categories}
In the situation of \cref{defn:trace_center}, the trace $\Tr(A) = A\otimes_{A\otimes A^\rev} A$ can be computed as the geometric realization (i.e., colimit) of the simplicial object obtained from the cyclic bar construction. Namely,
\[
  \Tr(A) = |A^{\otimes (\bullet + 1)}| \defeq \colim A^{\otimes(\bullet +1)},
\]
where the last face map is given by multiplying the last and first factors of $A$ and where the other face maps are given by multiplying adjacent $A$ factors. See also~\cite[\S5.1.1]{ben-zvi_integral_2010} and~\cite[Remark 5.5.3.13]{lurie_higher_2017}.

The small and big category settings are related by taking $\Ind$ via the following standard result.

\begin{prop}[{\cite[Corollary 3.6 and Proposition 3.3]{ben-zvi_character_2009}}] \label{prop:big_vs_small}
  Let $\mathcal{A} \in \Alg(\DGCat_{\pres,\cont})$ be a compactly generated rigid monoidal category (see~\cite[Definition 3.1]{ben-zvi_character_2009}). Taking $\Ind$-completion (resp. passing to right adjoints, resp. passing to the full subcategory $(-)^c$ spanned by compact objects) induces a canonical equivalence \ref{item:prop:big_vs_small:small}$\to$\ref{item:prop:big_vs_small:big_quasi_proper} (resp. \ref{item:prop:big_vs_small:big_quasi_proper}$\to$\ref{item:prop:big_vs_small:right_adjoint}, resp. \ref{item:prop:big_vs_small:big_quasi_proper}$\to$\ref{item:prop:big_vs_small:small}), where
  \begin{enumerate}
    \item\label{item:prop:big_vs_small:small} the category $\mathcal{A}^c\hphMod$ of $\mathcal{A}^c$-module categories in $\DGCat_{\idem,\ex}$, where $\mathcal{A}^c$ is full subcategory of $\mathcal{A}$ spanned by compact objects;
    \item\label{item:prop:big_vs_small:big_quasi_proper} the category $\mathcal{A}\hphMod_c$ of $\mathcal{A}$-module categories in $\DGCat_{\pres,\cont,c}$;
    \item\label{item:prop:big_vs_small:right_adjoint} the opposite of the category whose objects are $\mathcal{A}$-module categories in the category of compactly generated $\DG$-categories and whose morphisms are functors that are both co-continuous (i.e., limit preserving) and continuous.
  \end{enumerate}
\end{prop}

As mentioned earlier, we will work mostly in the setting of \cref{prop:big_vs_small}.\ref{item:prop:big_vs_small:big_quasi_proper} even though we are most interested in the ``small'' setting described in \cref{prop:big_vs_small}.\ref{item:prop:big_vs_small:small} because the former has more functoriality.

\begin{rmk}
  This result was also proved (but not explicitly stated) in~\cite[Volume I, Chapter 1, \S9]{gaitsgory_study_2017}. Note also that~\cite{ben-zvi_character_2009} proved the result above more generally for semi-rigid monoidal categories. We do not need this level of generality here.
\end{rmk}

\begin{rmk}
  \label{rmk:rigid_for_small}
  In the literature, the property of being rigid is usually applied to a ``small'' monoidal $\DG$-category. By definition, a category $\mathcal{C}_0 \in \Alg(\DGCatidemex)$ is \emph{rigid} if all objects $c\in \mathcal{C}_0$ have left and right duals. By~\cite[Volume I, Chapter 1, Lemma 9.1.5]{gaitsgory_study_2017}, such a category $\mathcal{C}_0$ is rigid if and only if $\mathcal{C} \defeq \Ind(\mathcal{C}_0)$ is rigid in the sense of~\cite[Volume I, Chapter 1, Definition 9.1.2]{gaitsgory_study_2017}. In this case, $\mathcal{C}$ is compactly generated by definition.

  We use the term rigid for both ``big'' and ``small'' categories. The context should make it clear in which sense we use the term.
\end{rmk}

\subsubsection{Computing colimits}
Recall that $\Ind: \DGCat_{\idem,\ex} \to \DGCat_{\pres,\cont}$ preserves all colimits by~\cite[Volume I, Chapter 1, Corollary 7.2.7]{gaitsgory_study_2017} and that the forgetful functor $\mathcal{A}\hphMod \to \DGCat_{\pres,\cont}$ preserves all colimits by~\cite[Corollaries 4.2.3.5]{lurie_higher_2017}. \cref{prop:big_vs_small} then implies that the colimit $\colim_{i\in I} \mathcal{C}_i^0$ of a diagram $I \to \mathcal{A}^c\hphMod$ can be computed as
\[
  \colim_{i\in I} \mathcal{C}_i^0 \simeq (\colim_{i\in I} \Ind(\mathcal{C}_i^0))^c,
\]
where $(-)^c$ denotes the procedure of taking the full subcategory spanned by compact objects. Moreover, the category underlying the colimit on the right-hand side can be computed in $\DGCat_{\pres,\cont}$.\footnote{Note, however, that colimits in the presentable setting are different from colimits in the category of all $\infty$-categories.}

\subsubsection{Computing limits}
Recall that the forgetful functors $\mathcal{A}\hphMod \to \DGCatprescont$ and $\DGCat_{\pres,\cont} \to \Cat$ preserves all limits by~\cite[Corollary 4.2.3.3]{lurie_higher_2017} and~\cite[Proposition 5.5.3.13]{lurie_higher_2017-1}, where $\Cat$ denotes the category of all $(\infty, 1)$-categories. The underlying category of a limit in $\mathcal{A}\hphMod$ can thus be computed in the category of all categories.

However, $\Ind$ does not preserve limits. In fact, this is the reason why the renormalized category of sheaves on a stack differs from the usual category. Moreover, as we shall see, this will also be responsible for the fact that the Drinfel'd center of the ``small'' version of a finite Hecke category differs from its trace.

\subsubsection{The strategy}
The trace of $\Hecke^{?,\ren}_G$ as an algebra object in $\Shv_?(\pt)\hphMod$ is given by the geometric realization of the cyclic bar construction
\[
  \Tr(\Hecke^{?,\ren}_G)_\bullet \defeq (\Hecke^{?,\ren}_G)^{\otimes_{\Shv_?(\pt)} (\bullet+1)} \simeq \Hecke^{?,\ren}_G \otimes_{\Shv_?(\pt)} \Hecke^{?,\ren}_G \otimes_{\Shv_?(\pt)} \cdots \otimes_{\Shv_?(\pt)} \Hecke^{?,\ren}_G \in (\Shv_?(\pt)\hphMod)^{\Delta^\opp}. \teq\label{eq:cyclic_bar_H_mixed_ren}
\]
In the next two subsections, we will discuss the two main technical ingredients used to understand the geometric realization of the simplicial object above. First, the categorical \Kuenneth{} formula for finite orbit stacks allows one to realize the terms in the simplicial object above as the categories of sheaves on certain stacks. Second, the Beck--Chevalley condition on a simplicial category is essentially a descent-type result that allows one to realize its geometric realization as a full subcategory of a more concrete category.

By a duality argument, we show that the Drinfel'd center of $\Hecke^{?, \ren}_G$ coincides with its trace. Finally, the traces and Drinfel'd centers of the small versions $\Hecke^?_G$ are obtained from those of the renormalized versions.

\subsection{Categorical \Kuenneth{} formula for finite orbit stacks}
\label{subsec:categorical_Kuenneth_finite_orbit}
This subsection is dedicated to the following result.

\begin{prop}\label{prop:categorical_Kuenneth_finite_orbit}
  Let $\mathcal{Y}_0$ and $\mathcal{Z}_0$ be Artin stacks over $\pt_0$. Suppose that $\mathcal{Y}_0$ is a finite orbit stack. Then, the external tensor product
  \[
    \Shv_?(\mathcal{Y})^\ren \otimes_{\Shv_?(\pt)} \Shv_?(\mathcal{Z})^\ren \xrightarrow[\simeq]{\boxtimes} \Shv_?(\mathcal{Y} \times \mathcal{Z})^\ren \teq\label{eq:box_tensor_mixed_ren}
  \]
  induces an equivalence of categories.
\end{prop}
\begin{proof}
  We will treat the mixed case, i.e., $? = \mixed$. The graded and ungraded cases can be treated similarly.

  By~\cite[(F.18)]{arinkin_stack_2020}, we know that \cref{eq:box_tensor_mixed_ren} is fully faithful. By~\cite[\crefnolink{mg:cor:compact_generators_relative_tensors}]{ho_revisiting_2022}, we know that the left-hand side of \cref{eq:box_tensor_mixed_ren} is generated by objects of the form $\mathcal{F}_0 \boxtimes \mathcal{G}_0$ where $\mathcal{F}_0$ and $\mathcal{G}_0$ are constructible. Since \cref{eq:box_tensor_mixed_ren} is continuous, it suffices to show that $\mathcal{F}_0 \boxtimes \mathcal{G}_0 \in \Shv_\mixed(\mathcal{Y}_0\times_{\pt_0} \mathcal{Z}_0)^\ren$ generates the category. Using excision and the fact that $\mathcal{Y}_0$ is a finite orbit stack, we reduce to the case where $\mathcal{Y}_0 = BH_0$ for a (smooth) algebraic group $H_0$ over $\Fq$. This case is treated in \cref{lem:categorical_Kuenneth_classifying} below.
\end{proof}

\begin{lem} \label{lem:categorical_Kuenneth_classifying}
  \cref{prop:categorical_Kuenneth_finite_orbit} holds in the case $\mathcal{Y}_0 = BH_0$, where $H_0$ is a smooth algebraic group over $\pt_0$.
\end{lem}
\begin{proof}
  It suffices to show that
  \[
    \Shv_{\mixed, c}(BH_0) \otimes_{\Shv_{\mixed, c}(\pt_0)} \Shv_{\mixed, c}(\mathcal{Z}_0) \xrightarrow{\boxtimes} \Shv_{\mixed, c}(BH_0 \times_{\pt_0} \mathcal{Z}_0)
  \]
  is an equivalence of categories.

  By descent along the surjective smooth map $\mathcal{Z}_0 \to BH_0 \times_{\pt_0} \mathcal{Z}_0$ (using the $(-)^!$-functor), we have
  \[
    \Shv_{\mixed,c}(BH_0 \times_{\pt_0} \mathcal{Z}_0) \simeq \Tot(\Shv_{\mixed, c}(H_0^{\times_{\pt_0}\bullet}\times_{\pt_0} \mathcal{Z}_0)).
  \]
  Then, \cite[Proposition 4.7.5.1]{lurie_higher_2017} provides an equivalence of categories
  \[
    \Shv_{\mixed, c}(BH_0 \times_{\pt_0} \mathcal{Z}_0) \simeq ((\pi_0)_!(\pi_0)^!\Qlbar)\hphMod(\Shv_{\mixed,c}(\mathcal{Z}_0)) \teq\label{eq:identification_Shv(ZxBH)}
  \]
  where $\pi_0: H_0 \to \pt_0$ and where $(\pi_0)_! (\pi_0)^! \Qlbar \in \Alg(\Shv_{\mixed, c}(\pt_0))$. Note that here, the algebra structure on $(\pi_0)_! (\pi_0)^! \Qlbar$ is obtained from the multiplication structure of $H$: indeed, $(\pi_0)_! (\pi_0)^! \Qlbar$ is the \emph{homology} $\Co_*(H) \defeq \pi_! \pi^! \Qlbar$ of $H$ along with a Frobenius action, where $\pi: H \to \pt$ is the pullback of $\pi_0$. \cref{eq:identification_Shv(ZxBH)} can thus be rewritten in the following more familiar form
  \[
    \Shv_{\mixed, c}(BH_0 \times_{\pt_0} \mathcal{Z}_0) \simeq \Co_*(H)\hphMod(\Shv_{\mixed,c}(\mathcal{Z}_0)).
  \]
  Note also that the action of $\Shv_{\mixed, c}(\pt_0)$ on $\Shv_{\mixed, c}(\mathcal{Z}_0)$ is used to make sense of the RHS of~\cref{eq:identification_Shv(ZxBH)}.

  Arguing similarly, we see that (see also~\cite[Volume I, Chapter 1, Proposition 8.5.4]{gaitsgory_study_2017})
  \begin{align*}
    \Shv_{\mixed, c}(BH_0) \otimes_{\Shv_{\mixed,c}(\pt_0)} \Shv_{\mixed, c}(\mathcal{Z}_0)
     & \simeq \Co_*(H)\hphMod(\Shv_{\mixed, c}(\pt_0)) \otimes_{\Shv_{\mixed, c}(\pt_0)} \Shv_{\mixed, c}(\mathcal{Z}_0) \\
     & \simeq \Co_*(H)\hphMod(\Shv_{\mixed, c}(\mathcal{Z}_0)).
  \end{align*}
  The proof thus concludes.
\end{proof}

\subsection{Beck--Chevalley condition}
\label{subsec:Beck_Chevalley}
We will now turn to the Beck--Chevalley condition, a technical condition that allows one to realize geometric realizations (i.e., colimits of simplicial objects) and totalizations (i.e., limits of co-simplicial objects) of categories in more concrete terms.

\begin{defn}[{\cite[Definition 4.7.4.13]{lurie_higher_2017}}]
  Given a diagram of $\infty$-categories
  \[
    \begin{tikzcd}
      \mathcal{C} \ar{r}{H} \ar{d}{V} & \mathcal{D} \ar{d}{V'} \\
      \mathcal{C}' \ar{r}{H'} & \mathcal{D}'
    \end{tikzcd}
  \]
  which commutes up to a specified equivalence $\alpha: V'\circ H \xrightarrow{\simeq} H'\circ V$.

  We say that this diagram is \emph{horizontally left} (resp. \emph{right}) \emph{adjointable} if $H$ and $H'$ admit left (resp. right) adjoints $H^L$ and $H'^R$ (resp. $H^R$ and $H'^R$), respectively, and if the composite transformation
  \begin{align*}
    H'^L \circ V' \to H'^L \circ V' \circ H \circ H^L \xrightarrow[\simeq]{\alpha} H'^L \circ H' \circ V \circ H^L \to V \circ H^L \\
    (\text{resp. } \quad\qquad V \circ H^R \to H'^R \circ H' \circ V \circ H^R \xrightarrow[\simeq]{\alpha^{-1}} H'^R \circ V' \circ H \circ H^R \to H'^R \circ V')
  \end{align*}
  is an equivalence.

  We say that this diagram is \emph{vertically left} (resp. \emph{right}) \emph{adjointable} if the transposed diagram
  \[
    \begin{tikzcd}
      \mathcal{C} \ar{r}{V} \ar{d}{H} & \mathcal{C}' \ar{d}{H'} \\
      \mathcal{D} \ar{r}{V'} & \mathcal{D}'
    \end{tikzcd}
  \]
  is \emph{horizontally left} (resp. \emph{right}) \emph{adjointable}.
\end{defn}

\begin{rmk}
  What we call horizontally left (resp. right) adjointable is simply called left (resp. right) adjointable in~\cite{lurie_higher_2017}. This condition on a commutative square of categories is also commonly called the Beck--Chevalley condition.
\end{rmk}

We are now ready to state the main statement of this subsection.
\begin{prop}[Lurie] \label{prop:Beck-Chevalley_fully_faithful}
  Let $\mathcal{C}_\bullet: \Delta_+^\opp \to \DGCat_{\pres,\cont}$ be an augmented simplicial diagram of $\DG$-categories. Let $\mathcal{C} = \mathcal{C}_{-1}$, $F: \mathcal{C}_0 \to \mathcal{C}$ be the obvious functor, and $G$ its right adjoint. Suppose that for every morphism $\alpha: [n] \to [m]$ in $\Delta_+$ which induces a morphism $\alpha_+: [n+1] \simeq [0]*[n] \to [0]*[m] \simeq [m+1]$, the diagram
  \[
    \begin{tikzcd}
      \mathcal{C}_{m+1} \ar{d} \ar{r}{d_0} & \mathcal{C}_m \ar{d} \\
      \mathcal{C}_{n+1} \ar{r}{d_0} & \mathcal{C}_n
    \end{tikzcd} \teq\label{eq:Beck-Chevalley_vert_right_adj}
  \]
  is vertically right adjointable. Then, the functor $\theta: |\mathcal{C}_\bullet|_{\Delta^\opp}| \to \mathcal{C}$ is fully faithful. When $G$ is conservative, then $\theta$ is an equivalence of categories.
\end{prop}
\begin{proof}
  We will deduce this proposition from its dual version~\cite[Corollary 4.7.5.3]{lurie_higher_2017}.

  Passing to right adjoints, we get an augmented co-simplicial object $\mathcal{C}^\bullet$. Note that for each $n$, $\mathcal{C}^n \simeq \mathcal{C}_n$ as only the functors change. By~\cite[Corollary 5.5.3.4]{lurie_higher_2017-1} (see also~\cite[Volume I, Chapter 1, Proposition 2.5.7]{gaitsgory_study_2017}), we have a canonical equivalence of categories
  \[
    |\mathcal{C}_\bullet|_{\Delta^\opp}| \simeq \Tot(\mathcal{C}^\bullet|_{\Delta}).
  \]
  Moreover, the canonical functor $\mathcal{C} \to \Tot(\mathcal{C}^\bullet|_{\Delta})$ is the right adjoint $\theta^R$ of $\theta$.

  Observe that
  \[
    \begin{tikzcd}
      \mathcal{C}^{m+1} & \mathcal{C}^m \ar{l}[swap]{d^0} \\
      \mathcal{C}^{n+1} \ar[u] & \mathcal{C}^n \ar{l}[swap]{d^0} \ar{u}
    \end{tikzcd}
  \]
  is horizontally left adjointable as this is equivalent to~\cref{eq:Beck-Chevalley_vert_right_adj} being vertically right adjointable. Thus, \cite[Corollary 4.7.5.3]{lurie_higher_2017} implies that $\theta$ is a fully faithful embedding, and moreover, it is an equivalence of categories when $G$ is conservative.
\end{proof}

\begin{rmk}
  The same proof goes through when we replace $\DGCat_{\pres,\cont}$ by $\mathcal{A}\hphMod$ where $\mathcal{A}$ is a rigid monoidal category, for example, when $\mathcal{A}$ is $\Shv_\mixed(\pt_0)$ or $\Vect^\gr \simeq \Shv_\gr(\pt)$. This is because the forgetful functor $\mathcal{A}\hphMod \to \DGCatprescont$ preserves all limits and colimits, by~\cite[Corollaries 4.2.3.3 and 4.2.3.5]{lurie_higher_2017}, and adjoints of an $\mathcal{A}$-linear functor are automatically $\mathcal{A}$-linear by~\cite[\crefnolink{mg:cor:lax_implies_strict_rigid}]{ho_revisiting_2022}.
\end{rmk}

\subsection{\texorpdfstring{$\Hecke^{?,\ren}_G$}{Hecke} as a functor out of $1$-manifolds}
\label{subsec:Hecke_mixed_ren_as_functor}

While it is straightforward how to apply \cref{prop:categorical_Kuenneth_finite_orbit}, the verification of the adjointability conditions (a.k.a. Beck--Chevalley conditions) necessary for the application of \cref{prop:Beck-Chevalley_fully_faithful} is more subtle. The argument can be made transparent by reformulating the simplicial object \cref{eq:cyclic_bar_H_mixed_ren} geometrically in terms of $1$-manifolds. In preparation for that, we will, in this subsection, upgrade $\Hecke^{?,\ren}_G$ to a right-lax symmetric monoidal functor coming out of the category of $1$-manifolds. The construction found in this subsection is a simplification of the one found in~\cite{ho_eisenstein_2022}.

\subsubsection{$1$-manifolds and their boundaries}
\label{subsubsec:1-mnfd_and_boundaries}
Let $\Mnfd_1$ denote the symmetric monoidal $\infty$-category of  $1$-dimensional manifolds with finitely many connected components whose morphisms are given by embeddings and whose symmetric monoidal structure is given by disjoint unions. Note that the objects are simply finite disjoint unions of lines and circles. Although these manifolds are, technically speaking, without boundary, we will make use of their ``boundaries'' in our construction. We will now explain what this means.

Let $\Mnfd_1'$ be the $\infty$-category of compact $1$-manifolds, usually denoted by $M$, with possibly non-empty boundary, denoted by $\partial M$, such that $\interior{M} \defeq M \setminus \partial M \in \Mnfd_1$. Moreover, morphisms are given by (necessarily closed) embeddings. Taking the interior gives a natural functor of $\infty$-categories
\[
  F: \Mnfd_1' \to \Mnfd_1, \qquad M \mapsto \interior{M}.
\]

\begin{lem}[{\cite[Lemma 2.4.10]{ho_eisenstein_2022}}]
  \label{lem:boundary_or_not}
  The functor $F$ is an equivalence of categories. We write $M \mapsto \lbar{M}$ to denote an inverse of $F$.
\end{lem}

Because of this equivalence, we will generally not make a distinction between $1$-manifolds without boundaries and compact $1$-manifolds with possibly non-empty boundaries, unless confusion is likely to occur. For instance, when $M \in \Mnfd_1$, by abuse of notation,
\[
  \partial M \defeq \partial \lbar{M} \defeq \lbar{M} \setminus M
\]
is used to denote the boundary of $\lbar{M}$.

We also use $\Disk_1$ to denote the full subcategory of $\Mnfd_1$ consisting of just lines. Moreover, the category $\Disk'_1 \defeq \Disk_1 \times_{\Mnfd_1} \Mnfd'_1$ consisting of compact line segments is equivalent to $\Disk_1$ itself.

\subsubsection{$\Hecke^{?,\ren}_G$ as a functor out of $\Mnfd_1$}
Following~\cite[\S3.1.2]{ho_eisenstein_2022}, we will now construct a right-lax symmetric monoidal functor
\[
  \Hecke^{?,\ren}_G: \Mnfd_1 \to \Shv_?(\pt)\hphMod,
\]
whose restriction to $\Disk_1$ is symmetric monoidal. The functor $\Hecke^{?,\ren}_G$ is obtained as a composition\footnote{\label{ftn:Stk_vs_StkFq} Even though we write $\Stk$, it should be understood as $\Stk_{\Fq}$ when we consider the mixed variant.}
\[
  \begin{tikzcd}
    \Mnfd_1 \ar{r}{M} & \Corr(\Stk)_{\proper;\sm} \ar{r}{\Shv^{*,\ren}_{?,!}} & \Shv_?(\pt)\hphMod.
  \end{tikzcd} \teq\label{eq:Hecke_as_a_functor}
\]

We will now explain the various functors and categories that appear in the diagram above. We start with the functor $M$, which is defined as the composition
\[
  \begin{tikzcd}
    \Mnfd_1 \ar{r}{\mathfrak{B}} & \Corr((\Spc^{\Delta^1}_\fin)^\opp) \ar{r}{\cMap} & \Corr(\Stk).
  \end{tikzcd} \teq\label{eq:M_as_composition}
\]

\subsubsection{Correspondences}
For any category $\mathcal{C}$, $\Corr(\mathcal{C})$ is the category of correspondences in $\mathcal{C}$. A morphism from $c_1$ to $c_2$ in $\Corr(\mathcal{C})$ is illustrated, equivalently, by diagrams of the form
\[
  \begin{tikzcd}
    c_1 & \ar{l}[swap]{h} c \ar{d}{v} \\
    & c_2
  \end{tikzcd} \qquad\text{or}\qquad
  \begin{tikzcd}[column sep=small]
    & \ar{dl}[swap]{h} c \ar{dr}{v} \\
    c_1 && c_2
  \end{tikzcd} \qquad\text{or}\qquad
  \begin{tikzcd}
    c_1 & \ar{l}[swap]{h} c \ar{r}{v} & c_2.
  \end{tikzcd}
\]
Here, $h$ and $v$ stand for horizontal and vertical respectively. As usual, compositions are given by pullbacks.

More generally, let $\vertc$ and $\horiz$ be two collections of morphisms in $\mathcal{C}$ such that $\vertc$ (resp. $\horiz$) is closed under pulling back along a morphism in $\horiz$ (resp. $\vertc$). Then, we let $\Corr(\mathcal{C})_{\vertc;\horiz}$ be the $1$-full subcategory of $\Corr(\mathcal{C})$ consisting of the same objects but for morphisms, we require that $v\in \vertc$ and $h\in \horiz$. See~\cite[\crefnolink{mg:subsubsec:correspondences}]{ho_eisenstein_2022} for more details.

\subsubsection{The functor $\mathfrak{B}$} \label{subsubsec:B_functor}
When $\mathcal{C} = \Spc_\fin$, the category of finite CW-complexes, the functor\footnote{$\mathfrak{B}$ stands for \emph{boundary}.}
\[
  \mathfrak{B}: \Spc_\fin \to \Corr((\Spc_\fin^{\Delta^1})^\opp),
\]
at the level of objects, is given by
\[
  \mathfrak{B}(M) = (\partial M \to \lbar{M}) \in \Spc_\fin^{\Delta^1}.
\]
Moreover, $\mathfrak{B}$ sends an open embedding $N \hookrightarrow M$ to the following morphism in $\Corr((\Spc_\fin^{\Delta^1})^\opp)$ %
\[
  \begin{tikzcd}
    \partial N \ar{d} \ar{r} & \lbar{M} \setminus N \ar{d} & \partial M \ar{d} \ar{l} \\
    \lbar{N} \ar{r} & \lbar{M} \arrow[phantom]{ul}[very near start]{\lrcorner} & \lbar{M}. \ar{l}[swap]{\simeq}
  \end{tikzcd} \teq\label{eq:co_corr_from_emb}
\]
We will often suppress the morphisms and simply use, for example, $(\partial M, \lbar{M})$ to denote an object in $\Spc_\fin^{\Delta^1}$ to make diagrams and formulas more compact.

\subsubsection{The functor $\cMap$}
We have a natural functor
\[
  \cMap: (\Spc^{\Delta^1})^\opp \to \Stk
\]
which assigns to each object $(N\to M) \in (\Spc^{\Delta^1})^\opp$ an object $(BB,BG)^{N,M} \defeq BB^N \times_{BG^N} BG^M$,\footnote{See also \cref{rmk:abuse_of_notation_forget_n_mixed,ftn:Stk_vs_StkFq}.} which is precisely the stack of commutative squares
\[
  \begin{tikzcd}
    N \ar{r} \ar{d} & BB \ar{d} \\
    M \ar{r} & BG.
  \end{tikzcd}
\]
This functor upgrades naturally to an eponymous functor
\[
  \cMap: \Corr((\Spc^{\Delta^1})^\opp) \to \Corr(\Stk),
\]
used in~\cref{eq:M_as_composition} above.

By construction, $\cMap$ turns colimits in $\Spc^{\Delta^1}$ to limits in $\Stk$.

\subsubsection{The functor $M$}
The functor $M$ in \cref{eq:Hecke_as_a_functor} is given by the following result.

\begin{lem} \label{lem:factors_through_prop_sm}
  The composition $\cMap\circ \mathfrak{B}$ factors through $\Corr(\Stk)_{\proper;\sm}$, i.e., the horizontal (resp. vertical) maps are smooth (resp. schematic and proper).
\end{lem}
\begin{proof}
  Let $N \hookrightarrow M$ be a morphism in $\Mnfd_1$, i.e., it is an open embedding of $1$-manifolds. We will now show that the following map
  \[
    (BB, BG)^{\lbar{M}\setminus N, \lbar{M}} \simeq BB^{\lbar{M} \setminus N} \times_{BG^{\lbar{M}\setminus N}} BG^{\lbar{M}} \to BB^{\partial N} \times_{BG^{\partial N}} BG^{\lbar{N}} \simeq (BB, BG)^{\partial N, \lbar{N}}, \teq\label{eq:left_map_smooth}
  \]
  which is induced by the left square of \cref{eq:co_corr_from_emb}, is smooth. Since the left square of \cref{eq:co_corr_from_emb} is a pushout, we have
  \[
    (BB, BG)^{\lbar{M}\setminus N, \lbar{M}} \simeq (BB, BG)^{\partial N, \lbar{N}} \times_{BB^{\partial N}} BB^{\lbar{M} \setminus N}
  \]
  and the map \cref{eq:left_map_smooth} identifies with the vertical map on the left of the following pullback diagram
  \[
    \begin{tikzcd}
      (BB, BG)^{\lbar{M} \setminus N, \lbar{M}} \ar{d} \ar{r} & BB^{\lbar{M} \setminus N} \ar{d} \\
      (BB, BG)^{\partial N, \lbar{N}} \ar{r} & BB^{\partial N}.
    \end{tikzcd}
  \]
  It thus suffices to show that the map
  \[
    BB^{\lbar{M}\setminus N} \to BB^{\partial N} \teq\label{eq:left_map_smooth_reduction}
  \]
  is smooth. Without loss of generality, we can assume that $M$ is connected, in which case, $M$ is either a circle or a line.

  When $N$ is empty, it is clear that \cref{eq:left_map_smooth_reduction} is smooth since it is simply
  \[
    BB^{\lbar{M}} \to \pt,
  \]
  where $BB^{\lbar{M}}$ is either $BB$ or $\frac{B}{B}$ depending on whether $M$ is a line or a circle. Thus, it remains to treat the case where $N$ is non-empty.

  When $M$ is a line, the desired statement follows from the fact that \cref{eq:left_map_smooth_reduction} is a product of (copies of) of the smooth maps $\id_{BB}$ and $\Delta_{BB}: BB \to BB\times BB$. When $M$ is a circle, since the only embedding of a circle into itself is a homeomorphism, $N$ can only be a circle or a disjoint union of lines. In the first case, the map under consideration is an equivalence and hence, it is smooth. In the second case, we also see that \cref{eq:left_map_smooth_reduction} is a product of (copies of) the diagonal map $\Delta_{BB}: BB \to BB\times BB$, which is smooth.

  Next, we will show that the following map
  \[
    (BB, BG)^{\lbar{M}\setminus N, \lbar{M}} \simeq BB^{\lbar{M}\setminus N} \times_{BG^{\lbar{M}\setminus N}} BG^{\lbar{M}} \to BB^{\partial M} \times_{BG^{\partial M}} BG^{\lbar{M}} \simeq (BB, BG)^{\partial M, \lbar{M}}, \teq\label{eq:right_map_proper}
  \]
  induced by the right square of \cref{eq:co_corr_from_emb}, is proper. As above, we can (and we will) assume that $M$ is connected. Note that when $N$ is empty, \cref{eq:right_map_proper} becomes
  \[
    BB^{\lbar{M}} \to BB^{\partial M} \times_{BG^{\partial M}} BG^{\lbar{M}},
  \]
  which is equivalent to
  \[
    \frac{B}{B} \simeq BB^{S^1} \to \frac{G}{G}
  \]
  when $M \simeq S^1$, and to
  \[
    BB \to BB\times_{BG} BB
  \]
  when $M \simeq \mathbb{R}$. Both of these are easily seen to be proper.

  When $N$ is not empty, then $\lbar{M} \setminus N$ and $\partial M$ are homotopically equivalent to a (possibly empty) disjoint union of points. Moreover, $\partial M \to \lbar{M} \setminus N$ is homotopically equivalent to an inclusion of connected components. We will thus treat them as disjoint unions of points in what follows. In particular, we can write $\lbar{M} \setminus N \simeq \partial M \sqcup ((\lbar{M}\setminus N) \setminus \partial M)$. The map \cref{eq:right_map_proper} factors as follows
  \begin{align*}
    BB^{\lbar{M}\setminus N} \times_{BG^{\lbar{M}\setminus N}} BG^{\lbar{M}}
     & \simeq (BB^{\partial M} \times_{BG^{\partial M}} BG^{\lbar{M}}) \times_{BG^{\lbar{M}}} (BB^{(\lbar{M}\setminus N) \setminus \partial M} \times_{BG^{(\lbar{M}\setminus N) \setminus \partial M}} BG^{\lbar{M}}) \\
     & \to BB^{\partial M} \times_{BG^{\partial M}} BG^{\lbar{M}}
  \end{align*}
  where the last map is induced by $BB^{(\lbar{M}\setminus N) \setminus \partial M} \to BG^{(\lbar{M}\setminus N) \setminus \partial M}$ which is proper since it is just a product of maps of the form $BB\to BG$. The proof thus concludes.
\end{proof}

\subsubsection{The functor $\Hecke^{?,\ren}_G$}
Applying the functor $\Shv_{?, !}^{\ren, *}$ of~\cite[\crefnolink{mg:eq:renormalized_mixed_sheaves_corr_!-pushforward,mg:thm:graded_sheaves_correspondence}]{ho_revisiting_2022}, which is right-lax symmetric monoidal, we complete \cref{eq:Hecke_as_a_functor}. Note that in our case, we only need $\Shv_{?, !}^{\ren, *}$ to encode (renormalized) $*$-pullbacks along smooth morphisms and (renormalized) $!$-pushforwards along proper morphisms rather than the more general case described in~\cite{ho_revisiting_2022}.

\subsubsection{Renormalized finite Hecke categories}
Being a right-lax symmetric monoidal functor, $\Hecke^{?,\ren}_G$ sends algebra objects to algebra objects. In particular,
\[
  \Hecke^{?,\ren}_G(\mathbb{R})
  \simeq \Shv_?((BB\times BB) \times_{BG\times BG} BG)^\ren
  \simeq \Shv_?(BB\times_{BG} BB)^\ren
\]
has an algebra structure, induced by the algebra structure on $\mathbb{R} \in \Mnfd_1$, whose multiplication is given by $\mathbb{R}^{\sqcup 2} \hookrightarrow \mathbb{R}$. Chasing through \cref{eq:co_corr_from_emb}, we see that this is precisely the monoidal structure on $\Shv_?(BB \times_{BG} BB)^\ren$ defined earlier via the correspondence \cref{eq:Hecke_convolution_correspondence}. This justifies the abuse of notation where we use $\Hecke^{?,\ren}_G$ to denote both the functor and its value on $\mathbb{R}$.

\subsubsection{Horocycle correspondence}
We will now explain how the horocycle correspondence appears naturally from this point of view. First, observe that as $\partial S^1 \simeq \emptyset$ and hence, $\Hecke^{?,\ren}_G(S^1) \simeq \Shv_?\left( \frac{G}{G} \right)^\ren$ since $BG^{S^1} \simeq \frac{G}{G}$, where the quotient is taken using the conjugation action.

Let $\psi: \mathbb{R} \to S^1$ be an embedding. The functor $\mathfrak{B}$ carries $\psi$ to the following diagram
\[
  \begin{tikzcd}
    \pt \sqcup \pt \ar{d} \ar{r} & \pt \ar{d} & \ar{l} \emptyset \ar{d} \\
    \mathbb{R} \ar{r} & S^1 & \ar{l} S^1
  \end{tikzcd}
\]
which is sent to
\[
  \begin{tikzcd}
    B\backslash G/B \simeq BB \times_{BG} BB & \ar{l} \frac{G}{B} \ar{r} & \frac{G}{G}
  \end{tikzcd}
\]
under $\cMap$, where we have used $BB \times_{BG} \frac{G}{G} \simeq \frac{G}{B}$. This is the horocycle correspondence appearing in \cref{thm:intro:trace_center_vs_character_sheaves}.

\subsection{Augmented cyclic bar construction via \texorpdfstring{$\Mnfd_1$}{Mnfd₁}}
\label{subsec:augmented_cyclic_bar_Mnfd1}
Using the geometric picture of \cref{subsec:Hecke_mixed_ren_as_functor}, we will now produce an augmented simplicial object that will be used as the input for \cref{prop:Beck-Chevalley_fully_faithful}. The underlying simplicial object is the same as the one that computes the categorical trace of $\Hecke^{?,\ren}_G$ introduced in \cref{eq:cyclic_bar_H_mixed_ren}.

\subsubsection{Circle geometry and augmented simplicial sets}
\label{subsubsec:cicle_geometry_vs_aug_simpl}
Consider the over-category $(\Mnfd_1)_{/S^1}$. Fix a morphism $\psi: \mathbb{R} \to S^1$, and consider $(\Mnfd_1)_{\psi//S^1} \defeq ((\Mnfd_1)_{/S^1})_{\psi/}$. We note that $(\Mnfd_1)_{/S^1}$ has a final object, by construction, and moreover, it is precisely the category obtained from $(\Disk_1)_{/S^1} \defeq \Disk_1 \times_{\Mnfd_1} (\Mnfd_1)_{/S^1}$ by adjoining a final object. Similarly, the category $(\Mnfd_1)_{\psi//S^1}$ also has a final object and is equivalent to the category obtained from $(\Disk_1)_{\psi//S^1}$ by adjoining a final object.

We will now relate $(\Mnfd_1)_{\psi // S^1}$ and $\Delta^{\opp}_+$ using a variation of the construction in~\cite[\S5.5.3]{lurie_higher_2017}. Note that an object of $(\Mnfd_1)_{\psi//S^1}$ is given by a diagram
\[
  \begin{tikzcd}[column sep=small]
    \mathbb{R} \ar{dr}[swap]{\psi} \ar{rr}{j} && U \ar{dl}{\psi'} \\
    & S^1
  \end{tikzcd}
\]
which commutes up to isotopy, where $U$ is either $S^1$ or a finite disjoint union of copies of $\mathbb{R}$. The set of components $\pi_0(S^1 \setminus \psi'(U))$ is finite: empty when $U$ is $S^1$ and equal to the number of components of $U$ when it is a disjoint union of copies of $\mathbb{R}$. Fix an orientation of the circle. We define a linear ordering $\leq$ $\pi_0(S^1 \setminus \psi'(U))$: if $x, y \in S^1$ belong to different components of $S^1 \setminus \psi'(U)$, then we write $x<y$ if the three points $(x, y,\psi'(j(0)))$ are arranged in clockwise order around the circle, and $y<x$ otherwise. This construction determines a functor from $(\Mnfd_1)_{\psi//S^1}$ to the opposite of the category of finite linearly ordered sets, which is $\Delta_+^\opp$.

\begin{lem} \label{lem:cyclic_bar_vs_disk_in_S^1}
  The above construction determines an equivalence of $\infty$-categories
  \[
    \theta_+: (\Mnfd_1)_{\psi//S^1} \to \Delta_+^\opp,
  \]
  which restricts to an equivalence of $\infty$-categories
  \[
    \theta \defeq \theta_+|_{(\Disk_1)_{\psi//S^1}}: (\Disk_1)_{\psi//S^1} \to \Delta^\opp.
  \]
\end{lem}
\begin{proof}
  The second equivalence is~\cite[Lemma 5.5.3.10]{lurie_higher_2017}. The first part is a direct extension of the second by adjoining a final object.
\end{proof}

\subsubsection{The augmented simplicial category}
We are now ready to construct the augmented simplicial category. Let $\lbar{\Tr}(\Hecke^{?,\ren}_G)_\bullet$ be the augmented simplicial category given by the following composition of functors
\[
  \begin{tikzcd}
    \Delta_+^\opp \ar{r}{\theta_+^{-1}} & (\Mnfd_1)_{\psi//S^1} \ar{r}{\varphi} & \Mnfd_1 \ar{r}{\Hecke^{?,\ren}_G} & \Shv_?(\pt)\hphMod.
  \end{tikzcd} \teq\label{eq:aug_simplicial_Tr_via_Mnfd_1}
\]
Here, $\varphi$ is the obvious forgetful functor.

Since $BB \times_{BG} BB \simeq B\backslash G/B$ is a finite orbit stack by the Bruhat decomposition, \cref{prop:categorical_Kuenneth_finite_orbit} implies that $\lbar{\Tr}(\Hecke^{?,\ren}_G)_\bullet|_{\Delta^\opp} \simeq \Tr(\Hecke^{?,\ren}_G)_\bullet$ of \cref{eq:cyclic_bar_H_mixed_ren}; see also~\cite[Remarks 5.5.3.13 and 5.5.3.14]{lurie_higher_2017}. In particular,
\[
  \lbar{\Tr}(\Hecke^{?,\ren}_G)_n \simeq \Tr(\Hecke^{?,\ren}_G)_n \simeq \Shv_?((BB\times_{BG} BB)^{n+1})^\ren, \qquad n\geq 0
\]
and
\[
  \lbar{\Tr}(\Hecke^{?,\ren}_G)_{-1} \simeq \Shv_?\left( \textstyle\frac{G}{G} \right)^\ren
\]
where, by convention, $\frac{G}{G}$ denotes the quotient of $G$ by itself under the adjoint action.

\subsubsection{``Small'' variants}  \label{subsubsec:small_variants_mixed}
We have the following ``small'' variants $\Hecke^?_G$ and $\lbar{\Tr}(\Hecke^?_G)_\bullet$ of $\Hecke^{?,\ren}_G$ and $\lbar{\Tr}(\Hecke^{?,\ren}_G)_\bullet$, given by the following diagram
\[
  \begin{tikzcd}
    \Delta_+^\opp \ar[bend right=9]{rrrr}[swap]{\lbar{\Tr}(\Hecke^?_G)_\bullet} \ar{r}{\theta_+^{-1}} & (\Mnfd_1)_{\psi//S^1} \ar{r}{\varphi} & \Mnfd_1 \ar{r}{M} \ar[bend left=19]{rr}{\Hecke^?_G}  & \Corr(\Stk)_{\proper;\sm} \ar{r}{\Shv_{?, c, !}^*} & \Shv_{?,c}(\pt)\hphMod,
  \end{tikzcd}
\]
The only difference between the two versions is that in the last step, we use $\Shv_{?,c,!}^*$ rather than $\Shv_{?, !}^{\ren, *}$.

\subsection{Adjointability}
\label{subsec:adjointability}
We are now ready to show that $\lbar{\Tr}(\Hecke^{?,\ren}_G)_\bullet$ satisfies the adjointability condition required by~\cref{prop:Beck-Chevalley_fully_faithful}.  The result follows from a simple geometric statement about topological $1$-manifolds.

\subsubsection{Adjointable squares in a category of correspondences}
For any category $\mathcal{C}$, a commutative square in $\Corr(\mathcal{C})$, which illustrates two morphisms from $x$ to $y'$ to be equivalent, has the following shape
\[
  \begin{tikzcd}
    x \arrow[phantom]{dr}{\scriptstyle{2}} & \ar{l} c_{xy} \arrow[phantom]{dr}{\scriptstyle{1}} \ar{r} & y \\
    c_{xx'} \arrow[phantom]{dr}{\scriptstyle{3}} \ar{u} \ar{d} & \ar{l} c_{xy'} \arrow[phantom]{dr}{\scriptstyle{4}} \ar{u} \ar{d} \ar{r} & c_{yy'} \ar{u} \ar{d} \\
    x' & \ar{l} c_{x'y'} \ar{r} & y'.
  \end{tikzcd} \teq\label{eq:commutative_diagram_corr}
\]
where $1$ and $3$ are pullback squares.

\begin{defn}
  The commutative diagram in $\Corr(\mathcal{C})$ \cref{eq:commutative_diagram_corr} is called \emph{adjointable} if $2$ and $4$ are also pullback squares, i.e., all squares are pullback squares.
\end{defn}

\subsubsection{}
Let $\alpha: [n] \to [m]$ be a morphism in $\Delta_+$ and $M \to N$ the corresponding morphism in $\Mnfd_1$ via $\varphi \circ \theta_+^{-1}$, see \cref{eq:aug_simplicial_Tr_via_Mnfd_1}. Let $M_+ \to N_+$ be the morphism associated to $\alpha_+: [n+1] \to [m+1]$. $M_+$ can be obtained from $M$ by deleting the image of $[-\varepsilon, \varepsilon] \subset \mathbb{R}$ in $M$ for some fixed $\varepsilon$, and similarly for $N$. This construction is functorial and hence, we obtain the map $M_+ \to N_+$.

We have the following commutative diagram in $\Mnfd_1$
\[
  \begin{tikzcd}
    M_+ \ar{r} \ar{d} & M \ar{d} \\
    N_+ \ar{r} & N
  \end{tikzcd}
\]
which induces, via the $\mathfrak{B}$ construction of \cref{subsubsec:B_functor}, the following commutative diagram in $\Corr((\Spc_\fin^{\Delta^1})^\opp)$
\[
  \begin{tikzcd}
    (\partial M_+, \lbar{M}_+) \ar{d} \ar{r} & (\lbar{M} \setminus M_+, \lbar{M}) \ar{d} & \ar{l} (\partial M, \lbar{M}) \ar{d} \\
    (\lbar{N}_+ \setminus M_+, \lbar{N}_+) \ar{r} & (\lbar{N} \setminus M_+, \lbar{N}) & \ar{l} (\lbar{N} \setminus M, \lbar{N}) \\
    (\partial N_+, \lbar{N}_+) \ar{r} \ar{u} & (\lbar{N} \setminus N_+, \lbar{N}) \ar{u} & \ar{l} (\partial N, \lbar{N}). \ar{u}
  \end{tikzcd} \teq\label{eq:adjointable_coCorr_Spc}
\]
Applying the $\cMap$ construction, we obtain the following commutative diagram in $\Corr(\Stk_{\Fq})_{\proper;\sm}$
\[
  \begin{tikzcd}
    (BB, BG)^{\partial M_+, \lbar{M}_+} \arrow[phantom]{dr}{\scriptstyle{2}} & \ar{l}[swap]{f} (BB, BG)^{\lbar{M} \setminus M_+, \lbar{M}} \arrow[phantom]{dr}{\scriptstyle{1}} \ar{r}{g} & (BB, BG)^{\partial M, \lbar{M}} \\
    (BB, BG)^{\lbar{N}_+ \setminus M_+, \lbar{N}_+} \arrow[phantom]{dr}{\scriptstyle{3}} \ar{u}{p} \ar{d}[swap]{q} & \ar{l}[swap]{f} (BB, BG)^{\lbar{N} \setminus M_+, \lbar{N}} \arrow[phantom]{dr}{\scriptstyle{4}} \ar{u}{p} \ar{d}[swap]{q} \ar{r}{g} & (BB, BG)^{\lbar{N} \setminus M, \lbar{N}} \ar{u}{p} \ar{d}[swap]{q} \\
    (BB, BG)^{\partial N_+, \lbar{N}_+} & \ar{l}[swap]{f} (BB, BG)^{\lbar{N} \setminus N_+, \lbar{N}} \ar{r}{g} & (BB, BG)^{\partial N, \lbar{N}}.
  \end{tikzcd} \teq\label{eq:adjointable_Corr_Stk}
\]

\begin{lem}
  The diagram \cref{eq:adjointable_coCorr_Spc} is adjointable in $\Corr((\Spc_\fin^{\Delta^1})^\opp)$. As a result, \cref{eq:adjointable_Corr_Stk} is adjointable in $\Corr(\Stk)_{\proper;\sm}$.
\end{lem}
\begin{proof}
  The second part follows from the first part because $\cMap$ turns colimits to limits, and hence, in particular, it sends pushout squares to pullback squares. The first part is an explicit and elementary statement about gluing 1-manifolds that is easier to check directly than to describe. We leave the details to the reader.
\end{proof}

\begin{prop} \label{prop:trace_Hecke_aug_adjointable}
  The augmented simplicial object $\lbar{\Tr}(\Hecke^{?,\ren}_G)_\bullet$ in $\Shv_?(\pt)\hphMod$ satisfies the condition of \cref{prop:Beck-Chevalley_fully_faithful}. Namely, for every $\alpha: [n] \to [m]$ in $\Delta_+$, the diagram
  \[
    \begin{tikzcd}
      \lbar{\Tr}(\Hecke^{?,\ren}_G)_{m+1} \ar{d} \ar{r}{d_0} & \lbar{\Tr}(\Hecke^{?,\ren}_G)_m \ar{d} \\
      \lbar{\Tr}(\Hecke^{?,\ren}_G)_{n+1} \ar{r}{d_0} & \lbar{\Tr}(\Hecke^{?,\ren}_G)_n
    \end{tikzcd} \teq\label{eq:adjointable_Tr_aug_simplicial}
  \]
  is vertically right adjointable.
\end{prop}
\begin{proof}
  Note that \cref{eq:adjointable_Tr_aug_simplicial} is obtained from \cref{eq:adjointable_Corr_Stk} by applying $\Shv_{?,!}^{*,\ren}$. By \cref{lem:factors_through_prop_sm}, we know that in \cref{eq:adjointable_Corr_Stk}, the morphisms $f$ and $p$ (resp. $g$ and $q$) are smooth (resp. proper). To prove that \cref{eq:adjointable_Tr_aug_simplicial} is vertically right adjointable, it suffices to show that we have the following natural equivalence
  \[
    g_! f^* p_* q^! \xrightarrow{\simeq} p_* q^! g_! f^*,
  \]
  where we start from the bottom left of \cref{eq:adjointable_Corr_Stk}. Indeed, we have
  \begin{align*}
    g_! f^* p_* q^! \simeq g_! p_* f^* q^! \simeq p_* g_! q^! f^* \simeq p_* q^! g_! f^*.
  \end{align*}
  where the first, second, and third equivalences are due to the following reasons, respectively
  \begin{itemize}
    \item smooth base change for square 2, using the fact that $f$ is smooth;
    \item commuting upper $!$ and upper $*$ using the fact that $f$ is smooth, and commuting lower $!$ and lower $*$ using the fact that $g$ is proper; and
    \item using the fact that $g$ is proper, the desired equivalence follows from $g_* q^! \simeq q^! g_*$, which is the Verdier dual of the usual proper base change result.
  \end{itemize}
\end{proof}

\begin{rmk}
  Unlike the ``big'' version, the ``small'' variant $\lbar{\Tr}(\Hecke^?_G)_\bullet$ discussed in \cref{subsubsec:small_variants_mixed} does not satisfy the adjointability condition of \cref{prop:Beck-Chevalley_fully_faithful}. This is because the right adjoint to the unit map does not preserve constructibility. Indeed, the right adjoint is given by $p_*q^!$ in the following diagram
  \[
    \begin{tikzcd}[column sep=small]
      & \ar{dl}[swap]{q} BB \ar{dr}{p} \\
      BB\times_{BG} BB & & \pt.
    \end{tikzcd}
  \]
  But now, note that $p_*$ does not preserve constructibility.
\end{rmk}

\subsection{Traces and Drinfel'd centers of finite Hecke categories}
\label{subsec:trace_center_Hecke}
We are now ready to complete the proof of the first main result. The trace case follows directly from the discussion above and thus will be handled at the beginning of this subsection. We will then deduce the Drinfel'd center case from the trace case.

\subsubsection{Traces of finite Hecke categories}
We will now prove the trace part of the first main result.

\begin{thm} \label{thm:trace_Hecke}
  The trace of $\Hecke^{?,\ren}_G$ (resp. $\Hecke^?_G$) where $\Hecke^{?,\ren}_G$ (resp. $\Hecke^?_G$) is viewed as an algebra object of $\Shv_?(\pt)\hphMod$ (resp. $\Shv_{?, c}(\pt)\hphMod$) coincides with the full subcategory of $\Shv_?\left( \frac{G}{G} \right)^\ren$ (resp. $\Shv_{?, c}\left( \frac{G}{G} \right)$) generated under colimits (resp. finite colimits and idempotent completion) by the image of $\Hecke^{?,\ren}_G$ (resp. $\Hecke^?_G$) via $q_!p^*$ in the horocycle correspondence\footnote{See \cref{rmk:abuse_of_notation_forget_n_mixed}.}
  \[
    \begin{tikzcd}[sep=tiny]
      & \frac{G}{B} \ar{dl}[swap]{p} \ar{dr}{q} \\
      B\backslash G/B && \frac{G}{G}.
    \end{tikzcd} \teq\label{eq:horocycle_correspondence}
  \]
  Moreover, under this identification, the natural functor $\tr: \Hecke^{?,\ren}_G \to \Tr(\Hecke^{?,\ren}_G)$ (resp. $\tr: \Hecke^?_G \to \Tr(\Hecke^?_G)$) is identified with $q_!p^*$.
\end{thm}
\begin{proof}
  We start with the ``big'' variant $\Hecke^{?,\ren}_G$. \cref{prop:trace_Hecke_aug_adjointable,prop:Beck-Chevalley_fully_faithful} imply that the natural functor $\Tr(\Hecke^{?,\ren}_G) \to \Shv_?\left( \frac{G}{G} \right)^\ren$ is fully faithful. Moreover, the functor $\Hecke^{?,\ren}_G \to \Tr(\Hecke^{?,\ren}_G) \hookrightarrow \Shv_?\left( \frac{G}{G} \right)^\ren$ identifies with $q_! p^*$ in the horocycle correspondence.

  By~\cite[Volume I, Chapter 1, Proposition 8.7.4]{gaitsgory_study_2017}, we know that $\Tr(\Hecke^{?,\ren}_G)$ is compactly generated by the essential image of $\Hecke^?_G \to \Hecke^{?,\ren}_G \to \Tr(\Hecke^{?,\ren}_G)$. Note that the conditions required to apply this result amount to the existence of the ``small'' variant introduced in \cref{subsubsec:small_variants_mixed}.

  Combining the two statements, we conclude that $\Tr(\Hecke^{?,\ren}_G)$ is identified with the full subcategory of $\Shv_?\left( \frac{G}{G} \right)^\ren$ compactly generated by the image of $\Hecke^?_G = \Shv_{?, c}(B\backslash G/B)$ under the horocycle correspondence. In particular, it is generated by the image of $\Hecke^{?, \ren}_G$ under colimits.

  The ``small'' variant is obtained from the first by taking the full subcategories of compact objects, using \cref{prop:big_vs_small}. The statement regarding generation under finite colimits and idempotent splittings follows from~\cite[Volume I, Chapter 1, Lemma 7.2.4.(1'')]{gaitsgory_study_2017}.
\end{proof}

\subsubsection{Drinfel'd centers of the finite Hecke categories}
The case of Drinfel'd center is slightly more subtle. Consider the following versions $\Hecke^{?,\ren,!}_G$ (resp. $\Hecke^{?,!}_G$) of the Hecke categories where instead of using $\Shv_{?,!}^{\ren,*}$ (resp. $\Shv_{?, c, !}^*$), we use $\Shv_{?,*}^{\ren,!}$ (resp. $\Shv_{?,c,*}^!$). More concretely, we use $g_*f^!$ in the correspondence \cref{eq:Hecke_convolution_correspondence} to define the convolution monoidal structure.

Since $f$ is smooth of relative dimension $\dim B$, $f^! \simeq f^*[2\dim B](\dim B)$ in the mixed case and $f^! \simeq f^*[2\dim B]\lrangle{2\dim B}$ in the graded case, where $(-)$ denotes the Tate twist and $\lrangle{-}$ denote the grading shift defined in~\cite[\crefnolink{mg:subsubsec:shift_of_gradings}]{ho_revisiting_2022}.\footnote{The factor $2$ is there because weight is twice the Tate twist.} Note that in the ungraded setting, we simply have $f^! \simeq f^*[2\dim B]$. Thus, by cohomologically shifting and Tate twisting $[-2\dim B](-\dim B)$ for the mixed case (resp. cohomologically shifting and grade shifting $[-2\dim B]\lrangle{-2\dim B}$ for the graded case, resp. cohomologically shifting $[-2\dim B]$ for the ungraded case), we obtain an equivalence of monoidal categories
\begin{align*}
  \Hecke^{?,\ren}_G \xrightarrow{\simeq} \Hecke^{?,\ren, !}_G \qquad & \text{and}\qquad \Hecke^{?}_G \xrightarrow{\simeq} \Hecke^{?,!}_G
\end{align*}

Even though the $!$-version is equivalent to the usual version, it is, as we shall see, technically more advantageous to use when studying the center.

\subsubsection{} \label{subsubsec:self_duality}
Observe that the category $\Hecke^{?,\ren,!}_G$ is self-dual as an object in $\Shv_?(\pt)\hphMod$. Indeed, first note that any object $\mathcal{Y}\in\Corr(\Stk)$ is self-dual, with duality datum given by
\[
  \pt \leftarrow \mathcal{Y} \rightarrow \mathcal{Y} \times \mathcal{Y} \qquad\text{and}\qquad \mathcal{Y} \times \mathcal{Y} \leftarrow \mathcal{Y} \rightarrow \pt.
\]
Thus, if $\mathcal{Y}$ is a finite orbit stack, such as $\mathcal{Y} = BB \times_{BG} BB$, $\Shv_{?,*}^{\ren,!}$ turns the self-duality datum above into a self-duality datum of $\Shv_?(X)^\ren$, using \cref{prop:categorical_Kuenneth_finite_orbit}. In particular, $\Hecke^{?, \ren, !}_G$ is self dualizable.

Suppose
\[
  \mathcal{X} \leftarrow \mathcal{Z} \rightarrow \mathcal{Y}
\]
is a morphism from $\mathcal{X}$ to $\mathcal{Y}$ in $\Corr(\Stk)$. The dual of this morphism is precisely
\[
  \mathcal{Y} \leftarrow \mathcal{Z} \rightarrow \mathcal{X}.
\]
Thus, if $\mathcal{X}$ and $\mathcal{Y}$ are finite orbit stacks, then the functors $\Shv_?(\mathcal{X})^\ren \to \Shv_?(\mathcal{Y})^\ren$ and $\Shv_?(\mathcal{Y})^\ren \to \Shv_?(\mathcal{X})^\ren$ associated to the two correspondences above are dual to each other. Note that when following the correspondences here, we use $(-)_*(-)^!$.

\subsubsection{}
We are now ready to compute the Drinfel'd centers of $\Hecke^{?,\ren}_G$.

\begin{thm}
  \label{thm:center_big_Hecke}
  The Drinfel'd center $\Z(\Hecke^{?,\ren}_G)$ of $\Hecke^{?,\ren}_G$, where the latter is viewed as an object in the category $\Alg(\Shv_?(\pt)\hphMod)$, coincides with its trace $\Tr(\Hecke^{?,\ren}_G)$. Moreover, under the identification of $\Tr(\Hecke_G^{?,\ren})$ as a full subcategory of $\Shv_?\left( \frac{G}{G} \right)^\ren$, the natural functor $\z: \Z(\Hecke^{?,\ren}_G) \to \Hecke^{?,\ren}_G$ can be identified with $p_*q^!$ in the horocycle correspondence \cref{eq:horocycle_correspondence}. In particular, we have a pair of adjoint functors $\tr \dashv \z$.
\end{thm}
\begin{proof}
  To simplify the notation, we will write $\mathcal{A} = \Hecke^{?,\ren,!}_G$ and $\mathcal{B} = \Hecke^{?,\ren}_G$. Moreover, unless otherwise specified, all tensors and $\cuHom$ are over $\Shv_?(\pt)$. We will therefore omit $\Shv_?(\pt)$ from the notation.

  Recalling from \cref{defn:trace_center}, we have
  \[
    \Z(\mathcal{A}) \simeq \cuHom_{\mathcal{A} \otimes \mathcal{A}^\rev}(\mathcal{A},\mathcal{A}),
  \]
  i.e., the category of continuous $\mathcal{A} \otimes \mathcal{A}^\rev$-linear functors from $\mathcal{A}$ to itself. As in~\cite[\S5.1.1]{ben-zvi_integral_2010}, we have an equivalence of categories $\mathcal{A} \simeq |\mathcal{A}^{\otimes (\bullet + 2)}|$. Thus,
  \begin{align*}
     & \alignsep\cuHom_{\mathcal{A} \otimes \mathcal{A}^\rev}(\mathcal{A},\mathcal{A})                             \\
     & \simeq \cuHom_{\mathcal{A} \otimes \mathcal{A}^\rev} (|\mathcal{A}^{\otimes (\bullet + 2)}|, \mathcal{A})   \\
     & \simeq \Tot(\cuHom_{\mathcal{A} \otimes \mathcal{A}^\rev}(\mathcal{A}^{\otimes (\bullet + 2)},\mathcal{A})) \\
     & \simeq \Tot(\cuHom(\mathcal{A}^{\otimes \bullet}, \mathcal{A}))                                             \\
     & \simeq \Tot(\mathcal{A}^{\otimes (\bullet + 1)}) \teq\label{eq:turning_to_dual}                             \\
     & \simeq |\mathcal{B}^{\otimes (\bullet + 1)}| \teq\label{eq:passing_to_left_adjoint}                         \\
     & \simeq \Tr(\mathcal{B}).
  \end{align*}
  Here, \cref{eq:turning_to_dual} is obtained by passing to the duals, using the discussion in \cref{subsubsec:self_duality}. Moreover, \cref{eq:passing_to_left_adjoint} is obtained by passing to left adjoints, see~\cite[Corollary 5.5.3.4]{lurie_higher_2017-1} (see also~\cite[Volume I, Chapter 1, Proposition 2.5.7]{gaitsgory_study_2017}), and noticing that the resulting diagram is precisely the diagram used to compute the trace of $\mathcal{B} \defeq \Hecke^{?,\ren}_G$. We deduce the corresponding statement for $\Z(\mathcal{B})$ using the equivalence of monoidal categories $\mathcal{A} \simeq \mathcal{B}$.

  The second part regarding the identification of the functor $\z$ follows from the argument above.
\end{proof}

\subsubsection{}
We will now study the centers of the ``small'' version $\Hecke^?_G$, which is more subtle than the case of traces above.

\begin{thm}
  \label{thm:center_small_Hecke}
  The Drinfel'd center $\Z(\Hecke^?_G)$ of $\Hecke^?_G\in\Alg(\Shv_{?,c}(\pt)\hphMod)$ is equivalent to the full subcategory of $\Z(\Hecke^{?,\ren}_G)$ consisting of objects whose images under the natural central functor $\z$ are constructible as sheaves on $B\backslash G/B$.
\end{thm}
\begin{proof}
  For brevity's sake, we write $\mathcal{A} = \Hecke^{?,\ren,!}_G$ and $\mathcal{A}_0 = \Hecke^{?,!}_G$. Moreover, as all tensors and $\cuHom$, unless otherwise specified, will be over $\Shv_?(\pt)$ or $\Shv_{?,c}(\pt)$ depending on whether we are working with ``big'' or ``small'' categories, which should be clear from the context. We will thus not include $\Shv_?(\pt)$ or $\Shv_{?,c}(\pt)$ in the notation.

  By definition,
  \[
    \Z(\mathcal{A}_0) \simeq \cuHom_{\mathcal{A}_0 \otimes \mathcal{A}_0^\rev}(\mathcal{A}_0, \mathcal{A}_0),
  \]
  the category of exact $\mathcal{A}_0 \otimes \mathcal{A}_0^\rev$-linear functors from $\mathcal{A}_0$ to itself. As in the proof of \cref{thm:center_big_Hecke}, we have
  \begin{align*}
    \Z(\mathcal{A}_0) \simeq \Tot(\cuHom(\mathcal{A}_0^{\otimes \bullet}, \mathcal{A}_0)),
  \end{align*}
  which naturally embeds into
  \[
    \Z(\mathcal{A}) \simeq \Tot(\cuHom(\mathcal{A}^{\otimes \bullet}, \mathcal{A}))
  \]
  whose essential image is the full subcategory of $\Z(\mathcal{A})$ consisting of objects whose images in $\cuHom(\mathcal{A}^{\otimes \bullet}, \mathcal{A})$ are compact preserving functors. Observe that all the functors between $\mathcal{A}^{\otimes \bullet}$ (which induces functors between $\Tot(\cuHom(\mathcal{A}^{\otimes \bullet}, \mathcal{A}))$) are compact preserving since we only push along schematic (in fact proper) morphisms, by \cref{lem:factors_through_prop_sm}. Thus, $\Z(\mathcal{A}_0)$ is, equivalently, the full subcategory of $\Z(\mathcal{A})$ spanned by objects whose images in
  \[
    \cuHom(\mathcal{A}^{\otimes 0}, \mathcal{A}) \simeq \cuHom(\Shv_?(\pt), \mathcal{A}) \simeq \mathcal{A} \simeq \Shv_\gr(B\backslash G/B)^\ren
  \]
  are compact, i.e., constructible. But this functor is precisely the central functor $\z$ and is identified with $p_*q^!$ in the horocycle correspondence~\cref{eq:horocycle_correspondence} as in the proof of \cref{thm:center_big_Hecke} above.

  The statement for $\Hecke^?_G$ follows since $\Hecke^?_G \simeq \Hecke^{?,!}_G$ and the proof concludes.
\end{proof}

\subsection{Character sheaves}
\label{subsec:char_sheaves}
We will now state our results in terms of character sheaves. We start with the following definition.

\begin{defn}
  The \emph{renormalized category of mixed (resp. graded, resp. ungraded) unipotent character sheaves} of $G$, denoted by $\Ch^{\unip,?,\ren}_G$ where $?=\mixed$ (resp. $?=\gr$, resp. $?$ is empty), is the full subcategory of $\Shv_?\left( \frac{G}{G} \right)^\ren$ generated by colimits by the image of $\Hecke^{?,\ren}_G$ under $q_!p^*$ in the horocycle correspondence \cref{eq:horocycle_correspondence}. Equivalently, it is compactly generated by the image of $\Hecke^?_G$ under this functor.

  The \emph{category of mixed (resp. graded, resp. ungraded) unipotent character sheaves} of $G$ is the full subcategory spanned by compact objects, i.e.,
  \[
    \Ch^{\unip,?}_G \defeq (\Ch^{\unip,?,\ren}_G)^c = \Ch^{\unip,?,\ren}_G \cap \Shv_{?, c}\left( \frac{G}{G} \right).
  \]

  Finally, the \emph{category of mixed (resp. graded, resp. ungraded) monodromic character sheaves}, denoted by $\ChT^{\unip,?}_G$, is the full subcategory of $\Ch^{\unip,?,\ren}_G$ consisting of objects who images under $p_* q^!$ are constructible, where $p$ and $q$ are defined in the horocycle correspondence \cref{eq:horocycle_correspondence}.
\end{defn}

\begin{rmk}
  The last equality in the definition of $\Ch^{\unip,?}_G$ uses the following fact: if we have a fully faithful and compact preserving functor $A \hookrightarrow B$, then $A^c = A \cap B^c$. Indeed, we have $A^c \subseteq A \cap B^c$ due to the compact preservation assumption. On the other hand, $A \cap B^c \subset A^c$ is always true by the very definition of compactness. In the case of interest, the fact that the functor involved is compact preserving is proved in \cref{thm:trace_Hecke} above.
\end{rmk}

With the new notation, \cref{thm:trace_Hecke,thm:center_big_Hecke,thm:center_small_Hecke} above can be summarized in the following theorem.

\begin{thm}
  \label{thm:trace_and_center_summary}
  The trace $\Tr(\Hecke^{?,\ren}_G)$ and center $\Z(\Hecke^{?,\ren}_G)$ of $\Hecke^{?,\ren}_G$ coincide with the full subcategory $\Ch^{\unip,?,\ren}_G$ of $\Shv_?\left(\frac{G}{G}\right)^\ren$ generated under colimits by the essential image of $\Hecke^{?,\ren}_G$ under $q_!p^*$ in the correspondence
  \[
    \begin{tikzcd}[sep=tiny]
      & \frac{G}{B} \ar{dl}[swap]{p} \ar{dr}{q} \\
      B\backslash G/B && \frac{G}{G}.
    \end{tikzcd}
  \]
  Moreover, under this identification, the canonical trace and center maps are adjoint $\tr \dashv \z$ and are identified with the adjoint pair $q_! p^* \dashv p_*q^!$.

  The trace $\Tr(\Hecke^?_G)$ coincides with the full subcategory $\Ch^{\unip,?}_G$ of $\Ch^{\unip,?,\ren}_G$ spanned by compact objects. Moreover, the center $\Z(\Hecke^?_G)$ is the full subcategory $\ChT^{\unip,?}_G$ of $\Ch^{\unip,?,\ren}_G$ spanned by the \emph{pre-image} of $\Hecke^?_G$ under the central functor $\z$.
\end{thm}

\begin{rmk} \label{rmk:small_center_is_bigger}
  It is easy to see that $\Ch^{\unip,\gr}_G$ is contained in $\ChT^{\unip,\gr}_G$. In fact, the latter is strictly larger, already when $G=\Gm$, the multiplicative group. Indeed, in this case, $\Hecke^{\gr,\ren}_{\Gm} \simeq \Shv_\gr(B\Gm)^\ren$ and $\Ch^{\unip,\gr,\ren}_{\Gm}$ is the full subcategory of $\Shv_\gr\left( \frac{\Gm}{\Gm} \right)^\ren \simeq \Shv_\gr \left( \Gm \times B\Gm \right)^\ren$ generated by the constant sheaf. Moreover, the natural central functor is simply given by $p_*$ where $p: \Gm \times B\Gm \to B\Gm$ is the projection onto the second factor. Thus, to see that $\ChT^{\unip,\gr}_G$ is larger than $\Ch^{\unip,\gr}_G$, it suffices to realize that $\pi_*: \LS^\unip_\gr(\Gm) \to \Vect^\gr$ does not reflect constructibility, where $\LS^\unip_\gr(\Gm)$ is the full subcategory of $\Shv_\gr(\Gm)$ consisting of unipotent local systems, i.e., it is generated by (grading shifts of) the constant sheaves, and where $\pi: \Gm \to \pt$ is the structure morphism.

  By construction,\footnote{The superscript $\gr$ in $\cHom^\gr$ and $\cEnd^\gr$ denotes the $\Vect^\gr$-enriched $\Hom$-spaces. See~\cite[\crefnolink{mg:subsubsec:enriched_Hom}]{ho_revisiting_2022}.}
  \[
    \LS^\unip_\gr(\Gm) \xrightarrow[\simeq]{\cHom^\gr(\Qlbar, -) \simeq \pi_*} \cEnd^\gr(\Qlbar)^{\opp}\hphMod(\Vect^\gr) \simeq \Qlbar[\alpha]\hphMod(\Vect^\gr),
  \]
  where $\alpha^2=0$ and $\alpha$ lives in cohomological degree $1$ and graded degree $2$. Under this equivalence, $\pi_*$ corresponds to the forgetful functor $\Qlbar[\alpha]\hphMod(\Vect^\gr) \to \Vect^\gr \simeq \Shv_\gr(\pt)$. But now, on the one hand, the object $\Qlbar \in \Qlbar[\alpha]\hphMod(\Vect^\gr)$ is not compact and hence, it corresponds to a non-constructible sheaf on the left. On the other hand, its image in $\Vect^\gr$ is compact.
\end{rmk}

\begin{rmk}
  \label{rmk:trace_center_adj_one_diagram}
  The restriction of the adjunction $\tr:\Hecke^{?,\ren}_G \rightleftarrows \Ch^{\unip,?,\ren}_G:\z$ to small categories yields the following commutative diagram
  \[
    \begin{tikzcd}[column sep=huge]
      & \Ch^{\unip,?}_G \ar[hookrightarrow]{d}{\iota} \ar[shift left=\arrdisp]{dl}{\lbar{\z}} \\
      \Hecke^?_G \ar[shift left=\arrdisp]{ur}{\tr} \ar[shift left=\arrdisp]{r}{\wtilde{\tr}} & \ar[shift left=\arrdisp]{l}{\z} \ChT^{\unip,?}_G
    \end{tikzcd}
  \]
  where $\wtilde{\tr} \simeq \iota\circ \tr$ and $\lbar{\z} \simeq \z \circ \iota$.
\end{rmk}

\subsection{Weight structure and perverse $t$-structure on \texorpdfstring{$\Ch^{\unip,\gr}_G$}{Ch^{u,gr}_G}}
\label{subsec:wt_and_t_on_Ch}
In this subsection, we will show that $\Ch^{\unip,\gr}_G$ inherits the weight structure and perverse $t$-structure from the ambient category $\Shv_{\gr,c}(G/G)$.

\subsubsection{The perverse $t$-structure on $\Ch^{\unip,?}_G$}
\label{subsubsec:perverse_t_structure_Ch_G}
We will now show that for $?\in \{\gr,\emptyset\}$, the category $\Ch^{\unip,?}_G$ inherits the perverse $t$-structure from $\Shv_{?,c}(G/G)$. Moreover, $\Ch^{\unip,\gr}_G$ inherits the weight structure from $\Shv_{\gr,c}(G/G)$.

We start with a general lemma concerning $t$-structures.

\begin{lem}
  \label{lem:t-structure_fullsubcat_Artinian}
  Let $\mathcal{D}$ be a triangulated category equipped with a bounded $t$-structure whose heart is Artinian. Let $\mathcal{C} \defeq \lrangle{\mathcal{L}_i}_{i\in I}$ be the smallest full $\DG$-subcategory of $\mathcal{D}$ containing a collection of simple objects $\{\mathcal{L}_i\}_{i\in I}$ for some indexing set $I$. Then $\mathcal{C}$ is stable under the truncation of $\mathcal{D}$ and hence, it inherits the $t$-structure on $\mathcal{D}$. Consequently, $\mathcal{C}$ is idempotent complete.
\end{lem}
\begin{proof}
  The argument is similar to that of~\cite[\crefnolink{mg:prop:characterization_Shv_infty}]{ho_revisiting_2022}. For $\mathcal{F} \in \Shv_{?,c}(\mathcal{Y})$, we will show that the following conditions are equivalent
  \begin{enumerate}
    \item \label{item:proof:lem:t-structure_fullsubcat_Artinian_taut} $\mathcal{F} \in \mathcal{C}$;
    \item \label{item:proof:lem:t-structure_fullsubcat_Artinian_t} the simple constituents of $\Ho^i(\mathcal{F})$ belongs to $\{\mathcal{L}_i\}_{i\in I}$.
  \end{enumerate}
  Indeed, suppose $\mathcal{F}$ satisfies \ref{item:proof:lem:t-structure_fullsubcat_Artinian_t}, then $\mathcal{F}$ can be built from successive extensions of $\mathcal{L}_i$'s, which then implies that $\mathcal{F} \in \mathcal{C}$, i.e., $\mathcal{F}$ satisfies \ref{item:proof:lem:t-structure_fullsubcat_Artinian_taut}. The other direction is obtained by observing that \ref{item:proof:lem:t-structure_fullsubcat_Artinian_t} is closed under finite direct sums, shifts, and cones.

  Condition \ref{item:proof:lem:t-structure_fullsubcat_Artinian_t} is clearly stable under the perverse truncations of $\Shv_{?,c}(\mathcal{Y})$ and hence, $\mathcal{C}$ inherits the perverse $t$-structure of $\Shv_{?,c}(\mathcal{Y})$. Since this $t$-structure is bounded, idempotent completeness of $\mathcal{C}$ follows from~\cite[Corollary 2.14]{antieau_k-theoretic_2019}.
\end{proof}

\begin{cor} \label{cor:t-structure_Ch^u_G}
  For $? \in \{\gr, \emptyset\}$, that is, we are working in the graded or ungraded setting, the category $\Ch^{\unip,?}_G$ is stable under the perverse truncations of $\Shv_{?,c}(G/G)$, and hence, it inherits the perverse $t$-structure on $\Shv_{?,c}(G/G)$.
\end{cor}
\begin{proof}
  We will argue for the graded case as the ungraded case is identical and is in fact well-known.

  By definition, $\Ch^{\unip,\gr}_G$ is the smallest idempotent complete full $\DG$-subcategory of $\Shv_{\gr,c}(G/G)$ containing the images $\tr(\mathcal{K})$ of $\mathcal{K} \in \Hecke^\gr_G$ that are irreducible perverse sheaves, i.e., they are (grading shifts of) Kazhdan--Lusztig elements. Since $\tr$ is obtained by pulling back along a smooth map and pushing forward along a proper map, $\tr(\mathcal{K})$ decomposes into a finite direct sum of simple perverse sheaves.\footnote{Note that this is not necessarily the case if we work with the mixed sheaves, i.e., when $? = \mixed$.} Let $\mathcal{C}$ be the smallest full $\DG$-subcategory containing these simple perverse sheaves. By \cref{lem:t-structure_fullsubcat_Artinian}, $\mathcal{C}$ inherits the perverse $t$-structure on $\Shv_{\gr,c}(G/G)$. It remains to show that $\mathcal{C} = \Ch^{\unip,\gr}_G$.

  Clearly, $\mathcal{C} \subseteq \Ch^{\unip,\gr}_G$. Moreover, $\Ch^{\unip,\gr}_G$ is the idempotent completion of $\mathcal{C}$. But by~\cref{lem:t-structure_fullsubcat_Artinian}, $\mathcal{C}$ is already idempotent complete. Thus, we are done.
\end{proof}

\subsubsection{Weight structure on $\Ch^{\unip,\gr}_G$}
We will now turn to the weight structure.

\begin{lem} \label{lem:weight-structure_Ch^u_gr_G}
  The category $\Ch^{\unip,\gr}_G$ inherits the weight structure on $\Shv_{\gr,c}(G/G)$. Consequently, the trace functor $\tr: \Hecke^\gr_G \to \Ch^{\unip,\gr}_G$ is weight exact.
\end{lem}
\begin{proof}
  As above, $\Ch^{\unip,\gr}_G$ is generated as a $\DG$-category by $\tr(\mathcal{K})$ where $\mathcal{K} \in \Hecke^\gr_G$ is pure of weight $0$. Note that $\tr(\mathcal{K})$ is pure of weight $0$ in $\Shv_{\gr,c}(G/G)$ since in the horocycle correspondence, we only pull back along a smooth map and push forward along a proper map. Thus, the objects $\tr(\mathcal{K})$ form a negatively self-orthogonal collection in $\Shv_{\gr,c}(G/G)$, and hence, also in $\Ch^{\unip,\gr}$. By~\cite[Corollary 2.1.2]{bondarko_constructing_2018} (see also~\cite[Remark 2.2.6]{elmanto_nilpotent_2021}), they form the weight heart of a weight structure on $\Ch^{\unip,\gr}_G$. It is clear from the definition of this weight structure that it is compatible with the one on $\Shv_{\gr,c}(G/G)$.
\end{proof}

\subsection{Rigidity and consequences}
\label{subsec:rigidity_etc}
In this subsection, we show that the finite Hecke categories are (compactly generated) rigid monoidal categories from which we deduce various interesting consequences. This subsection is mostly independent of the rest of the paper.

\subsubsection{Rigidity of finite Hecke categories}
The rigidity of Hecke categories has been established in~\cite{ben-zvi_character_2009}.

\begin{prop}[{\cite[Theorem 6.2]{ben-zvi_character_2009}}] \label{prop:Hecke_rigid}
  The category $\Hecke^{?,\ren}_G$ is compactly generated rigid monoidal.
\end{prop}
\begin{proof}
  $\Hecke^{?,\ren}$ is compactly generated by construction. Rigidity follows from the same proof as that of~\cite[Theorem 6.2]{ben-zvi_character_2009}. Indeed, the same argument as over there implies that $\Hecke^{?,\ren}_G$ is semi-rigid. But since we are working with the renormalized category of sheaves, all constructible sheaves are compact by definition. In particular, the monoidal unit is compact. Thus, they are rigid, by~\cite[Proposition 3.3]{ben-zvi_character_2009}.
\end{proof}

\begin{rmk}
  As a consequence of \cref{prop:Hecke_rigid}, $\Hecke^?_G$ is rigid as ``small'' monoidal categories in the sense of \cref{rmk:rigid_for_small}.
\end{rmk}

\subsubsection{Rigidity of Drinfel'd centers}
The rigidity of finite Hecke categories implies that for their Drinfel'd centers.

\begin{prop}
  The (braided) monoidal category $\Z(\Hecke^{?,\ren}_G)$ is semi-rigid. Moreover, the (braided) monoidal category $\Z(\Hecke^?_G)$ is rigid (in the sense of \cref{rmk:rigid_for_small}).
\end{prop}
\begin{proof}
  The ``small'' case follows from~\cite[Remark 2.4.2]{kong_center_2018}. For the ``big'' case, first note that by construction, we have a (braided) monoidal functor $\Z(\Hecke^?_G) \hookrightarrow \Z(\Hecke^{?,\ren}_G)$ whose image contains the compact generators of $\Z(\Hecke^{?,\ren}_G)$ (see also \cref{rmk:small_center_is_bigger}). But since all objects of $\Z(\Hecke^?_G)$ are dualizable (as the category is rigid), the compact generators of $\Z(\Hecke^{?,\ren}_G)$ are also dualizable. In other words, $\Z(\Hecke^{?,\ren}_G)$ is compactly generated by dualizable objects, and hence, is semi-rigid by definition (see~\cite[Definition 3.1]{ben-zvi_character_2009}).
\end{proof}

Directly from the \cref{defn:trace_center}, we see that the Drinfel'd center always acts on the original category and its trace. We will use $\otimes$ to denote this action. In the case of Hecke categories, we have the following compatibility.

\begin{cor}
  Using the same notation as in \cref{rmk:trace_center_adj_one_diagram}, let $a \in \Hecke^{?,\ren}_G$ (resp. $a\in \Hecke^?_G$) and $b \in \Ch^{\unip,?,\ren}_G$ (resp. $b\in \ChT^{\unip,?}_G$). Then, we have a natural equivalence
  \begin{align*}
    \tr(a \otimes \z(b))                             & \simeq \tr(a) \otimes b            \\
    (\text{resp.}\qquad \wtilde{\tr}(a\otimes \z(b)) & \simeq \wtilde{\tr}(a) \otimes b).
  \end{align*}
\end{cor}
\begin{proof}
  The statement is equivalent to stating that the functor $\tr$ (resp. $\wtilde{\tr}$) is a functor of $\Ch^{\unip,?,\ren}_G$-(resp. $\ChT^{\unip,?}_G$-)modules. By~\cite[Lemma 3.5]{ben-zvi_character_2009} (resp.~\cite[Remark 2.4.2]{kong_center_2018}), this follows from the fact that the central functor $\z$ is monoidal and that $\Ch^{\unip,\ren}_G$ (resp. $\Ch^\unip_G$) is semi-rigid (resp. rigid).
\end{proof}

\subsubsection{Hochschild homology}
Let $c$ be a dualizable object, with dual $c^\vee$, in a symmetric monoidal category $\mathcal{C}$ (see~\cite[Volume I, Chapter 1, \S4]{gaitsgory_study_2017} for an extended discussion on this topic).\footnote{Note that since we are working in a symmetric monoidal category, the left and right duals coincide.} Then, the trace of $c$, also known as the \emph{Hochschild homology} of $c$ and denoted by $\HH(c)$,\footnote{This is not to be confused with the notion of trace used above. Because of this reason, why we will exclusively use the term Hochschild homology for the type of trace discussed here.} is an element of $\End_\mathcal{C}(\munit_\mathcal{C})$, where $\munit_\mathcal{C}$ is the monoidal unit of $\mathcal{C}$, given by the following composition
\[
  \munit_\mathcal{C} \xrightarrow{\coev_c} c \otimes c^\vee \simeq c^\vee \otimes c \xrightarrow{\ev_c} \munit_\mathcal{C},
\]
where $\coev_c$ and $\ev_c$ are, respectively, the co-evaluation and evaluation maps coming from the duality datum between $c$ and $c^\vee$.

The category $\mathcal{C}$ of interest to us is $\Shv_?(\pt)\hphMod$. Since $\Shv_?(\pt)$ is a rigid symmetric monoidal category, \cite[Volume I, Chapter 1, Propositions 7.3.2 \& 9.4.4]{gaitsgory_study_2017} imply that any compactly generated category in $\Shv_?(\pt)\hphMod$ is dualizable. Given such an object $\mathcal{A}$ in $\Shv_?(\pt)\hphMod$, its Hochschild homology $\HH(\mathcal{A})$ is an object in $\Shv_?(\pt)$.

When $\mathcal{A} \in \Alg(\Shv_?(\pt)\hphMod)$, such that the underlying object $\mathcal{A} \in \Shv_?(\pt)\hphMod$ is dualizable, $\HH(\mathcal{A})$ acquires a natural algebra structure, i.e., $\HH(\mathcal{A}) \in \Alg(\Shv_?(\pt))$, by~\cite[\S3.3.2]{gaitsgory_toy_2022}.

\subsubsection{}
We have the following result from~\cite[Theorem 3.8.5]{gaitsgory_toy_2022}. Note that the meanings of $\HH$ and $\Tr$ in that paper are switched compared to ours. Note also that they prove it for the case where the ambient category is $\DGCatprescont$. However, the same proof carries through for $\Shv_?(\pt)\hphMod$.

\begin{thm}[{\cite[Theorem 3.8.5]{gaitsgory_toy_2022}}] \label{thm:Tr_vs_HH}
  Let $\mathcal{A}$ be a rigid monoidal category in $\Shv_?(\pt)\hphMod$ such that $\mathcal{A}$ is dualizable. Then, there is a canonical equivalence of associative algebras
  \begin{align*}
    \HH(\mathcal{A}) & \simeq \cEnd^?_{\Tr(\mathcal{A})}(\tr(\munit_\mathcal{A})) \in \Alg(\Shv_?(\pt)),
  \end{align*}
  where the superscript $?$ in $\cEnd$, which is a placeholder for $\mixed$ (resp. $\gr$, resp. nothing/$\emptyset$), denotes the $\Shv_\mixed(\pt_0)$-(resp. $\Vect^\gr$-, resp. $\Vect$-)enriched $\Hom$-spaces. See~\cite[\crefnolink{mg:subsubsec:enriched_Hom}]{ho_revisiting_2022}.
\end{thm}

\subsubsection{}
Let $\Spr^?_G = \tr(1_{\Hecke^?_G}) \in \Ch^{\unip,?}_G$ denote the image of the monoidal unit of the finite Hecke category via the natural functor. In the ungraded setting, this is known as the \emph{Grothendieck--Springer sheaf}. In the mixed (resp. graded) case, we will thus refer to this object as the \emph{mixed (resp. graded) Grothendieck--Springer sheaf}. \cref{prop:Hecke_rigid,thm:Tr_vs_HH} implies the following result.

\begin{thm}
  \label{thm:HH_Hecke_as_End}
  We have a natural equivalence of algebras
  \begin{align*}
    \HH(\Hecke^{?,\ren}_G) & \simeq \cEnd^?(\Spr^?_G) \in \Alg(\Shv_?(\pt)).
  \end{align*}
\end{thm}

\section{Formality of the Grothendieck--Springer sheaf} \label{sec:formality_gspringer_sheaf}
The first main result of this section, \cref{thm:formality_Spr^gr}, states that $\cEnd^\gr(\Spr^\gr_G) \in \Alg(\Shv_\gr(\pt)) \simeq \Alg(\Vect^\gr)$ is formal. It is obtained by a spreading argument, taking, as input, the ungraded case, proved by the second-named author using transcendental methods. Combining with a generation result of Lusztig in type $A$, the formality result provides a concrete realization of the category of graded unipotent character sheaves, \cref{thm:explicit_Char_type_A}.

Below, the ungraded case is recalled in \cref{subsec:ungraded_case}, followed by the proof of the graded case in \cref{subsec:graded_case} using a spreading argument. The section then concludes with \cref{subsec:explicit_Ch}, where everything can be made more explicit, especially in type $A$.

\subsection{The ungraded case via transcendental method}
\label{subsec:ungraded_case}
In this subsection, we will work over $\CC$; namely, the geometric objects are stacks over $\CC$ and the sheaves are in $\CC$-vector spaces. Let $p: B/B \to G/G$, $q: B/B \to T/T$, $\Ind_{T \subset B}^G:= p_!q^*$, and its right adjoint $\Res_{T \subset B}^G:= q_*p^!$. We have $\Spr_G = p_! \CC_{B/B} = \Ind_{T \subset B}^G \CC_{T/T} \in \Ch^\unip_G$. Here, all the quotients are obtained by using the adjoint actions.

The category of \emph{all} character sheaves $\Ch_G$ was explicitly computed in~\cite{li_derived_2018} (for simply connected groups) and~\cite{li_derived_2023} (for reductive groups) by using a complex analytic cover of the adjoint quotient $G/G$ to reduce the calculation to generalized Springer theory. We need the following particular statement.

\begin{thm}[{\cite{li_derived_2023}}] \label{thm:li_transcendental_char}
  Let $G$ be a reductive group.
  \begin{enumerate}
    \item \label{item:compute_alg_thm:li_transcendental_char} There is an equivalence of $\DG$-algebras
          \[
            \cEnd(\Spr_G) \simeq (\Ho^*(BT) \otimes \Ho^*(T)) \rotimes \CC[W].
          \]
          In particular, $\cEnd(\Spr_G)$ is a formal $\DG$-algebra.
    \item \label{item:compute_morphism_thm:li_transcendental_char} There is a natural equivalence $\Res_{T \subset B}^G (\Spr_G) \simeq \CC_{T/T}^{\oplus W}$. Moreover, the natural $\DG$-algebra homomorphism $\cEnd(\Spr_G) \to \cEnd(\Res_{T \subset B}^G (\Spr_G)) \simeq \cEnd( \CC_{T/T}^{\oplus W})$ can be expressed explicitly via the commutative diagram
          \[
            \begin{tikzcd}[column sep=small]
              \cEnd(\Spr_G) \ar{r} \ar{d}{\simeq} & \cEnd(\Res_{T\subset B}^G(\Spr_G)) \ar{d}{\simeq} \\
              (\Ho^*(BT) \otimes \Ho^*(T)) \rotimes \CC[W] \ar{r} \ar{d}{\simeq}	& (\Ho^*(BT) \otimes \Ho^*(T)) \ar{d}{\simeq} \otimes \cEnd_{\CC}(\CC[W]) \\
              \cEnd_{(\Ho^*(BT) \otimes \Ho^*(T)) \rotimes \CC[W]} ((\Ho^*(BT) \otimes \Ho^*(T)) \rotimes \CC[W]) \ar[r] & \cEnd_{\Ho^*(BT) \otimes \Ho^*(T)} ((\Ho^*(BT) \otimes \Ho^*(T)) \rotimes \CC[W])
            \end{tikzcd}
          \]
          where the bottom arrow is induced by the restriction of module structure along
          \[
            \Ho^*(BT) \otimes \Ho^*(T) \to (\Ho^*(BT) \otimes \Ho^*(T)) \rotimes \CC[W].
          \]
    \item Similarly, the natural $\DG$-algebra homomorphism $\cEnd(\CC_{T/T}) \to \cEnd(\Spr_G)$ can can be identified with
          \[
            \Ho^*(BT) \otimes \Ho^*(T) \to (\Ho^*(BT) \otimes \Ho^*(T)) \rotimes \CC[W].
          \]
  \end{enumerate}
\end{thm}

Picking an abstract isomorphism of fields $\CC \cong \Qlbar$, we see that \cref{thm:li_transcendental_char} holds equally well for cohomology with coefficients in $\Qlbar$. In the rest of this section, we will work with $\Qlbar$ coefficients.

\subsection{Formality of the graded Grothendieck--Springer sheaf}
\label{subsec:graded_case}
The main goal of the current subsection is a graded version of \cref{thm:li_transcendental_char} formulated in the following theorem whose proof will conclude in \cref{subsubsec:proof_thm:formality_Spr^gr}, after some preliminary preparation.

\begin{thm} \label{thm:formality_Spr^gr}
  Let $G$ be a reductive group over $\Fqbar$, $T\subset B$ and $W$ as before.
  \begin{enumerate}
    \item \label{item:compute_alg_thm:formality_Spr^gr} There is an equivalence of $\DG$-algebras
          \[
            \cEnd^\gr(\Spr^\gr_G) \simeq (\Ho_\gr^*(BT) \otimes \Ho_\gr^*(T)) \rotimes \Qlbar[W] \simeq \Qlbar[\uline{x}, \uline{\theta}] \rotimes \Qlbar[W]
          \]
          where $\uline{x}$ and $\uline{\theta}$ are generators of the cohomology rings $\Ho_\gr^*(BT)$ and $\Ho_\gr^*(T)$, respectively. In particular, $\cEnd^\gr(\Spr^\gr_G)$ is formal. Moreover, $\uline{x}$ and $\uline{\theta}$ have degrees $(2, 2)$ and $(2, 1)$ respectively, with the first (resp. second) index indicating graded (resp. cohomological) degrees.
    \item \label{item:compute_morphism_thm:formality_Spr^gr} There is a natural equivalence $\Res_{T\subset B}^G (\Spr^\gr_G) \simeq \QlbarA_{T/T}^{\oplus W}$. Moreover, after taking cohomology, the natural $\DG$-algebra homomorphism $\cEnd^\gr(\Spr^\gr_G) \to \cEnd^\gr(\Res_{T\subset B}^G(\Spr^\gr_G)) \simeq \cEnd^\gr(\QlbarA_{T/T}^{\oplus W})$ can be expressed explicitly via the commutative diagram
          \[
            \begin{tikzcd}[column sep=tiny]
              \Ho^*(\cEnd^\gr(\Spr^\gr_G)) \ar{r} \ar{d}{\simeq} & \Ho^*(\cEnd^\gr(\Res_{T\subset B}^G(\Spr^\gr_G))) \ar{d}{\simeq} \\
              (\Ho_\gr^*(BT) \otimes \Ho_\gr^*(T)) \rotimes \Qlbar[W] \ar{r} \ar{d}{\simeq}	& (\Ho_\gr^*(BT) \otimes \Ho_\gr^*(T)) \otimes \cEnd_{\Qlbar}(\Qlbar[W]) \ar{d}{\simeq} \\
              \cEnd^\gr_{(\Ho_\gr^*(BT) \otimes \Ho_\gr^*(T)) \rotimes \Qlbar[W]} ((\Ho_\gr^*(BT) \otimes \Ho_\gr^*(T)) \rotimes \Qlbar[W] ) \ar[r] & \cEnd^\gr_{\Ho^*(BT) \otimes \Ho_\gr^*(T)} ((\Ho_\gr^*(BT) \otimes \Ho_\gr^*(T)) \rotimes \Qlbar[W])
            \end{tikzcd}
          \]
          where the bottom arrow is induced by the restriction of module structure along
          \[
            \Ho_\gr^*(BT) \otimes \Ho_\gr^*(T) \to (\Ho_\gr^*(BT) \otimes \Ho_\gr^*(T)) \rotimes \Qlbar[W].
          \]
    \item \label{item:compute_ind_morphism_thm:formality_Spr^gr} Similarly, after taking cohomology, the natural $\DG$-algebra homomorphism $\cEnd^\gr(\QlbarA_{T/T}) \to \cEnd^\gr(\Spr_G)$ can be identified with
          \[
            \Ho_\gr^*(BT) \otimes \Ho_\gr^*(T) \to (\Ho_\gr^*(BT) \otimes \Ho_\gr^*(T)) \rotimes \Qlbar[W].
          \]
  \end{enumerate}
\end{thm}

Combining with \cref{thm:HH_Hecke_as_End}, we obtain the following result.

\begin{cor}
  \label{cor:formality_HH(Hecke)}
  The $\DG$-algebra $\HH(\Hecke^{\gr,\ren}_G)$ is formal, and we have an equivalence of $\DG$-algebras
  \[
    \HH(\Hecke^{\gr,\ren}_G) \simeq \Qlbar[\uline{x}, \uline{\theta}] \rotimes \Qlbar[W],
  \]
  wgere $\uline{x}$ and $\uline{\theta}$ have degrees $(2, 2)$ and $(2, 1)$ respectively, where the first (resp. second) index indicates graded (resp. cohomological) degrees.
\end{cor}

\subsubsection{The strategy for proving \cref{thm:formality_Spr^gr}} \label{subsec:strategy_formality}
The passage from \cref{thm:li_transcendental_char} to \cref{thm:formality_Spr^gr} is via a spreading argument that we will now explain. Let $R$ be a strictly Henselian discrete valuation ring between $\ZZ[1/(lN)]$ and $\CC$, where $N\gg 0$. Let $\imath: \Spec \Fqbar \to \Spec R$ and $\jmath: \Spec \CC \to \Spec R$ be geometric points over the special and generic points of $\Spec R$, respectively. Then, we have the following symmetric monoidal equivalences of categories
\[
  \begin{tikzcd}
    \Vect \simeq \Shv(\Spec \Fqbar) & \ar{l}{\simeq}[swap]{\imath^*} \LS(\Spec R) \ar{r}{\jmath^*}[swap]{\simeq} & \Shv(\Spec \CC) \simeq \Vect,
  \end{tikzcd}
\]
where $\LS(\Spec R)$ denotes the category of $\Qlbar$-local systems on $\Spec R$. This induces equivalences of categories
\[
  \begin{tikzcd}
    \Alg(\Vect) \simeq \Alg(\Shv(\Spec \Fqbar)) & \ar{l}{\simeq}[swap]{\imath^*} \Alg(\LS(\Spec R)) \ar{r}{\jmath^*}[swap]{\simeq} & \Alg(\Shv(\Spec \CC)) \simeq \Alg(\Vect).
  \end{tikzcd} \teq\label{eq:equivalence_of_AlgCats_spreadout}
\]
Thus, if we have an algebra in $\Alg(\Shv(\Spec \Fqbar))$ whose formality we would like to establish, it suffices to produce a natural candidate in $\Alg(\LS(\Spec R))$ whose image under $\jmath^*$ is known to be formal in $\Alg(\Shv(\Spec \CC))$.

The algebra in question is $\cEnd(\Spr_G)$ which we will now spread out to $\Spec R$.

\subsubsection{Spreading out}
Let $G_R$ denote the split reductive group over $\Spec R$ given by the same root datum as that of $G$. Fix $T_R \subset B_R$ a pair of a maximal torus and a Borel subgroup. All the objects considered above have natural relative versions over $R$. We will use the subscript $R$ (resp. $\Fqbar$, resp. $\CC$) in the notation, for example, $\Spr_{G_R}$ (resp. $\Spr_{G_{\Fqbar}}$, resp. $\Spr_{G_{\CC}}$), when it is necessary to emphasize where these objects live over, i.e., over $R$ (resp. $\Fqbar$, resp. $\CC$).

Let $\St_{G_R} = B_R/B_R \times_{G_R/G_R} B_R/B_R$. Then we have the following composition of correspondences
\[
  \begin{tikzcd}[sep=small]
    & & \St_{G_R} \ar{dl}[swap]{r} \ar{dr}{s} \\
    & B_R/B_R \ar{dl}[swap]{q} \ar{dr}{p} & & B_R/B_R \ar{dl}[swap]{p} \ar{dr}{q} \\
    T_R/T_R & & G_R/G_R & & T_R/T_R.
  \end{tikzcd} \teq\label{eq:steinberg_corr}
\]

\begin{lem}[Mackey filtration] \label{lem:Mackey_filtration}
  $\St_{G_R}$ has a locally closed stratification indexed by $w \in W$
  \[
    \St_{G_R} = \bigsqcup_{w\in W} \St_{G_R}^w
  \]
  such that on $\St_{G_R}^w$, \cref{eq:steinberg_corr} is identified with
  \[
    \begin{tikzcd}[sep=small]
      & & B_R^w/B_R^w \ar{dl}[swap]{i} \ar{dr}{\Ad_w} \\
      & B_R/B_R \ar{dl}[swap]{q} \ar{dr}{p} & & B_R/B_R \ar{dl}[swap]{p} \ar{dr}{q} \\
      T_R/T_R & & G_R/G_R & & T_R/T_R,
    \end{tikzcd}
  \]
  where $B_R^w = B_R \cap \Ad_w^{-1}(B_R)$.
\end{lem}
\begin{proof}
  Observe that for any stack $\mathcal{Y}$, $\Map(S^1, \mathcal{Y}) \simeq \mathcal{Y} \times_{\mathcal{Y}\times \mathcal{Y}} \mathcal{Y}$. Hence, if $\mathcal{Z}$ is a locally closed substack of $\mathcal{Y}$, then
  \[
    \Map(S^1, \mathcal{Z}) \simeq \mathcal{Z} \times_{\mathcal{Z} \times \mathcal{Z}} \mathcal{Z} \simeq \mathcal{Z} \times_{\mathcal{Y}\times \mathcal{Y}} \mathcal{Z}
  \]
  is a locally closed substack of $\mathcal{Y} \times_{\mathcal{Y} \times \mathcal{Y}} \mathcal{Y}$. Moreover, if $\mathcal{Y} = \sqcup_i \mathcal{Z}_i$ is a stratification of $\mathcal{Y}$ where the $\mathcal{Z}_i$'s are locally closed substack of $\mathcal{Y}$, then
  \[
    \Map(S^1, \mathcal{Y})
    \simeq \mathcal{Y} \times_{\mathcal{Y} \times \mathcal{Y}} \mathcal{Y}
    \simeq \sqcup_{i, j} \mathcal{Z}_i \times_{\mathcal{Y} \times \mathcal{Y}} \mathcal{Z}_j
    \simeq \sqcup_{i} \mathcal{Z}_i \times_{\mathcal{Y} \times \mathcal{Y}} \mathcal{Z}_i
    \simeq \sqcup_i \mathcal{Z}_i \times_{\mathcal{Z}_i \times \mathcal{Z}_i} \mathcal{Z}_i
    \simeq \sqcup_i \Map(S^1, \mathcal{Z}_i)
  \]
  is a stratification of $\mathcal{Y} \times_{\mathcal{Y} \times \mathcal{Y}} \mathcal{Y}$ by locally closed substacks. Here, the third equivalence is due to the fact that when $i\neq j$, $\mathcal{Z}_i \times_{\mathcal{Y}\times \mathcal{Y}} \mathcal{Z}_j$ is empty.

  Applying this to the case where $\mathcal{Y} = BB_R \times_{BG_R} BB_R$ and the Bruhat stratification, we obtain
  \begin{align*}
    \St_{G_R} & = \Map(S^1, BB_R \times_{BG_R} BB_R) = \Map(S^1, B_R\backslash G_R/B_R) = \sqcup_{w\in W} \Map(S^1, B_R\backslash B_R w B_R/B_R) \\
              & = \sqcup_{w\in W} \Map(S^1, B B_R^w) = \sqcup_{w\in W} B_R^w/B_R^w.
  \end{align*}

\end{proof}

\begin{cor} \label{cor:Res_of_Spr}
  $\Res_{T_R \subset B_R}^{G_R} \Spr_{G_R} \simeq \Res_{T_R \subset B_R}^{G_R} \Ind_{T_R\subset B_R}^{G_R} \QlbarA_{T_R/T_R} \simeq \QlbarA_{T_R/T_R}^{\oplus W}$.
\end{cor}
\begin{proof}
  We first show that $\Res_{T_R \subset B_R}^{G_R} \Spr_{G_R}$ is a local system on $T_R/T_R$ concentrated in cohomological degree $0$. The filtration in \cref{lem:Mackey_filtration} implies that $\Res_{T_R \subset B_R}^{G_R} \Spr_{G_R}$ has a filtration whose associated pieces are given by $(q\circ i)_* \QlbarA_{B_R^w/B_R^w}$ which is simply $\QlbarA_{T_R/T_R}$. Thus, $\Res_{T_R \subset B_R}^{G_R} \Spr_{G_R}$ is a complex of local systems. To see that it concentrates in degree $0$, it suffices to check the stalk at a point in $T_{\CC}/T_{\CC}$. But this is now a well-known statement (see, for example,~\cite[Proposition 3.2]{chen_conjecture_2021}).

  Generically on $T_R/T_R$, the map $\St_{G_R} \to T_R/T_R$ is a trivial W-cover. Thus, on this open dense subset, $\Res_{T_R \subset B_R}^{G_R} \Spr_{G_R}$ is a trivial local system of rank $|W|$. This implies the same statement for $\Res_{T_R \subset B_R}^{G_R} \Spr_{G_R}$ itself as it is the IC-extension of the local system on an open dense subset.
\end{proof}

\begin{cor} \label{cor:End_in_family_LS}
  $\pi_{G_R/G_R, *} \cuEnd(\Spr_{G_R}) \simeq \pi_{\St_{G_R}, *} \QlbarA_{\St_{G_R}}$ is a complex of local systems on $\Spec R$, i.e., an object in $\LS(\Spec R)$. Here, $\cuEnd$ denotes the sheaf of endomorphisms, and for any $R$-stack $\mathcal{Y}_R$, $\pi_{\mathcal{Y}_R}: \mathcal{Y}_R \to \Spec R$ denotes the structure map of $\mathcal{Y}_R$.
\end{cor}
\begin{proof}
  The equivalence $\pi_{G_R/G_R, *} \cuEnd(\Spr_{G_R}) \simeq \pi_{\St_{G_R}, *} \QlbarA_{\St_{G_R}}$ is a standard statement.
  The second claim follows from the first since
  \[
    \pi_{\St_{G_R}, *} \QlbarA_{\St_{G_R}} \simeq \pi_{T_R/T_R, *} q_* r_* \QlbarA_{\St_{G_R}} \simeq \pi_{T_R/T_R,*} \QlbarA_{T_R/T_R}^{\oplus W} \simeq (\pi_{T_R/T_R, *} \QlbarA_{T_R/T_R})^{\oplus W},
  \]
  where the second equivalence is due to \cref{cor:Res_of_Spr}, and
  \[
    \pi_{T_R/T_R, *} \QlbarA_{T_R/T_R, *} \simeq \pi_{T_R \times BT_R, *} \QlbarA_{T_R\times BT_R},
  \]
  which is a complex of local systems.
\end{proof}

\begin{cor} \label{cor:restricting_cEnd}
  $\imath^* \pi_{G_R/G_R, *} \cuEnd(\Spr_{G_R}) \simeq \cEnd(\Spr_{G_{\Fqbar}})$ and $\jmath^* \pi_{G_R/G_R, *} \cuEnd(\Spr_{G_R}) \simeq \cEnd(\Spr_{G_{\CC}})$.
\end{cor}
\begin{proof}
  We will prove the statement for $\jmath^*$ only; the proof for $\imath^*$ is identical. By adjunction, we have a natural map of algebras
  \[
    \jmath^* \pi_{G_R/G_R, *} \cuEnd(\Spr_{G_R}) \to \pi_{G_{\CC}/G_{\CC}, *} \cuEnd(\Spr_{G_{\CC}}) \simeq \cEnd(\Spr_{G_{\CC}}).
  \]
  It suffices to show that this is an equivalence of the underlying chain complexes. But now, we have
  \begin{align*}
    \jmath^* \pi_{G_R/G_R, *} \cuEnd(\Spr_{G_R})
    \simeq \jmath^* \pi_{\St_{G_R}, *} \QlbarA_{\St_{G_R}}
    \simeq \jmath^* \pi_{T_R/T_R, *} \QlbarA_{T_R/T_R}^{\oplus W}
    \simeq \pi_{T_{\CC}/T_{\CC}, *} \QlbarA_{T_{\CC}/T_{\CC}}^{\oplus W}
    \simeq \cEnd(\Spr_{G_{\CC}}).
  \end{align*}
  Here, the third equivalence is by direct computation (as it only involves the torus) and the other equivalences are due to \cref{cor:End_in_family_LS}, which holds equally over $\Fqbar$ and $\CC$.
\end{proof}

We are now ready to prove \cref{thm:formality_Spr^gr}.

\subsubsection{Proof of \cref{thm:formality_Spr^gr}}
\label{subsubsec:proof_thm:formality_Spr^gr}
Applying the discussion in \cref{subsec:strategy_formality} to the object $\pi_{G_R/G_R, *} \cuEnd(\Spr_{G_R})$ using \cref{cor:End_in_family_LS,cor:restricting_cEnd} and the result over $\CC$, \cref{thm:li_transcendental_char}.\ref{item:compute_alg_thm:li_transcendental_char}, we obtain the formality of the algebra $\cEnd(\Spr_{G_{\Fqbar}})$. But since the formality of a graded $\DG$-algebra can be detected at the ungraded level, we obtain the formality of $\cEnd^\gr(\Spr^\gr_G)$ as well. In particular, we have an equivalence of $\DG$-algebras
\[
  \cEnd^\gr(\Spr^\gr_G) \simeq \Ho^*(\cEnd^\gr(\Spr^\gr_G)) \simeq \Qlbar[\uline{x}, \uline{\theta}] \rotimes \Qlbar[W]
\]
for \emph{some} variables $\uline{x}$ and $\uline{\theta}$ whose cohomology degrees we know (i.e., $2$ and $1$, respectively) but whose graded degrees still need to be computed.

We note that in the case of a torus, $\Spr^\gr_T \simeq \QlbarA_{T/T}$ and we have
\[
  \cEnd^\gr(\Spr^\gr_T) \simeq \Ho_\gr^*(BT) \otimes \Ho_\gr^*(T) \simeq \Qlbar[\uline{x}, \uline{\theta}],
\]
where $\uline{x}$ and $\uline{\theta}$ have degrees $(2, 2)$ and $(2, 1)$, respectively. We will deduce the general case from this case.

Over $\Fqbar$, we still have the equivalence $\Res_{T\subset B}^G(\Spr_G) \simeq \QlbarA_{T/T}^{\oplus W}$, which is the precise statement of \cite[Proposition 3.2]{chen_conjecture_2021}. In fact, the ungraded version\footnote{That is, we still work over $\Fqbar$, but forget the grading. In other words, we work with $\Shv$ rather than $\Shv_\gr$.} of \cref{thm:formality_Spr^gr}.\ref{item:compute_morphism_thm:formality_Spr^gr} follows from \cref{eq:equivalence_of_AlgCats_spreadout} and \cref{cor:restricting_cEnd} and the complex version, which is \cref{thm:li_transcendental_char}.\ref{item:compute_morphism_thm:li_transcendental_char}. In particular, at the ungraded level, the natural algebra map
\[
  \cEnd(\Spr_G) \to \cEnd(\Res_{T\subset B}^G \Spr_G)
\]
is identified with
\[
  \Qlbar[\uline{x}, \uline{\theta}] \rotimes \Qlbar[W] \hookrightarrow \cEnd_{\Qlbar}(\Qlbar[W]) \otimes (\Ho^*(BT) \otimes \Ho^*(T)) \simeq \Qlbar[\uline{x}, \uline{\theta}] \otimes \cEnd_{\Qlbar}(\Qlbar[W]).
\]
The grading on the RHS is known, as this is the case of a torus. By degree considerations, $\uline{x}$ and $\uline{\theta}$ on the LHS must have degrees $(2, 2)$ and $(2, 1)$, respectively.

The identifications of the algebra morphisms in \cref{thm:formality_Spr^gr}.\ref{item:compute_morphism_thm:formality_Spr^gr} and \ref{item:compute_ind_morphism_thm:formality_Spr^gr} follow from the ungraded case.
\qed

\begin{rmk}
  In the above, we used the fact that the formality of a graded $\DG$-algebra (i.e., a $\DG$-algebra with an \emph{extra} formal grading) can be detected after forgetting the grading. This can be seen, for example, by applying the main result of~\cite{melani_formality_2019} and the fact that the formation of a certain spectral sequence in loc. cit. commutes with forgetting the formal grading.

  We expect that the same statement also holds for the \emph{formality of a morphism} between formal graded $\DG$-algebras. However, we know neither how to prove this statement nor a place where it is proved. Because of that, the statements appearing in \cref{thm:formality_Spr^gr}.\ref{item:compute_morphism_thm:formality_Spr^gr} and \ref{item:compute_ind_morphism_thm:formality_Spr^gr} are only at the level of cohomology groups.
\end{rmk}

\subsection{An explicit presentation of \texorpdfstring{$\Ch^{\unip,\gr}_G$}{Ch^{u,\gr}_G}}
\label{subsec:explicit_Ch}
The main goal of this subsection is to give an explicit presentation of $\Ch^{\unip,\gr}_G$. We will only fully work out type $A$ case.

\subsubsection{A generation statement in type $A$}
We first recall the following result which is well-known to experts. We include it here for the reader's convenience.

\begin{prop} \label{prop:classification_character_sheaves}
  Let $G$ be a reductive group of type $A$ over $\Fqbar$. Then $\Ch^\unip_G$ is generated by the Grothendieck--Springer sheaf $\Spr_G$.
\end{prop}
\begin{proof}
  Let $\mathcal{F} \in \Ch^\unip_G$. Then, by \cref{cor:t-structure_Ch^u_G}, we know that $\mathcal{F}$ can be built from successive extensions of irreducible character sheaves. It thus suffices to show that when $\mathcal{F}$ is irreducible, it is a direct summand of $\Spr_G$.

  By~\cite[Theorem 4.4.(a)]{lusztig_character_1985-1}, $\mathcal{F}$ is a summand of $\Ind_{L\subset P}^G(\mathcal{K})$ for some cuspidal character sheaf $\mathcal{K}$ on $L$, the Levi factor of some parabolic subgroup $P$ of $G$. The group $A(G) \defeq Z(G)/Z^0(G)$ acts on $\mathcal{F}$ via a character $\chi_\mathcal{F}$ which is the composition $A(G) \twoheadrightarrow A(L) \xrightarrow{\chi_\mathcal{K}} \Qlbar$, where $\chi_\mathcal{K}$ is the action of $A(L)$ on $\mathcal{K}$. Since $\mathcal{F}$ is a unipotent character sheaf, $\chi_\mathcal{F}$ is trivial, and hence, so is $\chi_\mathcal{K}$.

  We claim that if $L$ is not a torus, then $\mathcal{K}$ as above with trivial $\chi_K$ must be zero, which would force $\mathcal{F}$ to be a summand of $\Ind_{T\subset B}^G(\Qlbar) = \Spr_G$. Indeed, by~\cite[(7.1.3)]{lusztig_character_1985-2}, $\mathcal{K} \simeq \IC(\overline{\Sigma},\mathcal{E})[\dim \Sigma]$, where $(\Sigma,\mathcal{E})$ is a cuspidal pair of $L$ in the sense of \cite[Definition 2.4]{lusztig_intersection_1984}. Moreover, by \cite[\S 2.10]{lusztig_intersection_1984}, the classification of cuspidal pairs of $L$ is further reduced to the classification of unipotent cuspidal pairs of $H$, where $H= Z_{L'}(s)$, with the group $L'$ being the simply connected cover of $L/Z^0(L)$ and $s$ a semisimple element in an isolated conjugacy class of $L'$. Now $L$ is of type $A$, hence $H$ is also of type $A$ (and is not isomorphic to a torus). Therefore, by the classification of unipotent cuspidal pairs \cite[beginning of page 206]{lusztig_intersection_1984}, we see that $\chi_K$ must be non-trivial if $\mathcal{K}$ is nonzero (it is easy to see that $\chi_K$'s (non-)triviality is preserved under the above reductions).
\end{proof}

We now turn to the graded version of the result above.
\begin{cor} \label{cor:type_A_generation_graded}
  Let $G$ be a reductive group of type $A$ over $\Fqbar$. Then $\Ch^{\unip,\gr}_G$ is generated by grading shifts $\Spr^\gr_G\lrangle{-}$ of the graded Grothendieck--Springer sheaf. In other words, $\Ch^{\unip,\gr}$ is generated by $\Spr^\gr_G$ as a $\Vect^{\gr,c}$-module category. Consequently, $\Ch^{\unip,\gr,\ren}_G$ is compactly generated by $\Spr^\gr_G$ as a $\Vect^\gr$-module category.
\end{cor}
\begin{proof}
  As above, it suffices to show that any irreducible perverse graded character sheaf $\mathcal{K} \in \Ch^{\unip,\gr}_G$ appears as a direct summand of $\Spr^\gr_G$ up to a grading shift. However, this can be checked after forgetting the grading since simplicity implies purity and moreover, we have a complete description of simple pure graded perverse sheaves by~\cite[\crefnolink{mg:thm:explicit_description_Perv_gr_c_w=nu}]{ho_revisiting_2022}.

  The last statement follows immediately since $\Ch^{\unip,\gr,\ren}_G \simeq \Ind(\Ch^{\unip,\gr}_G)$.
\end{proof}

\begin{thm} \label{thm:explicit_Char_type_A}
  Let $G$ be a reductive group of type $A$ over $\Fqbar$. Then, taking $\cHom^\gr_{\Ch^{\unip,\gr,\ren}_G}(\Spr^\gr_G,-)$ induces an equivalence of $\Vect^\gr$-module categories
  \begin{align*}
    \Ch^{\unip,\gr,\ren}_G
     & \simeq \Qlbar[\uline{x},\uline{\theta}] \rotimes \Qlbar[W]\hphMod(\Vect^\gr)
    \eqdef (\Ho_\gr^*(BT) \otimes \Ho_\gr^*(T)) \rotimes \Qlbar[W]\hphMod^\gr       \\
     & \simeq \HH(\Hecke^{\gr,\ren}_G)\hphMod^\gr.
  \end{align*}
  Taking the full subcategory of compact objects, we get
  \begin{align*}
    \Ch^{\unip,\gr}_G
     & \simeq \Qlbar[\uline{x},\uline{\theta}] \rotimes \Qlbar[W]\hphMod(\Vect^{\gr,c})^{\perf}
    \eqdef (\Ho_\gr^*(BT) \otimes \Ho_\gr^*(T)) \rotimes \Qlbar[W] \hphMod^{\gr,\perf}          \\
     & \simeq \HH(\Hecke^{\gr,\ren}_G)\hphMod^{\gr,\perf}.
  \end{align*}
  Here, $\perf$ denotes the full subcategory consisting of perfect complexes.
\end{thm}
\begin{proof}
  The first statement can be obtained by a standard Barr--Beck--Lurie argument using the generation result in \cref{cor:type_A_generation_graded} and the identification of $\HH(\Hecke^{\gr,\ren}_G)$ in \cref{cor:formality_HH(Hecke)}. The second statement follows from the first by taking compact objects.
\end{proof}

\section{A coherent realization of character sheaves via Hilbert schemes of points on \texorpdfstring{$\CC^2$}{ℂ²}}
\label{sec:char_sheaves_vs_Hilb}
Working in type $A$, i.e., $G=\GL_n$ for some $n$, this section's main result, \cref{thm:2-periodized_hilb_vs_ch_sheaves}, relates the category of graded unipotent character sheaves $\Ch^{\unip,\gr}_G \simeq \Tr(\Hecke^\gr_G)$ and the Hilbert scheme of $n$ points on $\CC^2$. The passage from the categorical trace to the Hilbert scheme of points on $\CC^2$ is given by Koszul duality and a result of Krug~\cite{krug_remarks_2018}, both of which will be reviewed in \cref{subsec:Koszul_dual_recollection,subsec:Hilb_recollection_Krug}. In \cref{subsec:coh_shear_2_period}, we recall the categorical constructions of cohomological shearing and $2$-periodization. The relation between the two allows us to introduce an extra formal grading at the cost of having to work with $2$-periodic objects. All of these results are then used in \cref{subsec:sheared_and_2-periodic_Hilb} to prove the main result of this section, \cref{thm:2-periodized_hilb_vs_ch_sheaves}. Finally, in \cref{subsec:match_objs_char_vs_Hilb}, we match various objects on the two sides, to be used in \cref{sec:GNR_conjecture} to realize the \homflypt{} homology of a link geometrically via Hilbert schemes of points.

In this section, we will work exclusively with the case $G=\GL_n$ and adopt the notation $\Hecke^\gr_n \defeq \Hecke^\gr_{\GL_n} \simeq \Ch^b(\SBim_n)$.

\subsection{Koszul duality: a recollection}
\label{subsec:Koszul_dual_recollection}

We will now recall an equivalence of categories coming out of Koszul duality. The materials presented here are classical, but it can also be viewed as a particularly simple case of the theory developed by Arinkin--Gaitsgory in~\cite{arinkin_singular_2015}; for example, see \S1.3.5 therein.

\subsubsection{A $t$-structure on \texorpdfstring{$\Ch^{\unip,\gr,\ren}_G$}{Ch^{u,gr,ren}_G}}
By shearing, we obtain an equivalence of $\DG$-categories
\[
  \Ch^{\unip,\gr,\ren}_G \simeq \Qlbar[\uline{x},\uline{\theta}] \rotimes \Qlbar[\SymGrp_n] \hphMod^\gr \simeq \Qlbar[\tilde{\uline{x}},\tilde{\uline{\theta}}] \rotimes \Qlbar[\SymGrp_n] \hphMod^\gr
\]
where $\tilde{\uline{x}}$ and $\tilde{\uline{\theta}}$ live in degrees $(2, 0)$ and $(2, -1)$ respectively. Here, we follow the same convention as in \cref{thm:formality_Spr^gr}.

The latter category has a natural $t$-structure which makes the forgetful functor to $\Vect^\gr$ $t$-exact where $\Vect^\gr$ is equipped with the standard $t$-structure. The equivalence of categories above then endows $\Ch^{\unip,\gr,\ren}_G \simeq \Qlbar[\uline{x},\uline{\theta}] \rotimes \Qlbar
  [W]\hphMod^\gr$ with a $t$-structure compatible with the \emph{sheared} $t$-structure on $\Vect^\gr$. Namely, the $t$-heart of $\Vect^\gr$ in this $t$-structure consisting of complexes whose support lie in degrees $(k, k)$ for $k\in \mathbb{Z}$.

Note that the algebra $\Qlbar[\uline{x},\uline{\theta}] \rotimes \Qlbar[\SymGrp_n] $ is connective with respect to this sheared $t$-structure on $\Vect^\gr$, whereas the algebra $\Qlbar[\tilde{\uline{x}},\tilde{\uline{\theta}}] \rotimes \Qlbar[\SymGrp_n]$ is connective is the usual $t$-structure. In fact, by shearing, we have a $t$-exact equivalence of $\DG$-categories
\[
  \Ch^{\unip,\gr,\ren}_G \simeq \QCoh((\mathbb{A}^n \times \mathbb{A}^n[-1])/(\Gm\times \SymGrp_n)),
\]
where $\Gm$ scales all coordinates with degree $2$.

\subsubsection{Coherent objects}
Following~\cite{gaitsgory_ind-coherent_2013}, we define
\[
  \Qlbar[\uline{x},\uline{\theta}] \rotimes \Qlbar[\SymGrp_n] \hphMod^{\gr,\coh} \subseteq \Qlbar[\uline{x},\uline{\theta}]\rotimes \Qlbar[\SymGrp_n] \hphMod^\gr \simeq \Ch^{\unip,\gr,\ren}_G
\]
to be the full subcategory spanned by objects of bounded cohomological amplitude and coherent cohomologies (with respect to the $t$-structure defined in the previous subsubsection). This category is generated by $\Qlbar[\uline{x}] \rotimes \Qlbar[\SymGrp_n]$ under finite (co)limits, idempotent splittings, and grading shifts.

\subsubsection{Relative Koszul duality}
Consider the functor of taking $\uline{\theta}$-invariants
\[
  \inv_\theta \defeq \cHom^\gr_{\Qlbar[\uline{x},
      \uline{\theta}]\rotimes \Qlbar[\SymGrp_n] \hphMod^{\gr,\coh}}(\Qlbar[\uline{x}] \rotimes \Qlbar[\SymGrp_n], -): \Qlbar[\uline{x},\uline{\theta}] \rotimes \Qlbar[\SymGrp_n]\hphMod^{\gr,\coh} \to \Vect^\gr.
\]
Since the category on the left is generated by $\Qlbar[\uline{x}] \rotimes \Qlbar[\SymGrp_n]$, an application of the Barr--Beck--Lurie theorem as in \cref{thm:explicit_Char_type_A} above implies that we have an equivalence of categories
\begin{align*}
   & \alignsep\Qlbar[\uline{x},\uline{\theta}] \rotimes \Qlbar[\SymGrp_n]\hphMod^{\gr,\coh}                                                                                                             \\
   & \xrightarrow[\simeq]{\inv_\theta^\enh} \cEnd^\gr_{\Qlbar[\uline{x}, \uline{\theta}]\rotimes \Qlbar[\SymGrp_n] \hphMod^{\gr,\coh}}(\Qlbar[\uline{x}] \rotimes \Qlbar[\SymGrp_n])\hphMod^{\gr,\perf} \\
   & \simeq \Qlbar[\uline{x},\uline{y}] \rotimes \Qlbar[\SymGrp_n] \hphMod^{\gr,\perf},
\end{align*}
where $\uline{y}$ live in degrees $(-2, 0)$.

We refer to this equivalence as the (relative) Koszul duality.

\subsubsection{} \label{subsubsec:algebras_A_and_B}
To simplify the notation, we let $V$ denote a graded vector space of dimension $n$ living in graded degree $2$ and cohomological degree $0$, equipped with a basis and hence also a permutation representation of $W=\SymGrp_n$. Then,
\[
  \Qlbar[\uline{\theta}] \simeq \Sym V_\theta[-1], \qquad \Qlbar[\uline{x}] \simeq \Sym V_x[-2], \qquad \Qlbar[\uline{y}] \simeq \Sym V_y^\vee,
\]
where the subscripts $x$, $\theta$, and $y$ are there just to make it easy to keep track of the names of the variables. To further simply the notation, we denote
\[
  \mathcal{A}_n \defeq \Sym(V_x[-2] \oplus V_\theta[-1]) \rotimes \Qlbar[\SymGrp_n]
\]
and
\[
  \mathcal{B}_n \defeq \Sym(V_x[-2] \oplus V^\vee_y) \rotimes \Qlbar[\SymGrp_n].
\]
The equivalence of categories above then becomes
\[
  \inv_\theta^\enh: \mathcal{A}_n\hphMod^{\gr,\coh} \xrightarrow{\simeq} \mathcal{B}_n\hphMod^{\gr,\perf}.
\]

\subsubsection{}
\label{subsubsec:KD_sends_trivial_to_free}
Let
\[
  \triv_\theta \defeq \Sym V_x[-2]\rotimes \Qlbar[\SymGrp_n] \in \mathcal{A}_n\hphMod^{\gr,\coh}
\]
where $\uline{\theta}$ act trivially, i.e., by $0$. Directly from the construction, we have
\[
  \inv_\theta^\enh(\triv_\theta) \simeq \mathcal{B}_n.
\]
Moreover, the inverse functor to $\inv_\theta^\enh$ is given by taking $\uline{y}$-coinvariants
\[
  \coinv_y^\enh \defeq \Qlbar \otimes_{\Sym V^\vee_y} -.
\]

\subsubsection{Twisting}
For our purposes, it is convenient to twist $\inv_\theta^\enh$ and $\coinv_y^\enh$ by the sign representation of $\SymGrp_n$. We let
\begin{align*}
  \wtilde{\inv}_\theta^\enh               & \defeq  \Sym^n V^\vee[1] \otimes \inv_\theta^\enh, \\
  \text{and}\qquad \wtilde{\coinv}_y^\enh & \defeq \Sym^n V[-1]\otimes \coinv_y^\enh
\end{align*}
be mutually inverse functors
\[
  \wtilde{\inv}_\theta^\enh: \mathcal{A}_n\hphMod^{\gr,\coh} \rightleftarrows \mathcal{B}_n\hphMod^{\gr,\perf}: \wtilde{\coinv}_y^\enh. \teq\label{eq:twisted_KD_equivalence}
\]

\subsubsection{Restricting to \texorpdfstring{$\mathcal{A}_n\hphMod^{\gr,\perf}$}{A-Mod^{gr,perf}}}
\label{subsubsec:Koszul_restr_to_A_n-Mod^perf}
Let
\[
  \triv_y \defeq \Sym(V_x[-2])\rotimes \Qlbar[\SymGrp_n] \in \mathcal{B}_n\hphMod^{\gr,\perf}
\]
where $\uline{y}$ act trivially, i.e., by $0$. An easy computation using the self-duality of the Koszul complex shows that
\[
  \wtilde{\coinv}_y^\enh(\triv_y) \simeq \mathcal{A}_n,
\]
and hence,
\[
  \wtilde{\inv}_\theta^\enh(\mathcal{A}_n) \simeq \triv_y \simeq \Sym(V_x[-2]) \rotimes \Qlbar[\SymGrp_n].
\]
In other words, $\wtilde{\inv}_\theta^\enh$ simply ``kills'' the variables $\uline{\theta}$ in $\mathcal{A}_n$.

This implies that $\wtilde{\inv}_\theta^\enh$ and $\wtilde{\coinv}_y^\enh$ restrict to a pair of (eponymous) mutually inverse functors
\[
  \wtilde{\inv}_\theta^\enh: \mathcal{A}_n\hphMod^{\gr,\perf} \rightleftarrows \mathcal{B}_n\hphMod^{\gr,\perf}_{\nilp_{\uline{y}}}: \wtilde{\coinv}_y^\enh, \teq\label{eq:KD_A_perf}
\]
where the subscript $\nilp_{\uline{y}}$ denotes the fact that the variables $\uline{y}$ act nilpotently. This is because on the one hand, the LHS is generated by $\mathcal{A}_n$ under finite (co)limits, idempotent splittings, and grading shifts and on the other hand, $\mathcal{A}_n$ is sent to $\triv_y$ on which $\uline{y}$ act by $0$.

The equivalence \cref{eq:KD_A_perf} is in fact a particularly simple case of the theory of singular support developed in~\cite{arinkin_singular_2015} where perfect complexes have $0$-singular support.

\subsubsection{Unipotent character sheaves}
The discussion above gives the following Koszul dual descriptions of $\Ch^{\unip,\gr}_G$ and $\Ch^{\unip,\gr,\ren}_G$ when $G=\GL_n$.

\begin{prop} \label{prop:explicit_Char_type_A_KD}
  Let $G = \GL_n$ be the general linear group of rank $n$ over $\Fqbar$. Then, we have an equivalence of $\Vect^{\gr,c}$-module categories
  \[
    \Ch^{\unip,\gr}_G \simeq \mathcal{B}_n\hphMod^{\gr,\perf}_{\nilp_{\uline{y}}}.
  \]
  Taking $\Ind$-completion, we get an equivalence of $\Vect^\gr$-module categories
  \[
    \Ch^{\unip,\gr,\ren}_G \simeq \mathcal{B}_n\hphMod^{\gr}_{\nilp_{\uline{y}}}.
  \]
\end{prop}
\begin{proof}
  The second statement follows from the first by taking the $\Ind$-completion on both sides. The first statement follows from \cref{thm:explicit_Char_type_A,eq:KD_A_perf}.
\end{proof}

\subsection{Hilbert scheme of points on \texorpdfstring{$\CC^2$}{ℂ²}}
\label{subsec:Hilb_recollection_Krug}
We will now recall the main results of \cite{krug_remarks_2018} which gives an explicit presentation of the categories of quasi-coherent sheaves on the Hilbert schemes of points on $\CC^2$. This will allow us to relate Hilbert scheme of points on $\CC^2$ and $\Ch^{\unip,\gr}_G$ and $\Ch^{\unip,\gr,\ren}_G$ (when $G = \GL_n$), which will be discussed in \cref{subsec:coh_shear_2_period}.

As we have been working over $\Qlbar$ rather than $\CC$, our Hilbert schemes are in fact $\Hilb_n(\Qlbar^2)$ which live over $\Qlbar$. Although the two are isomorphic as abstract varieties due to the isomorphism $\Qlbar \simeq \CC$, we will keep using $\Qlbar$ for the sake of consistency. In fact, all varieties appearing in this subsection are over $\Qlbar$ and hence, we will employ base-field-agnostic notation as much as possible; for example, $\mathbb{A}^2$ will be used to denote the $2$-dimensional affine space over $\Qlbar$.

\subsubsection{The isospectral Hilbert scheme}
Let $\Hilb_n$ denote the Hilbert scheme of points on $\mathbb{A}^2$, $\mathbb{A}^{2n}/\SymGrp_n = (\mathbb{A}^2)^n/\SymGrp_n$ the stack quotient of $\mathbb{A}^{2n}$ by the permutation action of $\SymGrp_n$, and $\mathbb{A}^{2n}\sslash \SymGrp_n$ the GIT quotient. There is a natural map $g: \mathbb{A}^{2n}/\SymGrp_n \to \mathbb{A}^{2n}\sslash \SymGrp_n$ and the Hilbert--Chow morphism $H: \Hilb_n \to \mathbb{A}^{2n}\sslash \SymGrp_n$. The isospectral Hilbert scheme is the \emph{reduced} pullback
\[
  \begin{tikzcd}
    \Iso_n \ar{d} \ar{r} & \mathbb{A}^{2n} \ar{d} \\
    \Hilb_n \ar{r} & \mathbb{A}^{2n}\sslash \SymGrp_n.
  \end{tikzcd}
\]

$\Iso_n$ has a natural $\SymGrp_n$-action compatible with the permutation $\SymGrp_n$-action on $\mathbb{A}^{2n}$ and the trivial $\SymGrp_n$-action on $\Hilb_n$. Thus, we obtain the following commutative diagram
\[
  \begin{tikzcd}
    \Iso_n/\SymGrp_n \ar{d}[swap]{q} \ar{r}{p} & \mathbb{A}^{2n}/\SymGrp_n \ar{d}{g} \\
    \Hilb_n \ar{r}{H} & \mathbb{A}^{2n}\sslash \SymGrp_n.
  \end{tikzcd} \teq\label{eq:iso_spectral_Hilb_stack}
\]

\subsubsection{A derived equivalence}
In the notation above, one of the main results of~\cite{krug_remarks_2018} (which is itself a variant of~\cite{bridgeland_mckay_2001} but has a more convenient form for us) takes the following form.

\begin{thm}[{\cite[Proposition 2.8]{krug_remarks_2018}}]
  The functor $\Psi \defeq q_* p^*: \QCoh(\mathbb{A}^{2n}/\SymGrp_n) \to \QCoh(\Hilb_n)$ is an equivalence of categories, where all functors are derived and where $\QCoh(-)$ denotes the ($\infty$-)derived category of quasi-coherent sheaves.

  Passing to the full subcategories spanned by compact objects, we obtain an equivalence $\Psi: \Perf(\mathbb{A}^{2n}/\SymGrp_n) \xrightarrow{\simeq} \Perf(\Hilb_n)$.
\end{thm}

\subsubsection{$\Gm$-equivariant structures}
Each of the terms in \cref{eq:iso_spectral_Hilb_stack} has a natural $\Gm^2$-action induced by the action of $\Gm^2$ on $\mathbb{A}^2$ where the first (resp. second) factor of $\Gm^2$ scales the first (resp. second) factor of $\mathbb{A}^2$ by weight $1$ (resp. $2$). Moreover, all maps are compatible with this action. We thus obtain an equivariant form of the theorem above. By abuse of notation, we will employ the same notation for the equivariant case as the non-equivariant case above. Note also that the category $\QCoh(\Hilb_n/(\Gm^2))$ below also appears as $\QCoh_{\Gm^2}(\Hilb_n)$ (or some variant thereof) in the literature.

\begin{cor}
  The functor $\Psi \defeq q_* p^*: \QCoh(\mathbb{A}^{2n}/(\SymGrp_n \times \Gm^2)) \to \QCoh(\Hilb_n/(\Gm^2))$ is an equivalence of categories. The same statement applies when restricted to the full subcategories of compact objects (i.e., perfect complexes).
\end{cor}

This equivalence allows us to have an explicit presentation of $\QCoh(\Hilb_n/(\Gm^2))$ as a plain (as opposed to a symmetric monoidal) $\DG$-category.

\begin{prop}
  \label{prop:explicit_description_QCoh(Hilb)_McKay}
  We have an equivalence of categories
  \[
    \QCoh(\Hilb_n/(\Gm^2)) \xrightarrow{\Psi^{-1}} \QCoh(\mathbb{A}^{2n}/(\SymGrp_n\times \Gm^2)) \simeq \Qlbar[\tilde{\uline{x}}, \tilde{\uline{y}}]\rotimes \Qlbar[\SymGrp_n] \hphMod^{\gr,\gr},
  \]
  where the superscript $\gr,\gr$ indicates the fact that we are working with bi-graded complexes. Moreover, $\tilde{\uline{x}}$ and $\tilde{\uline{y}}$ have cohomological degree $0$ and bi-degrees (wrt. the superscript $\gr,\gr$) $(1, 0)$ and $(0, 2)$, respectively.
\end{prop}
\begin{proof}
  Only the last equivalence needs to be proved. Applying Barr--Beck--Lurie to the following adjunction
  \[
    \begin{tikzcd}
      \QCoh(\mathbb{A}^{2n}/(\SymGrp_n\times \Gm^2)) \ar[shift right=\arrdisp]{r}[swap]{p_*} & \ar[shift right=\arrdisp]{l}[swap]{p^*} \QCoh(B(\SymGrp_n \times \Gm^2)) \simeq \Qlbar[\SymGrp_n]\hphMod^{\gr,\gr},
    \end{tikzcd}
  \]
  where $p: \mathbb{A}^{2n}/(\SymGrp_n\times \Gm^2) \to B(\SymGrp_n\times \Gm^2)$ is the natural map, we obtain an equivalence of categories
  \[
    \QCoh(\mathbb{A}^{2n}/(\SymGrp_n \times \Gm^2)) \simeq \Qlbar[\tilde{\uline{x}},\tilde{\uline{y}}] \hphMod(\Qlbar[\SymGrp_n]\hphMod^{\gr,\gr}).
  \]
  But the latter is equivalent to
  \[
    \Qlbar[\tilde{\uline{x}},\tilde{\uline{y}}]\rotimes \Qlbar[\SymGrp_n]\hphMod^{\gr,\gr}
  \]
  and we are done.
\end{proof}

\begin{rmk}
  Instead of $(1, 0)$ and $(0, 2)$, we could have chosen any other bigradings for the $\tilde{\uline{x}}$ and $\tilde{\uline{y}}$ variables. As we will see, these specific choices are made because of the connection to character sheaves and \homflypt{} homology. For example, see \cref{subsubsec:X_Y_gradings_to_X'_Y'} for the relation to the grading on the category of graded unipotent character sheaves.
\end{rmk}

\subsubsection{Supports}
\label{subsubsec:supports_hilb_scheme_defn}
Let $i_{\tilde{\uline{x}}}: \mathbb{A}_{\tilde{\uline{x}}}^n/\SymGrp_n \to \mathbb{A}^{2n}/\SymGrp_n$ denote the closed subscheme defined by the vanishing of all the $\tilde{\uline{y}}$-coordinates. Let $\Hilb_{n,\tilde{\uline{x}}}$ be the pullback
\[
  \begin{tikzcd}
    \Hilb_{n,\tilde{\uline{x}}} \ar{d} \ar{r} & \Hilb_n \ar{d} \\
    \mathbb{A}_{\tilde{\uline{x}}}^n/\SymGrp_n \ar{r}{i_{\tilde{\uline{x}}}} & \mathbb{A}^{2n}/{\SymGrp_n}
  \end{tikzcd}
\]

For any closed embedding of schemes (or in fact, stacks) $Z \subset X$, we let $\QCoh(X)_Z$ denote the full subcategory of $\QCoh(X)$ consisting of objects whose supports lie in $Z$. Similarly, we let $\Perf(X)_Z$ denote the full subcategory of $\Perf(X)$ spanned by objects whose supports lie in $Z$.

\begin{lem}
  The equivalence $\Psi$ above restricts to the following equivalences of categories
  \begin{align*}
    \Qlbar[\tilde{\uline{x}}, \tilde{\uline{y}}] \rotimes \Qlbar[\SymGrp_n]\hphMod^{\gr,\gr}_{\nilp_{\tilde{\uline{y}}}}
     & \simeq \QCoh(\mathbb{A}^{2n}/(\SymGrp_n\times \Gm^2))_{\mathbb{A}_{\tilde{\uline{x}}}^n/(\SymGrp_n\times \Gm^2)} \\
     & \xrightarrow{\Psi} \QCoh(\Hilb_n/(\Gm^2))_{\Hilb_{n,\tilde{\uline{x}}}/(\Gm^2)}.
  \end{align*}
  The same statement applies when restricted to the full subcategories of compact objects (i.e., perfect complexes).
\end{lem}
\begin{proof}
  For brevity's sake, we suppress the $\Gm$-equivariant structures (i.e., the gradings) from the notation.

  From \cref{eq:iso_spectral_Hilb_stack}, we see that $\Psi$ is compatible with the action of $\QCoh(\mathbb{A}^{2n}/\SymGrp_n)$. Now, we note that the full subcategories of interest are cut out precisely by the condition that the variables $\tilde{\uline{y}}$ (from $\mathbb{A}^{2n}/\SymGrp_n$) act nilpotently.

  Alternatively (and more categorically), one can use the theory of support as discussed in, for example, \cite[\S3.5]{arinkin_singular_2015}, to conclude:
  \begin{align*}
    \QCoh(\mathbb{A}^{2n}/\SymGrp_n)_{\mathbb{A}^n_{\tilde{\uline{x}}}/\SymGrp_n}
     & \simeq \QCoh(\mathbb{A}^{2n}/\SymGrp_n) \otimes_{\QCoh(\mathbb{A}^{2n} \sslash \SymGrp_n)} \QCoh(\mathbb{A}^{2n}\sslash \SymGrp_n)_{\mathbb{A}^n_{\tilde{\uline{x}}} \sslash \SymGrp_n} \\
     & \simeq \QCoh(\Hilb_n) \otimes_{\QCoh(\mathbb{A}^{2n} \sslash \SymGrp_n)} \QCoh(\mathbb{A}^{2n}\sslash \SymGrp_n)_{\mathbb{A}^n_{\tilde{\uline{x}}} \sslash \SymGrp_n}                   \\
     & \simeq \QCoh(\Hilb_n)_{\Hilb_{n,\tilde{\uline{x}}}}.
  \end{align*}
\end{proof}

\subsection{Cohomological shearing and \texorpdfstring{$2$}{2}-periodization}
\label{subsec:coh_shear_2_period}
The algebra $\Qlbar[\tilde{\uline{x}}, \tilde{\uline{y}}] \rotimes \Qlbar[\SymGrp_n]$ appearing in \cref{prop:explicit_description_QCoh(Hilb)_McKay} and the algebra $\mathcal{B}_n = \Qlbar[\uline{x}, \uline{y}] \rotimes \Qlbar[\SymGrp_n]$ appearing in \cref{subsubsec:algebras_A_and_B} are almost the same except for the mismatch in the cohomological gradings the number of formal gradings.

In this subsection, we will discuss two general categorical constructions which will allow us, in \cref{subsec:sheared_and_2-periodic_Hilb}, to relate the categories of modules over these two rings, and hence, also to relate the category of (quasi-)coherent sheaves on $\Hilb_n$ and $\Ch^{\unip,\gr,\ren}_G$ for $G=\GL_n$. One of them, known as cohomological shearing, reduces the number of formal gradings whereas the other, $2$-periodization, increases the number of gradings. As it turns out, these two are equivalent, which is the content of \cref{prop:shear_vs_2-periodization} below.

The materials presented here are more or less standard, at least among the experts. We include it here since we cannot find a place where everything is written down in a way that is convenient for us.

The discussion in this subsection applies to both small and large categories. For brevity's sake, we will only discuss the large category case. As usual, the small case can be obtained by passing to compact objects.

\subsubsection{Shearing functors on $\Vect^\gr$}
Consider the shear functors
\begin{align*}
  \sh^\Leftarrow: \Vect^\gr & \to \Vect^\gr       & \sh^\Rightarrow: \Vect^\gr & \to \Vect^\gr         \\
  (V_i)_i                   & \mapsto (V_i[2i])_i & (V_i)_i                    & \mapsto (V_i[-2i])_i.
\end{align*}
Note that these are equivalences of symmetric monoidal categories.

For $A^\gr \in \Vect^\gr$, we let $A^{\gr,\Leftarrow} \defeq \sh^{\Leftarrow}(A^\gr) \in \Vect^\gr$ and $A^{\gr,\Rightarrow} \defeq \sh^{\Rightarrow}(A^\gr) \in \Vect^\gr$. Since $\sh^{\Leftarrow}$ and $\sh^{\Rightarrow}$ are symmetric monoidal equivalences, $A^{\gr,\Leftarrow}$ (resp. $A^{\gr,\Rightarrow}$) is equipped with a natural algebra structure if and only if $A^\gr$ is.

\subsubsection{Shearing a $\Vect^\gr$-module structure}
\label{subsec:shearing_Vectgr_mod_structure}
Since $\Vect^\gr$ is symmetric monoidal, any $\mathcal{C}^\gr \in \VectgrMod$ can be upgraded naturally to a $\Vect^\gr$-bimodule category $\mathcal{C} \in \VectgrBiMod$. We visualize them as $\Vect^\gr$ acting on the left and on the right even though they all come from the left module structure.

\begin{rmk} \label{rmk:ignoring_right_Vectgr_action}
  For us, the right $\Vect^\gr$-module structure is important only for the purpose of taking relative tensor products. We will generally ignore this structure unless the formation of relative tensor products is involved.
\end{rmk}

For each of the left/right module structures, we can precompose the action of $\Vect^\gr$ with $\sh^\Leftarrow$ to obtain a new module structure. For example, ${}^\Rightarrow\mathcal{C}^{\gr,\Leftarrow}$ has the following $\Vect^\gr$-bimodule structure\footnote{Note that the formula below is compatible with the definitions of $\sh^\Rightarrow$ and $\sh^\Leftarrow$ since, for example, $\Qlbar\lrangle{i}$ lives in graded degree $-i$.}
\[
  \Qlbar\lrangle{i} \boxtimes c \boxtimes \Qlbar\lrangle{j} \mapsto c\lrangle{i+j}[2i-2j].
\]
Here, the positions of the arrows with respect to the category indicate which of the left or right module structure we are shearing, and the directions of the arrows indicate which shear we use. A missing arrow indicates that shearing does not occur; for example, $\mathcal{C}^{\gr,\Leftarrow}$ has a sheared $\Vect^\gr$-action on the right while keeping the action on the left unchanged. Note that these modifications only change the action of $\Vect^\gr$ while keeping the underlying $\DG$-category intact.

\subsubsection{}
The following diagram demonstrates how the functors $\sh^\Leftarrow$ and $\sh^\Rightarrow$ interact with $\Vect^\gr$-bimodule structures
\[
  \begin{tikzcd}
    {}^\Rightarrow\Vect^{\gr,\Rightarrow} \ar[shift left=\arrdisp]{r}{\sh^\Leftarrow} & \ar[shift left=\arrdisp]{l}{\sh^\Rightarrow} \Vect^{\gr} \ar[shift left=\arrdisp]{r}{\sh^\Leftarrow} & \ar[shift left=\arrdisp]{l}{\sh^\Rightarrow} {}^\Leftarrow\Vect^{\gr,\Leftarrow},
  \end{tikzcd}
\]
where all functors are equivalences as morphisms in $\VectgrBiMod$. Similarly, we can mix and match the directions of the arrows; For example
\[
  \sh^{\Leftarrow}: \Vect^{\gr,\Rightarrow} \xrightarrow{\simeq} {}^\Leftarrow \Vect^\gr.
\]

Shearing $\Vect^\gr$-bimodule structures can be written naturally in terms of relative tensors. For example,
\[
  \mathcal{C}^{\gr,\Leftarrow} \simeq \mathcal{C}^\gr \otimes_{\Vect^\gr} \Vect^{\gr,\Leftarrow} \qquad\text{and}\qquad {}^\Rightarrow\mathcal{C}^\gr \simeq {}^\Rightarrow\Vect^\gr \otimes_{\Vect^\gr} \mathcal{C}^\gr.
\]

\subsubsection{Shearing $A^\gr\hphMod^\gr$}
Let $A^\gr$ be a graded $\DG$-algebra, i.e., $A^\gr \in \Alg(\Vect^\gr)$. Then, the category of graded $A^\gr$-modules, $A^\gr\hphMod^\gr \defeq A^\gr\hphMod(\Vect^\gr) \in \VectgrMod$, is naturally a $\Vect^\gr$-module category, and hence, also an object in $\VectgrBiMod$, as discussed above. A similar discussion as in the case of $\Vect^\gr$ gives the following equivalences of objects in $\VectgrMod$
\[
  \begin{tikzcd}
    {}^\Rightarrow (A^{\gr,\Rightarrow}\hphMod^{\gr,\Rightarrow}) \ar[shift left=\arrdisp]{r}{\sh^\Leftarrow} & \ar[shift left=\arrdisp]{l}{\sh^\Rightarrow} A^\gr\hphMod^\gr \ar[shift left=\arrdisp]{r}{\sh^\Leftarrow} & \ar[shift left=\arrdisp]{l}{\sh^\Rightarrow} {}^\Leftarrow (A^{\gr,\Leftarrow}\hphMod^{\gr,\Leftarrow}).
  \end{tikzcd}
\]
Moreover, as above, we can mix and match the directions of the arrows.

Replacing $A^\gr$ by $A^{\gr,\Leftarrow}$, we obtain the following equivalences of objects in $\VectgrBiMod$
\[
  \begin{tikzcd}
    A^\gr\hphMod^{\gr,\Rightarrow} \ar[shift left=\arrdisp]{r}{\sh^\Leftarrow} & \ar[shift left=\arrdisp]{l}{\sh^\Rightarrow} {}^\Leftarrow(A^{\gr,\Leftarrow}\hphMod^\gr)
  \end{tikzcd}
\]
and
\[
  \begin{tikzcd}
    A^\gr\hphMod^{\gr,\Leftarrow} \ar[shift left=\arrdisp]{r}{\sh^\Rightarrow} & \ar[shift left=\arrdisp]{l}{\sh^\Leftarrow} {}^{\Rightarrow}(A^{\gr,\Rightarrow}\hphMod^\gr).
  \end{tikzcd}
\]

\subsubsection{Sheared degrading}
\label{subsubsec:sheared_degrading}
Let $\mathcal{C}^\gr \in \VectgrMod$, then we can form the degraded version $\mathcal{C}$ of $\mathcal{C}^\gr$
\[
  \mathcal{C} \defeq \Vect \otimes_{\Vect^\gr} \mathcal{C}^\gr,
\]
where the $\Vect^\gr$-module structure on $\Vect$ is given by the symmetric monoidal functor
\begin{align*}
  \oblv_\gr: \Vect^\gr & \to \Vect               \\
  (V_i)_i              & \mapsto \bigoplus_i V_i
\end{align*}
which forgets the grading on $\Vect^\gr$.

We note that the functor of forgetting the grading is also defined for $\mathcal{C}^\gr$
\[
  \oblv_\gr: \mathcal{C}^\gr \simeq \Vect^\gr \otimes_{\Vect^\gr} \mathcal{C}^\gr \xrightarrow{\oblv_\gr \otimes \id_{\mathcal{C}^\gr}} \Vect\otimes_{\Vect^\gr} \mathcal{C}^\gr \simeq \mathcal{C}.
\]

\subsubsection{}
The sheared degradings of $\mathcal{C}^{\gr}$, denoted by $\mathcal{C}^\Lleftarrow$ and $\mathcal{C}^\Rrightarrow$, are defined to be the usual degradings of ${}^\Rightarrow\mathcal{C}^\gr$ and ${}^\Leftarrow\mathcal{C}^\gr$, respectively (note the reversal of the arrows!)
\[
  \mathcal{C}^{\Rrightarrow} \defeq \Vect \otimes_{\Vect^\gr} {}^\Leftarrow\mathcal{C}^\gr
  \quad\text{and}\quad
  \mathcal{C}^{\Lleftarrow} \defeq \Vect \otimes_{\Vect^\gr} {}^\Rightarrow \mathcal{C}^{\gr}.
\]

\begin{rmk}
  In general, $\mathcal{C}$, $\mathcal{C}^\Lleftarrow$, and $\mathcal{C}^\Rrightarrow$ are different as $\DG$-categories. However, the case where $\mathcal{C}^\gr = \Vect^\gr$ is special as we always have an equivalence of $\DG$-categories
  \[
    \Vect \xrightarrow[\simeq]{V \mapsto V \boxtimes \Qlbar} \Vect\otimes_{\Vect^\gr} {}^\Leftarrow\Vect^\gr \eqdef \Vect^{\Rrightarrow}
    \quad\text{and}\quad
    \Vect \xrightarrow[\simeq]{V\mapsto V\boxtimes \Qlbar}
    \Vect \otimes_{\Vect^\gr} {}^\Rightarrow\Vect^\gr \eqdef \Vect^{\Lleftarrow}. \teq\label{eq:equivalence_cat_shear_degrading_Vect}
  \]

  Since $\Vect$ has a natural $\Vect^\gr$-module category structure, we can shear this action and obtain, for example, $\Vect^\Leftarrow$ and $\Vect^\Rightarrow$. These two are in fact equivalent to $\Vect^\Lleftarrow$ and $\Vect^\Rrightarrow$, respectively, which explains the reversal of the arrows in the definition of the sheared degrading procedure.

  Indeed, the right action of $\Vect^\gr$ on $\Vect^{\Rrightarrow}$ is as follows: $\Qlbar\lrangle{i}[j]$ sends $V\boxtimes \Qlbar$ to
  \[
    V \boxtimes (\Qlbar\lrangle{i}[j]) \simeq V \boxtimes (\Qlbar\lrangle{i}[-2i][j+2i]) \simeq (V[j+2i]) \boxtimes \Qlbar,
  \]
  which corresponds to $V[j+2i] \in \Vect$ under the equivalence of categories \cref{eq:equivalence_cat_shear_degrading_Vect}. But this is precisely the $\Rightarrow$-sheared action of $\Vect^\gr$ on $\Vect$ on the right.
\end{rmk}

\subsubsection{} \label{subsubsec:sheared_degrading_vs_tensors}
We can write sheared degradings in a natural way as relative tensors of categories
\[
  \mathcal{C}^{\Lleftarrow} \simeq \Vect \otimes_{\Vect^\gr} {}^\Rightarrow\mathcal{C}^\gr \simeq \Vect \otimes_{\Vect^\gr} {}^\Rightarrow\Vect^\gr \otimes_{\Vect^\gr} \mathcal{C}^\gr \simeq \Vect^{\Leftarrow} \otimes_{\Vect^\gr} \mathcal{C}^{\gr} \simeq \Vect^{\Lleftarrow} \otimes_{\Vect^\gr} \mathcal{C}^\gr
\]
Similarly,
\[
  \mathcal{C}^{\Rrightarrow} \simeq \Vect^{\Rightarrow} \otimes_{\Vect^\gr} \mathcal{C}^\gr \simeq \Vect^{\Rrightarrow} \otimes_{\Vect^\gr} \mathcal{C}^\gr.
\]

\subsubsection{Sheared degrading $A^\gr\hphMod$}
Let $A^\gr \in \Alg(\Vect^\gr)$ be as above and consider the category of (un-graded) $A^\gr$-modules, denoted by $\oblv_\gr(A^\gr)\hphMod$. By abuse of notation, when confusion is unlikely, we also write $A^\gr \hphMod \defeq \oblv_\gr(A^\gr) \hphMod$.

Clearly, $A^\gr\hphMod \simeq \Vect \otimes_{\Vect^\gr} A^\gr\hphMod^\gr$ is a degrading of $A^\gr\hphMod^\gr$. But we also have sheared degradings as defined in \cref{subsubsec:sheared_degrading}. Namely,
\[
  A^\gr\hphMod^{\Lleftarrow} \defeq \Vect\otimes_{\Vect^\gr} {}^{\Rightarrow} (A^\gr\hphMod^\gr) \xrightarrow[\simeq]{\id_{\Vect} \otimes \sh^{\Leftarrow}} \Vect\otimes_{\Vect^\gr} A^{\gr,\Leftarrow} \hphMod^{\gr,\Leftarrow} \simeq \oblv_\gr(A^{\gr,\Leftarrow}) \hphMod^{\Leftarrow},
\]
and similarly for $A^\gr\hphMod^{\Rightarrow}$.

\emph{Ignoring the right $\Vect^\gr$-actions} (see also \cref{rmk:ignoring_right_Vectgr_action}), we have the following equivalences of $\DG$-categories
\[
  A^\gr\hphMod^\Lleftarrow \simeq \oblv_\gr(A^{\gr,\Leftarrow})\hphMod
  \quad\text{and}\quad
  A^\gr\hphMod^\Rrightarrow \simeq \oblv_\gr(A^{\gr,\Rightarrow})\hphMod. \teq\label{eq:inner_vs_outer_shear}
\]

\subsubsection{$2$-periodization}
\label{subsubsec:2_periodization}
We will now turn to $2$-periodization. We start with the simplest example of $2$-periodization and its interaction with sheared degrading.

Let $\Qlbar[u,u^{-1}] \in \ComAlg(\Vect^\gr)$ be a commutative algebra object where $u$ is in graded degree $1$ and cohomological degree $0$. It is easy to see that the functor of extracting the graded degree $0$ part
\[
  \Qlbar[u, u^{-1}]\hphMod^\gr \xrightarrow[\simeq]{(M_i)_i \mapsto M_0} \Vect
\]
is an equivalence of $\Vect^\gr$-module categories where the actions of $\Vect^\gr$ are the obvious ones (i.e., no shearing). The inverse functor is given by $V \mapsto \Qlbar[u, u^{-1}] \otimes V$.

The discussion above thus gives us the following equivalences
\[
  \begin{tikzcd}[column sep=huge]
    \Vect^\Lleftarrow \simeq \Vect^\Leftarrow \ar[bend right=9]{rr}[swap]{\Qlbar[u, u^{-1}]^\Rightarrow \oldtimes -} \ar{r}{\Qlbar[u, u^{-1}] {\oldtimes} -}[swap]{\simeq} &\Qlbar[u, u^{-1}]\hphMod^{\gr,\Leftarrow} \ar{r}{\sh^\Rightarrow}[swap]{\simeq} & {}^{\Rightarrow}(\Qlbar[u, u^{-1}]^{\Rightarrow} \hphMod^{\gr})
  \end{tikzcd}
\]
and,
\[
  \begin{tikzcd}[column sep=huge]
    \Vect^\Rrightarrow \simeq \Vect^\Rightarrow \ar[bend right=9]{rr}[swap]{\Qlbar[u, u^{-1}]^\Leftarrow \oldtimes -} \ar{r}{\Qlbar[u, u^{-1}] {\oldtimes} -}[swap]{\simeq} &\Qlbar[u, u^{-1}]\hphMod^{\gr,\Rightarrow} \ar{r}{\sh^\Leftarrow}[swap]{\simeq} & {}^{\Leftarrow}(\Qlbar[u, u^{-1}]^{\Leftarrow} \hphMod^{\gr}).
  \end{tikzcd}
\]
Here, by construction, $\Qlbar[u,u^{-1}]^{\Rightarrow}$ (resp. $\Qlbar[u, u^{-1}]^{\Leftarrow}$) is the algebra of Laurent polynomials generated by $1$ element living in graded degree $1$ and cohomological degree $2$ (resp. $-2$).

In particular, we have
\[
  {}^\Leftarrow \Vect^\Lleftarrow \simeq {}^\Leftarrow \Vect^\Leftarrow \simeq \Qlbar[u, u^{-1}]^\Rightarrow\hphMod^\gr
  \quad\text{and}\quad
  {}^\Rightarrow \Vect^\Rrightarrow \simeq {}^\Rightarrow \Vect^\Rightarrow \simeq \Qlbar[u, u^{-1}]^\Leftarrow\hphMod^\gr. \teq\label{eq:2-periodization_special_case}
\]

\subsubsection{}
We will now turn to the general case.

Let $\mathcal{C}^\gr \in \VectgrMod$ be as above. Then, the positive and negative $2$-periodizations of $\mathcal{C}^\gr$ are defined to be
\begin{align*}
  \mathcal{C}^{\gr, \per_+} & \defeq \Qlbar[u,u^{-1}]^\Rightarrow\hphMod^\gr \otimes_{\Vect^\gr} \mathcal{C}^\gr \simeq \Qlbar[u, u^{-1}]^\Rightarrow\hphMod(\mathcal{C}^\gr)
  \shortintertext{and}
  \mathcal{C}^{\gr, \per_-} & \defeq \Qlbar[u,u^{-1}]^\Leftarrow\hphMod^\gr \otimes_{\Vect^\gr} \mathcal{C}^\gr \simeq \Qlbar[u, u^{-1}]^\Leftarrow\hphMod(\mathcal{C}^\gr),
\end{align*}
respectively.

Note that $\Qlbar[u, u^{-1}]^\Leftarrow$ (resp. $\Qlbar[u, u^{-1}]^\Rightarrow$) is equivalent to $\Qlbar[\beta, \beta^{-1}]$ where $\beta$ lives in graded degree $1$ and cohomological degree $-2$ (resp. $2$). The category of module objects over $\Qlbar[\beta,\beta^{-1}]$ is thus the category of (graded) $2$-periodic objects appearing in the literature.

The grading allows one to absorb $2$-periodicity. In fact, the $2$-periodization construction agrees with the sheared degrading construction from above.

\begin{prop} \label{prop:shear_vs_2-periodization}
  Let $\mathcal{C}^\gr \in \VectgrMod$. Then, we have natural equivalences of $\Vect^\gr$-module categories
  \[
    \mathcal{C}^{\gr,\per_-} \simeq {}^\Rightarrow (\mathcal{C}^\Rrightarrow)
    \quad\text{and}\quad
    \mathcal{C}^{\gr,\per_+} \simeq {}^\Leftarrow (\mathcal{C}^{\Lleftarrow}).
  \]
\end{prop}
\begin{proof}
  We treat the case of $\per_-$ as the other case is similar. We have
  \begin{align*}
    \mathcal{C}^{\gr, \per_-}
     & \defeq \Qlbar[u, u^{-1}]^\Leftarrow\hphMod^\gr \otimes_{\Vect^\gr} \mathcal{C}^{\gr}
    \simeq {}^\Rightarrow \Vect^{\Rrightarrow} \otimes_{\Vect^\gr} \mathcal{C}^\gr \simeq {}^\Rightarrow \mathcal{C}^{\Rrightarrow},
  \end{align*}
  where the second and third equivalences are due to \cref{eq:2-periodization_special_case} and \cref{subsubsec:sheared_degrading_vs_tensors}, respectively.
\end{proof}

\begin{rmk} \label{rmk:shear_vs_2-periodization_2-periodic_action}
  We have the following commutative diagrams of symmetric monoidal categories
  \[\begin{tikzcd}
      \Vect^\gr \ar{r}{\oblv_\gr \circ \sh^\Rightarrow} \ar{d}[swap]{\Qlbar[u,u^{-1}]^\Leftarrow \oldtimes -} & \Vect \\
      \Qlbar[u, u^{-1}]^\Leftarrow\hphMod^\gr \ar[equal]{r} \ar{ur}[swap]{M \mapsto M_0} & \Vect^{\gr, \per_-} \ar{u}
    \end{tikzcd}
    \quad
    \begin{tikzcd}
      \Vect^\gr \ar{r}{\oblv_\gr \circ \sh^\Leftarrow} \ar{d}[swap]{\Qlbar[u,u^{-1}]^\Rightarrow \oldtimes -} & \Vect \\
      \Qlbar[u, u^{-1}]^\Rightarrow\hphMod^\gr \ar[equal]{r} \ar{ur}[swap]{M \mapsto M_0} & \Vect^{\gr, \per_+} \ar{u}
    \end{tikzcd}
  \]
  Thus, by construction, the action of $\Vect^\gr$ on ${}^\Rightarrow(\mathcal{C}^\Rrightarrow)$ (resp. ${}^\Leftarrow(\mathcal{C}^\Lleftarrow)$) factors through $\Vect^{\gr,\per_-}$ (resp. $\Vect^{\gr,\per_+}$). On the other hand, $\mathcal{C}^{\gr,\per_-}$ (resp. $\mathcal{C}^{\gr, \per_+}$) are equipped with natural actions of $\Vect^{\gr,\per_-}$ (resp. $\Vect^{\gr,\per_+}$).

  Chasing through the definitions, it is easy to see that the two actions are compatible. In other words, the two equivalences in \cref{prop:shear_vs_2-periodization} are equivalences of $\Vect^{\gr,\per_-}$-(resp. $\Vect^{\gr,\per_+}$-)module categories.
\end{rmk}

\subsubsection{More formal gradings}
In the above, the ``background category'' is $\Vect$ in the sense that all categories are $\DG$-categories, i.e., they are $\Vect$-module categories. Then, we work with categories with one extra grading, i.e., with $\Vect^\gr$-module categories, and consider various categorical operations using this structure. We could have started with $\Vect^\gr$-module categories but still added another grading and worked with $\Vect^{\gr,\gr}$-module categories. Everything discussed above still goes through except that now, everything has one more grading. For example, degrading goes from $\Vect^{\gr,\gr}$-module categories to $\Vect^\gr$-categories.\footnote{In fact, we could, more generally, work with $(\Vect^\gr)^{\otimes n}$-module categories and then add one more grading to get $n+1$ gradings. We will not need this generality in the current paper.} This is the setting that we will work with below.

\subsection{$2$-periodic Hilbert schemes and character sheaves}
\label{subsec:sheared_and_2-periodic_Hilb}

We will now relate Hilbert schemes of points and character sheaves, using the procedure of $2$-periodization discussed above to make the precise statement.

\subsubsection{Introduce an extra grading} \label{subsubsec:introduce_extra_grading}
By default, all of our categories are singly graded in the sense that they are $\Vect^{\gr}$-(or, if we work with the small variant, $\Vect^{\gr,c}$-)module categories. Note that any $\DG$-category can be viewed as a $\Vect^\gr$-module category via the symmetric monoidal functor $\oblv_\gr: \Vect^\gr \to \Vect$. In particular, any $\Vect^{\gr}$-module category is a $\Vect^{\gr,\gr}$-module category, where the second $\Vect^\gr$ acts by forgetting the formal grading. By convention, we use $X$ to denote the default grading and $Y$ the extra grading we just introduced. We will also use the notation $\Vect^{\gr_X, \gr_Y}$ if we want to make it clear which grading convention we are using.

In this subsection, unless otherwise specified, all the shearing and $2$-periodization constructions will be with respect to the $Y$-grading even when we start with $\mathcal{C}^\gr \in \Vect^{\gr_X}\hphMod$, which has only \emph{one} grading. In this case, we simply view $\mathcal{C}$ as an object in $\Vect^{\gr_X,\gr_Y}\hphMod$ where $\Vect^{\gr_Y}$ acts via the symmetric monoidal functor $\oblv_{\gr_Y}: \Vect^{\gr_Y} \to \Vect$. In particular, if we forget the $\Vect^{\gr_Y}$-action, $\mathcal{C}^\gr$, ${}^\Leftarrow \mathcal{C}^\gr$, and ${}^\Rightarrow \mathcal{C}^\gr$ are equivalent as objects in $\Vect^{\gr_X}\hphMod$ (and hence, also as $\DG$-categories). The difference in the $\Vect^{\gr_Y}$-module structures can be seen by looking at $\Vect^{\gr_X,\gr_Y}$-enriched $\Hom$, denoted by $\cHom^{\gr_X, \gr_Y}$.

More precisely, for $c_1, c_2 \in \mathcal{C}^\gr$, on the one hand, we have
\[
  \cHom^{\gr_X,\gr_Y}_{\mathcal{C}^\gr}(c_1, c_2) \simeq \bigoplus_{k} \cHom^{\gr_X}_{\mathcal{C}^\gr}(c_1, c_2)\lrangle{k}_Y \simeq \Qlbar[u, u^{-1}] \otimes \cHom^{\gr_X}_{\mathcal{C}^\gr}(c_1, c_2) \in \Vect^{\gr_X, \gr_Y}.
\]
Here, in the tensor formula, the $\gr_X$-enriched $\Hom$'s are put in $Y$-degree $0$. In other words, it is simply copies of the $\Vect^{\gr_X}$-enriched $\Hom$, put in all $Y$-degrees. On the other hand, a simple argument using adjunctions gives the following identification, which is essentially the same as \cref{prop:shear_vs_2-periodization} but with one extra grading (see also \cref{rmk:shear_vs_2-periodization_2-periodic_action}).

\begin{lem} \label{lem:sheared_mod_cat_str_vs_enriched_Hom}
  In the situation above, we have natural equivalences
  \[
    \cHom^{\gr_X,\gr_Y}_{{}^\Leftarrow\mathcal{C}^\gr}(c_1, c_2) \simeq \sh^{\Rightarrow} (\cHom^{\gr_X,\gr_Y}_{\mathcal{C}^\gr}(c_1, c_2)) \simeq \Qlbar[u, u^{-1}]^\Rightarrow \otimes \cHom_{\mathcal{C}^\gr}^{\gr_X}(c_1, c_2) \in \Vect^{\gr_X, \gr_Y,\per_+}
  \]
  and
  \[
    \cHom^{\gr_X,\gr_Y}_{{}^\Rightarrow\mathcal{C}^\gr}(c_1, c_2) \simeq \sh^{\Leftarrow} (\cHom^{\gr_X,\gr_Y}_{\mathcal{C}^\gr}(c_1, c_2)) \simeq \Qlbar[u, u^{-1}]^\Leftarrow \otimes \cHom_{\mathcal{C}^\gr}^{\gr_X}(c_1, c_2) \in \Vect^{\gr_X, \gr_Y, \per_-},
  \]
  where $\sh^{\Rightarrow}$ and $\sh^{\Leftarrow}$ are with respect to the $Y$-grading. Moreover, in the tensor formula, the $\gr_X$-enriched $\Hom$'s are put in $Y$-degree $0$.
\end{lem}

\subsubsection{The algebra $\mathcal{B}_n$ and variants}
Let $\mathcal{B}_n^{\gr} = \Qlbar[\uline{x}, \uline{y}] \rotimes \Qlbar[\SymGrp_n] \in \Alg(\Vect^{\gr_X,\gr_Y})$ be a bi-graded $\DG$-algebra, where the variables $\uline{x}$ (resp. $\uline{y}$) live in cohomological degree $2$ (resp. $0$) and graded degrees $(2, 1)$ (resp. $(-2, 0)$). Here, the first (resp. second) coordinate represents the $X$-(resp. $Y$-)grading. Note that the $Y$-grading is the extra grading whereas the $X$-grading came from the above.

By construction, we have an equivalence of (singly graded) $\DG$-algebras $\oblv_{\gr_Y}(\mathcal{B}_n^\gr) \simeq \mathcal{B}_n \in \Alg(\Vect^{\gr_X})$, where $\mathcal{B}_n$ is the graded $\DG$-algebra defined in \cref{subsubsec:algebras_A_and_B}. Let $\wtilde{\mathcal{B}}_n^{\gr} \defeq \mathcal{B}_n^{\gr, \Leftarrow}$ and $\wtilde{\mathcal{B}}_n \defeq \oblv_{\gr_Y}(\wtilde{\mathcal{B}}_n^{\gr})$, where the shear is with respect to the $Y$-grading. Then,
\[
  \wtilde{\mathcal{B}}_n^{\gr} \simeq \Qlbar[\tilde{\uline{x}},\tilde{\uline{y}}] \rotimes \Qlbar[\SymGrp_n] \in \Alg(\Vect^{\gr_X, \gr_Y})
\]
where the variables $\tilde{\uline{x}}$ (resp. $\tilde{\uline{y}}$) live in cohomological degree $0$ and graded degrees $(2, 1)$ (resp. $(-2, 0)$). Moreover,
\[
  \wtilde{\mathcal{B}}_n \defeq \oblv_{\gr_Y}(\mathcal{B}_n^\gr) \simeq \Qlbar[\tilde{\uline{x}},\tilde{\uline{y}}]\rotimes \Qlbar[\SymGrp_n] \in \Alg(\Vect^{\gr_X}),
\]
where $\tilde{\uline{x}}$ (resp. $\tilde{\uline{y}}$) live in cohomological degree $0$ and graded degree $2$ (resp. $-2$).

\begin{lem} \label{lem:shear_B_n_modules}
  We have the following equivalence of $\Vect^{\gr_X}$-categories
  \[
    \mathcal{B}_n\hphMod^{\gr_X} \simeq \wtilde{\mathcal{B}}_n^{\gr}\hphMod^{\gr_X,\Rrightarrow}.
  \]
\end{lem}
\begin{proof}
  We have
  \[
    \mathcal{B}_n \simeq \oblv_{\gr_Y}(\mathcal{B}_n^\gr) \simeq \oblv_{\gr_Y} (\wtilde{\mathcal{B}}_n^{\gr,\Rightarrow})
  \]
  and hence, by \cref{eq:inner_vs_outer_shear}, we have the following equivalence of $\Vect^{\gr_X}$-module categories
  \[
    \mathcal{B}_n\hphMod^{\gr_X} \simeq \oblv_{\gr_Y}(\wtilde{\mathcal{B}}_n^{\gr,\Rightarrow}) \hphMod^{\gr_X} \simeq \wtilde{\mathcal{B}}_n^{\gr}\hphMod^{\gr_X,\Rrightarrow}
  \]
  where, by convention, all shears are with respect to the $Y$-grading.
\end{proof}

\begin{cor} \label{cor:sheared_B_n-Mod_vs_2-periodic}
  We have a natural equivalence of $\Vect^{\gr_X,\gr_Y}$-module categories
  \[
    {}^\Rightarrow(\mathcal{B}_n\hphMod^{\gr_X}) \simeq \wtilde{\mathcal{B}}_n^\gr\hphMod^{\gr_X, \gr_Y, \per_-},
  \]
  where the $\Vect^{\gr_Y}$-module structure on the right-hand side is the obvious one. Moreover, the shear on the left-hand side is with respect to the $\Vect^{\gr_Y}$-module structure. Namely,  $\Vect^{\gr_Y}$ acts on the left via
  \[
    \Vect^{\gr_Y} \xrightarrow{\sh^\Leftarrow} \Vect^{\gr_Y} \xrightarrow{\oblv_{\gr_Y}} \Vect,
  \]
  as in \cref{prop:shear_vs_2-periodization}. Consequently, by \cref{rmk:shear_vs_2-periodization_2-periodic_action}, the equivalence above is an equivalence of $\Vect^{\gr_X, \gr_Y, \per_-}$-module categories, where $\per_-$ is with respect to the $Y$-grading.
\end{cor}
\begin{proof}
  We have
  \[
    {}^\Rightarrow(\mathcal{B}_n\hphMod^{\gr_X}) \simeq {}^\Rightarrow(\wtilde{\mathcal{B}}_n^{\gr}\hphMod^{\gr_X,\Rrightarrow}) \simeq \wtilde{\mathcal{B}}_n^{\gr}\hphMod^{\gr_X, \gr_Y, \per_-},
  \]
  where the first and second equivalences are due to \cref{lem:shear_B_n_modules} and \cref{prop:shear_vs_2-periodization}, respectively. By convention, shearing and $2$-periodization are with respect to the $Y$-grading.

  The last statement is due to \cref{rmk:shear_vs_2-periodization_2-periodic_action}.
\end{proof}

\subsubsection{Change of gradings}
\label{subsubsec:X_Y_gradings_to_X'_Y'}
The gradings $X$ and $Y$ used above need to be changed to match with the one used in~\cite{gorsky_algebra_2022}. We denote the new gradings $\wtilde{X}$ and $\wtilde{Y}$ where
\[
  \wtilde{X} = X^2 Y \quad\text{and}\quad \wtilde{Y} = X^{-1}
\]
or equivalently,
\[
  X = \wtilde{Y}^{-1} \quad\text{and}\quad Y = \wtilde{X}\wtilde{Y}^2.
\]

Under this change of gradings, we have
\[
  \wtilde{\mathcal{B}}_n^{\gr} \simeq \Qlbar[\tilde{\uline{x}}, \tilde{\uline{y}}] \rotimes \Qlbar[\SymGrp_n] \in \Alg(\Vect^{\gr_{\wtilde{X}}, \gr_{\wtilde{Y}}}),
\]
where the variables $\tilde{\uline{x}}$ (resp. $\tilde{\uline{y}}$) live in cohomological degree $0$ and graded degrees $(1, 0)$ (resp. $(0, 2)$). We note that this algebra matches precisely with the algebra appearing in \cref{prop:explicit_description_QCoh(Hilb)_McKay}.

\begin{rmk}
  In the rest of this paper, unless otherwise specified, all the shears and $2$-periodizations are still with respect to the $Y$-grading (in the $X,Y$-grading system), even when we are working with the $\wtilde{X},\wtilde{Y}$-grading. As different grading conventions in the \homflypt{} homology theory cause a lot of confusion, at least, to the authors of the current paper, we write down explicit examples below.
\end{rmk}

\subsubsection{$\per_-$ in terms of $\wtilde{X},\wtilde{Y}$-grading} \label{subsubsec:per_-_in_terms_of_wtilde(X)_wtilde(Y)}
In terms of the $\wtilde{X},\wtilde{Y}$-grading, negative periodization $\per_-$ in the $Y$-direction with respect to the $X,Y$-grading (such as the one appearing in \cref{cor:sheared_B_n-Mod_vs_2-periodic}) takes the following form. For $\mathcal{C} \in \Vect^{\gr_X,\gr_Y}\hphMod \simeq \Vect^{\gr_{\wtilde{X}}, \gr_{\wtilde{Y}}}\hphMod$,
\[
  \mathcal{C}^{\per_-} \simeq \Qlbar[\beta, \beta^{-1}]\hphMod(\mathcal{C}),
\]
where $\beta$ lives in cohomological degree $-2$ and graded degrees $(0, 1)$ in terms of the $X,Y$-grading or equivalently, $(1, 2)$ in terms of the $\wtilde{X},\wtilde{Y}$-grading. Since this is the only $2$-periodization that we will use in connection to the $\wtilde{X},\wtilde{Y}$-grading, we will adopt the following notation in the rest of the paper
\[
  \mathcal{C}^{\per} \defeq \mathcal{C}^{\per_-} \simeq \Qlbar[\beta,\beta^{-1}]\hphMod(\mathcal{C}),
\]
where, as above, $\beta$ lives in cohomological degree $-2$ and graded degrees $(0, 1)$ in terms of the $X,Y$-grading or, equivalently, $(1, 2)$ in terms of the $\wtilde{X},\wtilde{Y}$-grading.

Below are a couple of examples that are important to us.

\begin{defn}[$2$-periodizing an object]
  \label{defn:2-periodizing_objects}
  For $\mathcal{C} \in \Vect^{\gr_X,\gr_Y}\hphMod \simeq \Vect^{\gr_{\wtilde{X}},\gr_{\wtilde{Y}}}\hphMod$ and $c\in \mathcal{C}$, the corresponding $2$-periodized object $c^{\per} \in \mathcal{C}^{\per}$ is defined to be $\Qlbar[\beta,\beta^{-1}]\otimes c$.
\end{defn}

\begin{defn}[$2$-periodized Hilbert schemes] \label{defn:2-periodized_hilb}
  Let $\Gm^2$ act on $\mathbb{A}^2$ by scaling the coordinate axes $\tilde{x}$ and $\tilde{y}$ with weights $(1, 0)$ and $(0, 2)$, respectively (in the $\wtilde{X},\wtilde{Y}$-grading). This induces an action of $\Gm^2$ on $\Hilb_n \defeq \Hilb_n(\mathbb{A}^2)$ for any $n$. The $2$-periodized category of quasi-coherent sheaves on $\Hilb_n$ is defined to be
  \[
    \QCoh(\Hilb_n/\Gm^2)^{\per} \defeq \Qlbar[\beta, \beta^{-1}]\hphMod(\QCoh(\Hilb_n/\Gm^2))
  \]
  where $\beta$ lives in cohomological degree $-2$ and bi-graded degrees $(1, 2)$. Taking the compact objects, we define\footnote{Note that since the Hilbert scheme is smooth, we could replace $\Perf$ by $\Coh$.}
  \[
    \Perf(\Hilb_n/\Gm^2)^{\per} \defeq (\Qlbar[\beta, \beta^{-1}]\hphMod(\QCoh(\Hilb_n/\Gm)))^{\perf}.
  \]
  The full subcategories $\QCoh(\Hilb_n/\Gm^2)_{\Hilb_{n,\tilde{\uline{x}}}/\Gm^2}^{\per}$ and $\Perf(\Hilb_n/\Gm^2)_{\Hilb_{n,\tilde{\uline{x}}}/\Gm^2}^{\per}$ are defined analogously (see also \cref{subsubsec:supports_hilb_scheme_defn}).
\end{defn}

\begin{defn}[$2$-periodized bi-graded chain complexes] \label{defn:2_periodized_bigraded_chain}
  We define
  \[
    \Vect^{\gr_{\wtilde{X}}, \gr_{\wtilde{Y}}, \per} \defeq \Qlbar[\beta,\beta^{-1}]\hphMod^{\gr_{\wtilde{X}},\gr_{\wtilde{Y}}}
  \]
  where $\beta$ lives in cohomological degree $-2$ and bi-graded degrees $(1, 2)$. The full subcategory spanned by compact objects is defined analogously and is denoted by
  \[
    \Vect^{\gr_{\wtilde{X}}, \gr_{\wtilde{Y}}, \per, c} \defeq \Qlbar[\beta, \beta^{-1}]\hphMod^{\gr_{\wtilde{X}}, \gr_{\wtilde{Y}}, \perf}.
  \]

  More generally, for any algebra $\mathcal{B} \in \Alg(\Vect^{\gr_{\wtilde{X}}, \gr_{\wtilde{Y}}})$, we let
  \[
    \mathcal{B}\hphMod^{\gr_{\wtilde{X}},\gr_{\wtilde{Y}},\per} \defeq \Qlbar[\beta,\beta^{-1}]\hphMod(\mathcal{B}\hphMod^{\gr_{\wtilde{X}},\gr_{\wtilde{Y}}}),
  \]
  and similarly for the full subcategory spanned by the compact objects
  \[
    \mathcal{B}\hphMod^{\gr_{\wtilde{X}},\gr_{\wtilde{Y}},\per, \perf} \defeq \Qlbar[\beta,\beta^{-1}]\hphMod(\mathcal{B}\hphMod^{\gr_{\wtilde{X}},\gr_{\wtilde{Y}}})^\perf.
  \]
\end{defn}

\subsubsection{Shearing in terms of $\wtilde{X}, \wtilde{Y}$-grading}
Let $\mathcal{C} \in \Vect^{\gr_X,\gr_Y}\hphMod \simeq \Vect^{\gr_{\wtilde{X}}, \gr_{\wtilde{Y}}}\hphMod$. Then, unless otherwise specified, all shears are considered to be with respect to the $Y$-grading. For example, the action of $\Vect^{\gr_{\wtilde{X}}, \gr_{\wtilde{Y}}}$ on ${}^\Rightarrow\mathcal{C}$ is given by
\[
  \Qlbar\lrangle{k}_{\wtilde{X}} \boxtimes \Qlbar\lrangle{l}_{\wtilde{Y}} \boxtimes c \mapsto \Qlbar\lrangle{2k-l}_X \boxtimes \Qlbar\lrangle{k}_Y \boxtimes c \mapsto c\lrangle{2k-l}_X \lrangle{k}_Y [2k].
\]
Here $c\lrangle{m}_X$ denotes the result obtained by $\Qlbar\lrangle{m}_X$ acting on $c$ (and similarly for $c\lrangle{m}_Y$). Moreover, the square bracket denotes cohomological shift, which appears here due to the shear.

By \cref{rmk:shear_vs_2-periodization_2-periodic_action}, ${}^\Rightarrow \mathcal{C}$ is equipped with a natural action of $\Vect^{\gr_{\wtilde{X}}, \gr_{\wtilde{Y}},\per}$ (see \cref{defn:2_periodized_bigraded_chain}).

\subsubsection{Hilbert scheme of points and graded unipotent character sheaves}
The relation between Hilbert schemes of points and graded unipotent character sheaves can now be deduced in a straightforward way from the discussion above.

We start with a $2$-periodized form of \cref{prop:explicit_description_QCoh(Hilb)_McKay}.
\begin{lem} \label{lem:2-periodized_McKay}
  We have the following equivalences of $\Vect^{\gr_{\wtilde{X}}, \gr_{\wtilde{Y}}, \per}$-module categories
  \[
    \QCoh(\Hilb_n/\Gm^2)^{\per} \simeq \wtilde{\mathcal{B}}_n^\gr \hphMod^{\gr_{\wtilde{X}}, \gr_{\wtilde{Y}}, \per}
    \quad\text{and}\quad
    \QCoh(\Hilb_n/\Gm^2)_{\Hilb_{n, \tilde{\uline{x}}}}^{\per} \simeq \wtilde{\mathcal{B}}_n^\gr\hphMod^{\gr_{\wtilde{X}},\gr_{\wtilde{Y}},\per}_{\nilp_{\tilde{\uline{y}}}},
  \]
  and similarly for the full subcategories of perfect complexes.
\end{lem}
\begin{proof}
  The second equivalence follows from the first. The first one follows from applying the $2$-periodization construction $\per$ described in \cref{subsubsec:per_-_in_terms_of_wtilde(X)_wtilde(Y)} to \cref{prop:explicit_description_QCoh(Hilb)_McKay}.
\end{proof}

\begin{thm} \label{thm:2-periodized_hilb_vs_ch_sheaves}
  When $G=\GL_n$, we have the following equivalence of $\Vect^{\gr_{\wtilde{X}}, \gr_{\wtilde{Y}}, \per}$-module categories
  \[
    {}^\Rightarrow \Ch^{\unip,\gr, \ren}_G \simeq
    \wtilde{\mathcal{B}}_n^\gr\hphMod^{\gr_{\wtilde{X}},\gr_{\wtilde{Y}},\per}_{\nilp_{\tilde{\uline{y}}}} \xrightarrow{\Psi^{\per}}
    \QCoh(\Hilb_n/\Gm^2)_{\Hilb_{n, \tilde{\uline{x}}}}^{\per}.
  \]
  And similarly, we have the small variant, which is an equivalence of $\Vect^{\gr_{\wtilde{X}}, \gr_{\wtilde{Y}}, \per, c}$-module categories
  \[
    {}^\Rightarrow \Ch^{\unip,\gr}_G \simeq \wtilde{\mathcal{B}}^\gr_n\hphMod^{\gr_{\wtilde{X}}, \gr_{\wtilde{Y}}, \per}_{\nilp_{\tilde{\uline{y}}}} \xrightarrow{\Psi^{\per}} \Perf(\Hilb_n / \Gm^2)^{\per}_{\Hilb_{n, \tilde{\uline{x}}}}.
  \]
\end{thm}
\begin{proof}
  The second statement follows from the first one by taking the full subcategories spanned by compact objects. For the first one, we have
  \[
    {}^\Rightarrow \Ch^{\unip,\gr,\ren}_G
    \simeq {}^\Rightarrow (\mathcal{B}_n \hphMod^{\gr_X}_{\nilp_{\uline{y}}})
    \simeq \wtilde{\mathcal{B}}_n^\gr \hphMod^{\gr_X,\gr_Y,\per_-}_{\nilp_{\tilde{\uline{y}}}}
    \simeq \wtilde{\mathcal{B}}_n^\gr \hphMod^{\gr_{\wtilde{X}}, \gr_{\wtilde{Y}}, \per}_{\nilp_{\tilde{\uline{y}}}}
    \simeq \QCoh(\Hilb_n/\Gm^2)_{\Hilb_{n, \tilde{\uline{x}}}}^{\per},
  \]
  where
  \begin{itemize}
    \item the first equivalence is due to \cref{prop:explicit_Char_type_A_KD} (after adding a sheared action of $\Vect^{\gr_Y}$);
    \item the second equivalence is due to \cref{cor:sheared_B_n-Mod_vs_2-periodic};
    \item the third equivalence is by definition (see also \cref{subsubsec:per_-_in_terms_of_wtilde(X)_wtilde(Y)});
    \item and finally, the last equivalence is due to \cref{lem:2-periodized_McKay}.
  \end{itemize}
\end{proof}

\subsection{Matching objects}
\label{subsec:match_objs_char_vs_Hilb}
Let $\mathcal{T}$ denote the tautological bundle over  $\Hilb_n$. By abuse of notation, we will also use the same symbol to denote its equivariant version, i.e., a vector bundle of rank $n$ over $\Hilb_n/\Gm^2$. We will now match its exterior powers $\wedge^\alpha \mathcal{T}$ under the equivalences of categories stated in \cref{prop:explicit_description_QCoh(Hilb)_McKay,lem:2-periodized_McKay}. This subsection is merely a recollection of the results proved in~\cite{krug_remarks_2018}, stated in our notation.

\subsubsection{Exterior powers of the permutation representation} \label{subsubsec:exterior_powers_perm}
Let $T \in \wtilde{\mathcal{B}}_n\hphMod^{\gr_{\wtilde{X}}, \gr_{\wtilde{Y}}} \simeq \QCoh(\mathbb{A}^{2n}/(\SymGrp_n\times \Gm^2))$ denote the tensor of the structure sheaf with the permutation representation. More explicitly,
\[
  T \defeq \Qlbar[\tilde{\uline{x}}, \tilde{\uline{y}}] \otimes P \in \wtilde{\mathcal{B}}_n\hphMod^{\gr_{\wtilde{X}}, \gr_{\wtilde{Y}}},
\]
where $P \in \Rep \SymGrp_n$ is the permutation representation. More geometrically, if we let
\[
  p: \mathbb{A}^{2n}/(\SymGrp_n\times \Gm^2) \to B(\SymGrp_n\times \Gm^2),
\]
then $T \defeq p^* (P)$, where $P \in \QCoh(B(\SymGrp_n\times \Gm^2))$ is the permutation of $\SymGrp_n$ put in bi-graded degrees $(0, 0)$ and cohomological degree $0$. More generally, we have $\wedge^\alpha T \simeq p^*(\wedge^\alpha P)$.

We have the following elementary lemma.
\begin{lem} \label{lem:wedge_perm_as_Ind}
  For any $0\leq \alpha\leq n$, we have a natural isomorphism of $\SymGrp_n$-representations
  \[
    \Ind_{\SymGrp_\alpha \times \SymGrp_{n-\alpha}}^{\SymGrp_n}(\sign_\alpha \boxtimes \triv_{n-\alpha}) \simeq \wedge^\alpha P,
  \]
  where $\sign_\alpha$ is the sign representation of $\SymGrp_\alpha$ and $\triv_{n-\alpha}$ is the $1$-dimensional trivial representation of $\SymGrp_{n-\alpha}$.
\end{lem}
\begin{proof}
  Let $v_1, \dots, v_n$ be the natural basis of $P$ as the permutation representation of $\SymGrp_n$ and $v$ a basis of $\sign_\alpha \boxtimes \triv_{n-\alpha}$. We have a natural morphism between $\SymGrp_\alpha \times \SymGrp_{n-\alpha}$-representations
  \begin{align*}
    \sign_\alpha \boxtimes \triv_{n-\alpha} & \to \wedge^\alpha P                      \\
    v                             & \mapsto v_1\wedge \dots \wedge v_\alpha.
  \end{align*}
  Since $v_1 \wedge \dots \wedge v_\alpha$ generates $\wedge^\alpha P$ as a representation of $\SymGrp_n$, we obtain a surjective map
  \[
    \Ind_{\SymGrp_\alpha \times \SymGrp_{n-\alpha}}^{\SymGrp_n}(\sign_\alpha \boxtimes \triv_{n-\alpha}) \to \wedge^\alpha P.
  \]
  Comparing the dimensions, we see that this must be an isomorphism.
\end{proof}

\subsubsection{Wedges of the tautological bundle and wedges of the permutation representation}
We will now match the wedges on both sides.

\begin{thm}[{\cite[Theorem 3.9]{krug_remarks_2018}}] \label{thm:Krug_match_wedges_taut}
  Under the equivalence of categories stated in \cref{prop:explicit_description_QCoh(Hilb)_McKay}, $\wedge^\alpha T \in \wtilde{\mathcal{B}}_n\hphMod^{\gr_{\wtilde{X}}, \gr_{\wtilde{Y}}} \simeq \QCoh(\mathbb{A}^{2n}/(\SymGrp_n\times \Gm^2))$ corresponds to $\wedge^\alpha \mathcal{T} \in \QCoh(\Hilb_n/\Gm^2)$.
\end{thm}
\begin{proof}
  This is \cite[Theorem 3.9]{krug_remarks_2018} in the case where the line bundle involved is just the structure sheaf. \cref{lem:wedge_perm_as_Ind} allows us to match the wedges $\wedge^\alpha T$ with the $W^\alpha(-)$ construction in~\cite[Definition 3.4]{krug_remarks_2018}.
\end{proof}

\subsubsection{$2$-periodization} \label{subsubsec:2-periodization_taut_bundle_matching}
We will need the $2$-periodized version of \cref{thm:Krug_match_wedges_taut} below.

\begin{cor} \label{cor:Krug_match_wedges_taut_2-periodized}
  Under the equivalence of categories stated in \cref{lem:2-periodized_McKay},
  \[
    (\wedge^\alpha T)^{\per} \in \wtilde{\mathcal{B}}_n\hphMod^{\gr_{\wtilde{X}}, \gr_{\wtilde{Y}},\per} \simeq \QCoh(\mathbb{A}^{2n}/(\SymGrp_n\times \Gm^2))^{\per}
  \]
  corresponds to $(\wedge^\alpha \mathcal{T})^{\per} \in \QCoh(\Hilb_n/\Gm^2)^{\per}$. Here, the superscript $\per$ denote the procedure of $2$-periodizing an object (see \cref{defn:2-periodizing_objects}).
\end{cor}
\begin{proof}
  This follows from $2$-periodizing \cref{thm:Krug_match_wedges_taut}.
\end{proof}

\section{\texorpdfstring{\homflypt{}}{HOMFLY-PT} homology via Hilbert schemes of points on \texorpdfstring{$\CC^2$}{ℂ²}}
\label{sec:GNR_conjecture}

We will now use the results above to obtain a proof of a version of a conjecture of Gorsky--\Negut--Rasmussen~\cite{gorsky_flag_2021} which predicts a remarkable way to realize \homflypt{} homology as the cohomology of a certain coherent sheaf on the Hilbert scheme of points on $\CC^2$. The precise version of the conjecture that we prove appeared as~\cite[Conjecture 7.2]{gorsky_algebra_2022} with some modifications emphasized in \cref{rmk:differences_with_GNR}. In particular, we will prove parts (a), (b), and (c) of~\cite[Conjecture 7.2]{gorsky_algebra_2022}.

Below, in \cref{subsec:restr_ext_along_weight_heart}, we use weight structures to describe a general mechanism to turn a cohomological grading into a formal grading and interpret the construction of the \homflypt{} homology theory in these terms (this is what happens with the $a$-degree). It is followed by \cref{subsec:homflypt_via_trace}, where we describe how \homflypt{} can be obtained from the categorical trace of $\Hecke^\gr_n$. In \cref{subsec:explicit_computation_of_functors}, we explicitly compute all the functors involved in the factorization of \homflypt{} homology through the categorical trace of $\Hecke^\gr_n$ in terms of the explicit description of the category of the latter, which was obtained in \cref{thm:explicit_Char_type_A}. In \cref{subsec:homflypt_via_Hilb}, transporting these computations to the Hilbert scheme side using Koszul duality and Krug's result (which are encapsulated in \cref{sec:char_sheaves_vs_Hilb}), we arrive at \cref{thm:GNR_conj_a} (which is part (a) of~\cite[Conjecture 7.2]{gorsky_algebra_2022}) which associates to each braid $\beta$ on $n$ strands a ($2$-periodic) coherent sheaf $\mathcal{F}_\beta$ on $\Hilb_n$ whose global sections (after being tensored by the dual of various exterior powers of the tautological bundle) recover the \homflypt{} homology of the link associated by $\beta$. Finally, in~\cref{subsec:match_actions_Qlbar[x]^{S_n},subsec:support_F_beta} we study the actions of symmetric functions on \homflypt{} homology and the support of $\mathcal{F}_\beta$ on $\Hilb_n$ in relation to the number of connected components the braid closure of $\beta$ has. These appear as \cref{thm:GNR_conj_b,thm:GNR_conj_b_c}, which are parts (b) and (c) of~\cite[Conjecture 7.2]{gorsky_algebra_2022}.

\subsection{Restricting to and extending from the weight heart}
\label{subsec:restr_ext_along_weight_heart}
Let $R_\beta \in \Ch^b(\SBim_n)$ be the Rouquier complex associated to a braid $\beta$. Recall that the \homflypt{} homology of the associated braid closure of $\beta$ is obtained by taking Hochschild homology \emph{term-wise} and then taking the cohomology of the resulting complexes. The final answer thus has \emph{three} gradings: internal grading, complex grading (from $\Ch^b(\SBim_n)$), and the Hochschild grading. From the point of view of homological algebra, this is quite unusual: one does not normally apply a derived functor term-wise to a complex.

There have been many interpretations of this construction, for example,~\cites{webster_geometric_2017, trinh_hecke_2021, bezrukavnikov_monodromic_2022}. In this subsection, we offer a general paradigm to understand this type of construction, using weight structures and then apply it to the case of \homflypt{} homology. This allows one to see more clearly what is happening conceptually, especially on the Koszul dual side, which will be discussed in \cref{subsec:homflypt_via_Hilb} below.

As we will use weight structures in a crucial way, the reader might find it beneficial to take a quick look at~\cite[\crefnolink{mg:subsec:weight_structures_review}]{ho_revisiting_2022} for a quick review of the theory.

\subsubsection{Extending from the weight heart}
One salient feature of the theory of weights is that given an idempotent complete stable $\infty$-category $\mathcal{C}$ with a bounded weight structure, the category $\mathcal{C}$ along with the weight structure can be reconstructed from its weight heart $\mathcal{C}^{\weightheart}$.\footnote{A $\DG$-structure is not necessary for this discussion in this subsection. However, the reader should feel free to replace all occurrences of stable $\infty$-categories with $\DG$-categories if they prefer.} In a precise sense, $\mathcal{C}$ is the free stable $\infty$-category generated by its weight heart. Namely, $\mathcal{C} \simeq (\mathcal{C}^{\weightheart})^\fin$, where for any additive category $\mathcal{A}$, $\mathcal{A}^\fin$ is the stabilization of the sifted-completion of $\mathcal{A}$ (see~\cite{elmanto_nilpotent_2021} for more details). We note that the $(-)^\fin$ procedure generalizes the construction $\Ch^b(-)$ in the following sense: if $\mathcal{A}$ is a classical additive category, then $\mathcal{A}^\fin \simeq \Ch^b(\mathcal{A})$.

\begin{lem}
  \label{lem:restr_to_wheart_eq_cats}
  Let $\mathcal{C}$ and $\mathcal{D}$ be idempotent complete stable $\infty$-categories such that $\mathcal{C}$ is equipped with a bounded weight structure. Then, restricting along $\iota: \mathcal{C}^{\weightheart} \to \mathcal{C}$ gives an equivalence of categories
  \[
    \Fun^\ex(\mathcal{C}, \mathcal{D}) \xrightarrow{\iota^!} \Fun^\add(\mathcal{C}^{\weightheart}, \mathcal{D}), \teq\label{eq:restr_to_wheart_eq_cats}
  \]
  where $\Fun^\ex$ denotes the category of exact functors (i.e., those that preserve all finite (co)limits) and $\Fun^\add$ denotes the category of additive functors (i.e., those that preserve all finite direct sums).
\end{lem}
\begin{proof}
  This is essentially~\cite[Proposition 3.3, part 2]{sosnilo_theorem_2019} (but see also~\cite[Theorem 2.2.9]{elmanto_nilpotent_2021} for a variant). We will thus only indicate the main steps.

  We will first show the following equivalence of categories
  \[
    \Fun^\ex(\mathcal{C}, \Ind(\mathcal{D})) \xrightarrow{\iota^!} \Fun^\add(\mathcal{C}^{\weightheart}, \Ind(\mathcal{D})), \teq\label{eq:restr_to_wheart_eq_cats_ind}
  \]
  where $\Ind(\mathcal{D})$ is the $\Ind$-completion of $\mathcal{D}$. $\iota^!$ admits a left adjoint $\iota_!$ given by left Kan extending along $\iota$. Since left Kan extending along a fully faithful embedding is fully faithful, we see that $\iota_!$ is fully faithful. It remains to show that the $\iota^!$ is also fully faithful. Namely, we want to show that the natural map $\iota_! \iota^! F \to F$ is an equivalence for any $F \in \Fun^\ex(\mathcal{C}, \mathcal{D})$.

  Since $\iota^! \iota_! \iota^! F \simeq \iota^! F$ by full faithfulness of $\iota_!$, we see that $\iota_! \iota^! F(c) \simeq F(c)$ for all $c\in \mathcal{C}^{\weightheart}$. But since both sides of $\iota_! \iota^! F \to F$ preserve all finite (co)limits and since $\mathcal{C}$ is generated by $\mathcal{C}^{\weightheart}$ under finite (co)limits, we are done.

  To obtain \cref{eq:restr_to_wheart_eq_cats} from \cref{eq:restr_to_wheart_eq_cats_ind} we only need to show that $\iota^!$ and $\iota_!$ restrict to the corresponding full subcategories on both sides. But this is immediate using the fact that $\mathcal{C}$ is generated by $\mathcal{C}^{\weightheart}$ under finite (co)limits.
\end{proof}

We will also need the following result, which is a consequence of the lemma above.
\begin{cor}[{\cite[Theorem 2.2.9]{elmanto_nilpotent_2021}}]
  \label{cor:restr_weight_exact_functors}
  Let $\mathcal{C}$ and $\mathcal{D}$ be idempotent complete stable $\infty$-categories such that both $\mathcal{C}$ and $\mathcal{D}$ are each equipped with a bounded weight structure. Then, restricting along $\iota: \mathcal{C}^{\weightheart} \to \mathcal{C}$ gives an equivalence of categories
  \[
    \Fun^{\ex,\wex}(\mathcal{C}, \mathcal{D}) \xrightarrow{\iota^!} \Fun^\add(\mathcal{C}^{\weightheart}, \mathcal{D}^{\weightheart}),
  \]
  where $\wex$ denotes weight exactness, i.e., we consider the category of exact functors which send weight heart to weight heart.

  Consequently, when $\mathcal{D}^{\weightheart}$ is classical, we have the following equivalences,
  \[
    \Fun^{\ex,\wex}(\mathcal{C}, \mathcal{D}) \simeq \Fun^\add(\mathcal{C}^{\weightheart}, \mathcal{D}^{\weightheart}) \simeq \Fun^\add(\ho\mathcal{C}^{\weightheart}, \mathcal{D}^{\weightheart}) \simeq \Fun^{\ex,\wex}(\Ch^b(\ho \mathcal{C}^\weightheart), \mathcal{D}).
  \]
\end{cor}

\subsubsection{Turning cohomological grading into a formal grading}
\label{subsubsec:turning_coh_gr_to_formal_gr}
In situations where $\mathcal{D}$ also has a $t$-structure, the procedure of restricting and extending along $\mathcal{C}^{\weightheart} \hookrightarrow \mathcal{C}$ also allows one to create an extra grading. Indeed, let $\mathcal{C}$ and $\mathcal{D}$ be as in \cref{lem:restr_to_wheart_eq_cats} and $F: \mathcal{C} \to \mathcal{D}$ be an exact functor. Suppose that $\mathcal{D}$ is equipped with a $t$-structure. Then, we obtain a functor
\[
  \Ho^*(F|_{\mathcal{C}^{\weightheart}}): \mathcal{C}^{\weightheart} \xrightarrow{F|_{\mathcal{C}^{\weightheart}}} \mathcal{D} \xrightarrow{\Ho^*} \mathcal{D}^{\theart,\gr},
\]
where $\mathcal{D}^{\theart,\gr}$ is the category of $\mathbb{Z}$-graded objects in $\mathcal{D}^{\theart}$. Applying the $(-)^\fin$ construction and using \cref{cor:restr_weight_exact_functors}, we obtain a functor
\[
  \tilde{F}: \mathcal{C} \to (\mathcal{D}^{\theart,\gr})^\fin \simeq \Ch^b(\mathcal{D}^{\theart,\gr}) \to \Ch^b(\mathcal{D}^{\theart})^\gr,
\]
where the equivalence is due to the fact that $\mathcal{D}^{\theart,\gr}$ is \emph{classical}, i.e., it does not have negative $\Ext$'s. In fact, since $\mathcal{D}^{\theart,\gr}$ is classical, $\mathcal{C}^{\weightheart} \to \mathcal{D}^{\theart,\gr}$ factors through $\mathcal{C}^{\weightheart} \to \ho \mathcal{C}^{\weightheart}$ and hence, $\tilde{F}$ factors as follows
\[
  \begin{tikzcd}[column sep=small]
    \mathcal{C} \ar{dr}[swap]{\wt} \ar{rr}{\tilde{F}} && \Ch^b(\mathcal{D}^{\theart})^\gr \\
    & \Ch^b(\ho \mathcal{C}^{\weightheart}). \ar{ur}[swap]{\bar{F}}
  \end{tikzcd} \teq\label{eq:t_to_gr_factors_through_h_weight_heart}
\]
Here $\ho \mathcal{C}^{\weightheart}$ is the homotopy category of $\mathcal{C}^{\weightheart}$, obtained from $\mathcal{C}^{\weightheart}$ by killing negative $\Ext$'s,\footnote{For an additive $\infty$-category $\mathcal{A}$, the homotopy category $\ho \mathcal{A}$ of $\mathcal{A}$ is obtained from $\mathcal{A}$ by taking $\pi_0$ of the $\Hom$ spaces. This is not to be confused with the homotopy category $\K^b(\mathcal{B}) \defeq \ho \Ch^b(\mathcal{B})$ of bounded chain complexes of a (classical) additive category $\mathcal{B}$, which is also sometimes referred to as, somewhat confusingly, the homotopy category of $\mathcal{B}$.} and $\wt$ denotes the weight complex functor (see~\cite[\crefnolink{mg:rmk:weight_complex_functor_explicit}]{ho_revisiting_2022} for a quick review). Note that $\ho\mathcal{C}^{\weightheart}$ is always a classical category whereas $\mathcal{C}^{\weightheart}$ might have non-trivial $\infty$-categorical structures. Moreover, $\bar{F}$ is given by applying $\Ho^*(F)$ term-wise.

\subsubsection{}
When $\mathcal{C}^{\weightheart}$ is already classical, i.e., $\mathcal{C}^{\weightheart} \simeq \ho\mathcal{C}^{\weightheart}$, the weight complex functor is an equivalence of categories
\[
  \wt: \mathcal{C} \xrightarrow{\simeq} \Ch^b(\mathcal{C}^{\weightheart}).
\]
Under this identification, the functor $\tilde{F}$ is identified with $\bar{F}$, which is given by applying $\Ho^*(F)$ term-wise to the chain complexes in $\mathcal{C}^{\weightheart}$.

\subsubsection{}
The category $\mathcal{D}$ that is of interest to us is $\mathcal{D} \simeq \Vect^\gr$. Given an exact functor $F: \mathcal{C} \to \Vect^\gr$, the construction above thus produces a functor $\tilde{F}: \mathcal{C} \to \Vect^{\gr, \gr}$ whose target is the $\DG$-category of chain complexes in doubly graded vector spaces. We can take cohomology another time and obtain triply graded vector spaces. As we will soon see, further specializing to the case where $\mathcal{C} = \Hecke^\gr_n \simeq \Ch^b(\Hecke^{\gr,\weightheart}_n) \simeq \Ch^b(\SBim_n)$ and $F$ is the functor of taking Hochschild homology, we recover the usual construction of the triply graded \homflypt{} homology theory.

\subsubsection{}
The last observation is a tautology, but a powerful one when combined with \cref{lem:restr_to_wheart_eq_cats} which allows one to compare functors by restricting to the weight heart of the source. This point of view is especially useful when the weight structure is not evident, as is the case on the Hilbert scheme side. Indeed, as we shall see below, \cref{lem:restr_to_wheart_eq_cats} allows us to identify $\HH_\alpha$ and the functor co-represented by the $\alpha$-th exterior of the tautological bundle on the Hilbert scheme.

\subsection{\texorpdfstring{\homflypt{}}{HOMFLY-PT} homology via \texorpdfstring{$\Tr(\Hecke^\gr_n)$}{Tr(Hᵍʳₙ)}}
\label{subsec:homflypt_via_trace}
We will now describe how \homflypt{} homology can be obtained from $\Tr(\Hecke^\gr_n) \simeq \Ch^{\unip,\gr}_G$. Everything in this subsection has essentially been proved in~\cite{webster_geometric_2017}, but in a different language: they use the chromatographic complex construction (which is not known to be functorial) of which the weight complex functor is a functorial upgrade. See also~\cite[\S6.4]{shende_legendrian_2017} for a presentation that is closer to ours (but still does not use weight structures).

\subsubsection{A geometric interpretation of \homflypt{} homology} \label{subsubsec:geometric_interpretation_homflypt}
We will start with a formulation of \homflypt{} homology in terms of graded sheaves on $BB\times_{BG} BB = B\backslash G/B$. Consider the following commutative diagram
\[
  \begin{tikzcd}
    BB \times_{BG} BB \ar{d}{\pi} & \ar{l}[swap]{p} \frac{G}{B} \ar{d}{\pi} \ar{r}{q} & \frac{G}{G} \ar{r} & \pt \\
    BB \times BB & \ar{l}[swap]{p} BB
  \end{tikzcd}
\]
where the square is Cartesian. We have the following functor
\[
  \HH: \Hecke^\gr_n \simeq \Shv_{\gr, c}(BB\times_{BG} BB) \xrightarrow{\pi_*} \Shv_{\gr, c}(BB\times BB) \xrightarrow{p^*} \Shv_{\gr, c}(BB) \xrightarrow{\Co_\gr^*(BB, -)} \Vect^{\gr}. \teq\label{eq:HH_first_way}
\]
Applying the construction in \cref{subsubsec:turning_coh_gr_to_formal_gr}, we obtain a functor
\[
  \wHH: \Hecke^\gr_n \to \Vect^{\gr,\gr}.
\]
Taking cohomology one more time, we get
\[
  \HHH \defeq \Ho^*(\wHH): \Hecke^\gr_n \to \Vect^{\theart, \gr,\gr,\gr},
\]
the category of triply graded vector spaces.

\begin{lem}
  Let $R_\beta \in \Hecke^\gr_n$ be associated to a braid $\beta$. Then up to a change of grading (to be discussed in \cref{subsubsec:AQT_aqt_gradings} below), $\Ho^*(\wHH(R_\beta))$ is the triply graded \homflypt{} link homology of $\beta$.
\end{lem}
\begin{proof}[Proof (Sketch)]
  This is the main theorem of~\cite{webster_geometric_2017}, phrased in our language. See also~\cite[\S6.4]{shende_legendrian_2017} for a presentation that is closer to ours. We will thus only indicate the main steps here.

  The construction in \cref{subsubsec:turning_coh_gr_to_formal_gr} amounts to saying that $\wHH(R_\beta)$ is obtained by applying $p^* \pi_*$ term-wise to $R_\beta \in \Hecke^\gr_n \simeq \Ch^b(\Hecke_n^{\gr,\weightheart})$. By~\cite[\crefnolink{mg:prop:bimod_vs_graded_sheaves_BG_x_BG}]{ho_revisiting_2022}, we know that $\Shv_{\gr,c}(BB\times BB)$ (resp. $\Shv_{\gr,c}(BB)$) corresponds to the category of (perfect complexes of) graded bimodules (resp. modules) over $\Co_\gr^*(BB) \simeq \Co^\gr_*(BT) \simeq \Sym (V[-2])$, where $V$ is the graded vector space of dimension $n$ living in graded degree $2$ and cohomological degree $0$ (see also~\cref{subsubsec:algebras_A_and_B}). Moreover, by~\cite[\crefnolink{mg:subsubsec:functoriality_sheaves_BG}]{ho_revisiting_2022}, pulling back along $p: BB \to BB \times BB$ corresponds to taking Hochschild homology. Finally, by~\cite[\crefnolink{mg:eq:factorization_comparison_functor_SBim}]{ho_revisiting_2022}, $\pi_*|_{\Hecke_n^{\gr,\weightheart}}$ is identified with the forgetful functor
  from Soergel bimodules to bimodules over $\Co^\gr_*(BB) \simeq \Co^\gr_*(BT)$. Thus, the proof concludes.
\end{proof}

\subsubsection{Gradings}
\label{subsubsec:AQT_aqt_gradings}
We let $Q'$ and $A'$ denote the two formal gradings on $\Vect^{\gr,\gr}$ (the target of the functor $\wHH$) and $T'$ the cohomological grading. To keep track of these gradings, we will, from now on, adopt the notation $\Vect^{\gr_{Q'}, \gr_{A'}}$ rather than simply $\Vect^{\gr,\gr}$ as above. We use $[n]_{Q'}$, $\lrangle{n}_{A'}$, and $\{n\}_{T'}$ to denote the grading shifts.

The relations between the $Q'$, $A'$, and $T'$ gradings and the $Q$, $A$, and $T$ of~\cite{gorsky_algebra_2022} are as follows (written multiplicatively)
\[
  Q' = A^{-1}Q, \qquad A' = A, \qquad T' = T.
\]
Note, however, that \cite{gorsky_algebra_2022} uses Hochschild \emph{cohomology} rather than homology to define the \homflypt{} link homology. Thus, for the unknot we get $\frac{1+Q'^2 A'}{1-Q'^2 A'^2} = \frac{1 + Q^2 A^{-1}}{1-Q^2}$ whereas they get $\frac{1+Q^{-2}A}{1-Q^2}$.

We also use the $q$, $t$, and $a$ gradings, where
\[
  q = Q^2 = Q'^2 A'^2, \qquad a=Q^2 A^{-1}=Q'^2 A', \qquad t=T^2Q^{-2}=T'^2 Q'^{-2} A'^{-2}.
\]
In terms of $q$, $a$, and $t$, the unknot gives $\frac{1+a}{1-q}$.

\subsubsection{}
On $\Hecke^\gr_n$, $\Ch^{\unip,\gr}_G \simeq \mathcal{A}_n\hphMod^{\gr,\perf}$, and $\mathcal{A}_n\hphMod^{\gr,\coh}$ there are two grading shift operators for each $n\in \mathbb{Z}$, the cohomological shift $[n]$ and the grading shift $\lrangle{n}$, and both are all compatible with the various functors between these categories. In terms of the grading convention of \cref{subsec:sheared_and_2-periodic_Hilb}, this formal grading is the $X$-grading over there.

\subsubsection{}
Unwinding the construction, we see that the functor $\wHH$  relates the grading shift operators on $\Hecke^\gr_n$ and those on $\Vect^{\gr_{Q'},\gr_{A'}}$ as follows\footnote{The first and last are easiest to see from the definition, from which the middle one could be deduced. Alternatively, the middle one can also be seen by looking at the weight complex functor.}
\[
  [n] \rightsquigarrow \{n\}_{T'}, \qquad \lrangle{n} \rightsquigarrow [n]_{A'}\lrangle{n}_{Q'}\{-n\}_{T'}, \qquad [n]\lrangle{n} \rightsquigarrow [n]_{A'} \lrangle{n}_{Q'}.
\]

\begin{rmk} \label{rmk:difference_QAT_Q'A'T'}
  The difference between $Q',A',T'$ and $Q,A,T$ comes from the fact that the identification $\Hecke^\gr_n \defeq \Shv_{\gr,c}(B\backslash G/B) \simeq \Ch^b(\SBim_n)$ involves a shear. More concretely, objects in $\Hecke^\gr_n$ are most naturally viewed as graded $\Co_\gr^*(BB)$-bimodules, where the generators of $\Co_\gr^*(BB) \simeq \Co_\gr^*(BT) \simeq \Sym V[-2]$ live in graded (resp. cohomological) degree $2$ (resp. $2$). On the other hand, the generators of the polynomial ring used to define Soergel bimodules are, by convention, put in graded (resp. cohomological) degrees $2$ (resp. $0$). See~\cite[\crefnolink{mg:sec:Hecke_categories}]{ho_revisiting_2022} for a more in-depth discussion.
\end{rmk}

\subsubsection{\texorpdfstring{\homflypt{}}{HOMFLY-PT} homology via truncated categorical trace} \label{subsubsec:homflypt_via_truncated_trace}
By smooth base change, we see that the functor $\HH$ of \cref{eq:HH_first_way} admits an alternative construction
\[
  \HH: \Hecke^\gr_n \simeq \Shv_{\gr,c}(BB\times_{BG} BB) \xrightarrow{p^*} \Shv_{\gr,c}\left( \frac{G}{B} \right) \xrightarrow{q_* \simeq q_!} \Shv_{\gr,c}\left( \frac{G}{G} \right) \xrightarrow{\Co_\gr^*(G/G,-)} \Vect^\gr,
\]
which is the same as
\[
  \HH: \Hecke^\gr_n \xrightarrow{\tr} \Tr(\Hecke^\gr_n) \simeq \Ch^{\unip,\gr}_G \hookrightarrow \Shv_{\gr,c}\left(G/G\right) \xrightarrow{\Co_\gr^*(G/G, -)} \Vect^\gr. \teq\label{eq:HH_via_trace}
\]

Since the functor $\tr: \Hecke^\gr_n \to \Ch^{\unip,\gr}_G$ is weight exact, it commutes with the weight complex functors. Thus, to construct $\wHH$ from $\HH$, we can apply the construction at \cref{subsubsec:turning_coh_gr_to_formal_gr} to the last step of \cref{eq:HH_via_trace}. As a result, we obtain the following factorization of $\wHH$
\[
  \begin{tikzcd}
    \Hecke^\gr_n \ar[bend left=12]{rr}{\wHH} \ar{r}[swap]{\tr} & \Ch^{\unip,\gr}_G \ar{r}[swap]{\Gamma} \ar{d}{\wt} & \Vect^{\gr_{Q'},\gr_{A'}}, \\
    & \Ch^b(\ho \Ch^{\unip,\gr, \weightheart}_G) \ar{ur}[swap]{\lbar{\Gamma}}
  \end{tikzcd} \teq\label{eq:wHH_via_tr}
\]
where $\Gamma$ is obtained by applying the construction at \cref{subsubsec:turning_coh_gr_to_formal_gr} to the functor $\Co_\gr^*(G/G,-)$ and $\lbar{\Gamma}$ is the induced functor (see \cref{eq:t_to_gr_factors_through_h_weight_heart}). The category $\Ch^b(\ho \Ch^{\unip,\gr,\weightheart}_G)$ (resp. $\ho \Ch^{\unip,\gr,\weightheart}_G$) can be thought of as a truncation of $\Ch^{\unip,\gr}_G$ (resp. $\Ch^{\unip,\gr,\weightheart}$), which justifies the name \emph{truncated trace}. It is also referred to as the \emph{underived trace} and is denoted by $\Tr_0(\Hecke^\gr_G)$ in~\cite{gorsky_derived_2021}.

\subsection{Explicit computation of functors}
\label{subsec:explicit_computation_of_functors}
In what follows, we will compute the functors $\Gamma$, $\lbar{\Gamma}$, and $\wt$ in \cref{eq:HH_via_trace} explicitly in terms of the identification
\[
  \Ch^{\unip,\gr}_G \simeq \Qlbar[\uline{x},\uline{\theta}]\rotimes \Qlbar[\SymGrp_n] \hphMod^{\gr,\perf} = \mathcal{A}_n\hphMod^{\gr,\perf}
\]
of \cref{thm:explicit_Char_type_A}. In particular, we apply \cref{lem:restr_to_wheart_eq_cats} to deduce a co-representability statement similar to the main result of~\cite{bezrukavnikov_monodromic_2022,trinh_hecke_2021}.\footnote{It is in fact possible to reprove the co-representability statement of~\cite{bezrukavnikov_monodromic_2022} using the computation presented here. The details will appear in a forthcoming paper.} Moreover, our computation identifies the weight complex functor $\wt$ with the functor of restricting to the nilpotent cone, explaining conceptually why the latter appears in~\cite{trinh_hecke_2021}; see \cref{subsubsec:restricting_nilp_cone}.

\subsubsection{Identifying the weight complex functor $\wt$ for \texorpdfstring{$\Ch^{\unip,\gr}_G$}{Ch^u,gr_G}}
\label{subsubsec:identify_wt_cmplx_for_Ch_unip_gr}
We will now give a concrete interpretation of the functor $\wt$ appearing in \cref{eq:wHH_via_tr} in terms of modules over $\mathcal{A}_n$. Consider the following functor
\[
  \wtilde{\inv}_\theta \defeq \Sym^n V^\vee[1] \otimes \inv_\theta: \mathcal{A}_n\hphMod^{\gr,\perf} \to \lbar{\mathcal{A}}_n\hphMod^{\gr,\perf},
\]
where $\inv_\theta$ is defined in \cref{subsec:Koszul_dual_recollection} and where we define
\[
  \lbar{\mathcal{A}}_n \defeq \Qlbar[\uline{x}]\rotimes \Qlbar[\SymGrp_n] \simeq \Sym V_x[-2]\rotimes \Qlbar[\SymGrp_n].
\]
See also \cref{subsubsec:algebras_A_and_B} for the definitions of $V$ and $\mathcal{A}_n$.

Note that $\inv_\theta$ is originally defined on $\mathcal{A}_n\hphMod^{\gr,\coh}$. So strictly speaking, in the above, we restrict this functor to the full subcategory of perfect complexes. It is easy to see that $\wtilde{\inv}_\theta$ is given precisely by applying $\wtilde{\inv}_\theta^\enh$ followed by the functor of forgetting the action by the variables $\uline{y}$.

The goal now is to identify $\wt: \Ch^{\unip,\gr}_G \to \Ch^b(\ho \Ch^{\unip,\gr,\weightheart}_G)$ with $\wtilde{\inv}_\theta$.

\subsubsection{}
By \cref{lem:weight-structure_Ch^u_gr_G}, we see that in type $A$, $\Ch^{\unip,\gr}_G \simeq \mathcal{A}_n\hphMod^{\gr,\perf}$ is equipped with a weight structure whose weight heart is spanned under finite direct sums of direct summands of $\Spr^\gr_G[k]\lrangle{k}$ (for all $k\in \mathbb{Z}$), which corresponds to $\mathcal{A}_n[k]\lrangle{k}$ (for all $k\in \mathbb{Z}$) under this equivalence of categories. Here, $[k]$ and $\lrangle{k}$ denote cohomological and grading shifts respectively.

A similar statement is true for $\lbar{\mathcal{A}}_n\hphMod^{\gr,\perf}$.

\begin{lem}
  The category $\lbar{\mathcal{A}}_n\hphMod^{\gr,\perf}$ has a natural weight structure whose weight heart is spanned by finite direct sums of direct summands of $\lbar{\mathcal{A}}_n[k]\lrangle{k}$ for all $k\in \mathbb{Z}$.

  Moreover, the weight heart $\lbar{\mathcal{A}}_n\hphMod^{\gr,\perf,\weightheart}$ is classical and hence the weight complex functor induces an equivalence of categories
  \[
    \wt: \lbar{\mathcal{A}}_n\hphMod^{\gr,\perf} \xrightarrow{\simeq} \Ch^b(\lbar{\mathcal{A}}_n\hphMod^{\gr,\perf,\weightheart}).
  \]
\end{lem}
\begin{proof}
  It is easy to see that the direct summands of $\lbar{\mathcal{A}}_n[k]\lrangle{k}$ generate the category under finite (co)limits. Moreover, as there are no positive $\Ext$'s between them, the weight structure is obtained by invoking~\cite[Theorem 4.3.2.II]{bondarko_weight_2010}. The second part follows by observing that there are no negative $\Ext$'s between these objects either (see also~\cite[\crefnolink{mg:subsubsec:weight_complex_functor}]{ho_revisiting_2022}).
\end{proof}

We are now ready to identify $\wt$ and $\wtilde{\inv}_\theta$.

\begin{prop} \label{prop:wt_cmplx_functor_Ch_u_gr_G_vs_wtilde(inv)_theta}
  The functor $\wtilde{\inv}_\theta$ induces an equivalence of categories $\lbar{\inv}_\theta$ as given in the following commutative diagram
  \[
    \begin{tikzcd}[column sep=tiny]
      \mathcal{A}_n\hphMod^{\gr,\perf} \ar{dr}[swap]{\wt} \ar{rr}{\wtilde{\inv}_\theta} && \lbar{\mathcal{A}}_n\hphMod^{\gr,\perf}. \\
      & \Ch^b(\ho \mathcal{A}_n\hphMod^{\gr,\perf,\weightheart}) \ar{ur}{\simeq}[swap]{\lbar{\inv}_\theta}
    \end{tikzcd}
  \]
\end{prop}
\begin{proof}
  The computation in \cref{subsubsec:Koszul_restr_to_A_n-Mod^perf} shows that $\wtilde{\inv}_\theta(\mathcal{A}_n) \simeq \lbar{\mathcal{A}}_n$ and hence, $\wtilde{\inv}_\theta$ is weight exact. By \cref{cor:restr_weight_exact_functors}, $\wtilde{\inv}_\theta$ fits into the following commutative diagram
  \[
    \begin{tikzcd}
      \mathcal{A}_n\hphMod^{\gr,\perf} \ar{d}{\wt} \ar{r}{\wtilde{\inv}_\theta} & \lbar{\mathcal{A}}_n\hphMod^{\gr,\perf} \ar{d}{\wt}[swap]{\simeq} \\
      \Ch^b(\ho \mathcal{A}_n\hphMod^{\gr,\perf,\weightheart}) \ar{r}{\lbar{\inv}_\theta} & \Ch^b(\lbar{\mathcal{A}}_n\hphMod^{\gr,\perf,\weightheart}).
    \end{tikzcd}
  \]
  It remains to show that $\lbar{\inv}_\theta$ is an equivalence of categories.

  The diagram above is determined by the following commutative diagram involving the weight hearts
  \[
    \begin{tikzcd}
      \mathcal{A}_n\hphMod^{\gr,\perf,\weightheart} \ar{d}{\pi_0} \ar{r}{\wtilde{\inv}_\theta^{\weightheart}} & \lbar{\mathcal{A}}_n\hphMod^{\gr,\perf,\weightheart} \ar[equal]{d} \\
      \ho \mathcal{A}_n\hphMod^{\gr,\perf,\weightheart} \ar{r}{\lbar{\inv}_\theta^{\weightheart}} & \lbar{\mathcal{A}}_n\hphMod^{\gr,\perf,\weightheart}.
    \end{tikzcd}
  \]
  By \cref{cor:restr_weight_exact_functors}, it suffices to show that $\lbar{\inv}_\theta^{\weightheart}$ is an equivalence of categories. But this is clear since both
  \[
    \cHom_{\mathcal{A}_n\hphMod^{\gr,\perf,\weightheart}}(\mathcal{A}_n, \mathcal{A}_n[k]\lrangle{k}) \to \cHom_{\ho \mathcal{A}_n\hphMod^{\gr,\perf,\weightheart}}(\pi_0 \mathcal{A}_n, \pi_0 \mathcal{A}_n[k]\lrangle{k})
  \]
  and
  \[
    \cHom_{\mathcal{A}_n\hphMod^{\gr,\perf,\weightheart}}(\mathcal{A}_n, \mathcal{A}_n[k]\lrangle{k}) \to \cHom_{\lbar{\mathcal{A}}_n\hphMod^{\gr,\perf,\weightheart}}(\lbar{\mathcal{A}}_n, \lbar{\mathcal{A}}_n[k]\lrangle{k})
  \]
  realize the RHS by killing off $\uline{\theta}$ from the LHS. Thus, we are done.
\end{proof}

\subsubsection{Restricting to the nilpotent cone}
\label{subsubsec:restricting_nilp_cone}
The discussion above has a geometric interpretation in terms of character sheaves and sheaves on the nilpotent cone of $G$. It gives an explanation for why the nilpotent cone should appear at all in~\cite{trinh_hecke_2021}. The materials presented here are of independent interest and are not needed anywhere else in the paper. We include them here mostly for the sake of completeness. The reader should feel free to skip to \cref{subsubsec:identifying_lbar_Gamma}.

Let $\mathcal{N}$ denote the nilpotent cone of the group $G = \GL_n$, equipped with the conjugation action of $G$. Let $\iota_\mathcal{N}: \mathcal{N} \hookrightarrow G$ denote the closed embedding, then we have the following functor
\[
  \iota_\mathcal{N}^*: \Ch^{\unip,\gr}_G \to \Shv_{\gr,c}(\mathcal{N}/G).
\]
We let $\lbar{\Spr}^\gr_G \defeq \iota_\mathcal{N}^* \Spr^\gr_G$ denote the (graded) Springer sheaf.

The starting point is the following formality result of Rider, which was formulated in a different language.
\begin{thm}[{\cite[Theorem 7.9]{rider_formality_2013}}]
  The functor $\cHom^{\gr}_{\Shv_{\gr,c}(\mathcal{N}/G)}(\lbar{\Spr}^\gr_G, -)$ induces an equivalence of $\Vect^{\gr,c}$-categories
  \[
    \Shv_{\gr,c}(\mathcal{N}/G) \simeq \cEnd^\gr(\lbar{\Spr}^\gr_G)\hphMod^{\gr, \perf} \simeq \lbar{\mathcal{A}}_n\hphMod^{\gr,\perf}.
  \]
\end{thm}

The next input is a theorem of Trinh which compares the graded derived endomorphism rings of $\Spr^\gr_G$ and $\lbar{\Spr}^\gr_G$.

\begin{thm}[{\cite[Theorem 1]{trinh_hecke_2021}}]
  The morphism of algebras
  \[
    \cEnd^\gr(\Spr^\gr_G) \to \cEnd^\gr(\iota_\mathcal{N}^* \Spr^\gr_G) \simeq \cEnd^\gr(\lbar{\Spr}^\gr_G)
  \]
  identifies with the natural quotient $\mathcal{A}_n \to \lbar{\mathcal{A}}_n$.
\end{thm}

The following is thus a direct consequence of the two theorems above and \cref{thm:explicit_Char_type_A}.
\begin{cor} \label{cor:restr_nilp_cone_vs_inv_theta}
  We have a commutative diagram
  \[
    \begin{tikzcd}
      \Ch^{\unip,\gr}_G \ar{r}{\iota_\mathcal{N}^*} \ar{d}{\simeq}[swap]{\cHom^\gr_{\Ch^{\unip,\gr}_G}(\Spr^\gr_G,-)} & \Shv_{\gr,c}(\mathcal{N}/G) \ar{d}[swap]{\simeq}{\cHom^\gr_{\Shv_{\gr,c}(\mathcal{N}/G)}(\lbar{\Spr}^\gr_G,-)} \\
      \mathcal{A}_n\hphMod^{\gr,\perf} \ar{r}{\wtilde{\inv}_\theta} & \lbar{\mathcal{A}}_n\hphMod^{\gr,\perf}.
    \end{tikzcd}
  \]
\end{cor}

\subsubsection{}
The category $\Shv_{\gr, c}(\mathcal{N}/G)$ has a natural weight structure whose weight heart is spanned by finite direct sums of summands of $\lbar{\Spr}^\gr_G[k]\lrangle{k}, k\in \mathbb{Z}$. This corresponds to the natural weight structure on $\lbar{\mathcal{A}}_n\hphMod^{\gr,\perf}$ which is spanned by finite direct sums of summands of $\lbar{\mathcal{A}}_n[k]\lrangle{k}, k\in \mathbb{Z}$. Since $\iota_\mathcal{N}^*(\Spr^\gr_G) \simeq \lbar{\Spr}^\gr_G$, $\iota_\mathcal{N}^*$ is weight exact (as can also be seen at the level of $\wtilde{\inv}_\theta$).

\begin{prop} \label{prop:wt_cmplx_functor_Ch_u_gr_G_vs_restr_nilp_cone}
  We have a commutative diagram
  \[
    \begin{tikzcd}[column sep=small]
      \Ch^{\unip,\gr}_G \ar{dr}[swap]{\wt} \ar{rr}{\iota_\mathcal{N}^*} && \Shv_{\gr,c}(\mathcal{N}/G). \\
      & \Ch^b(\ho \Ch^{\unip,\gr,\weightheart}_G) \ar{ur}[swap]{\simeq}
    \end{tikzcd}
  \]
  where all functors are weight-exact.
\end{prop}
\begin{proof}
  This follows directly from \cref{cor:restr_nilp_cone_vs_inv_theta} and the identification of $\wt$ with $\wtilde{\inv}_\theta$ in \cref{prop:wt_cmplx_functor_Ch_u_gr_G_vs_wtilde(inv)_theta}.
\end{proof}

In~\cite{webster_geometric_2017}, the \homflypt{} homology groups are constructed via the chromatographic complex construction, which is a triangulated (as opposed to $\DG$) version of the weight complex functor (see also~\cite[\crefnolink{mg:rmk:weight_complex_functor_explicit}]{ho_revisiting_2022}). In~\cite{trinh_hecke_2021}, the \homflypt{} homology construction is then proved to factor through the nilpotent cone. \cref{prop:wt_cmplx_functor_Ch_u_gr_G_vs_restr_nilp_cone} above puts these two results into context by identifying the weight complex functor itself and the functor $\iota_\mathcal{N}^*$ of restricting to the nilpotent cone, explaining the reason why the nilpotent cone appears. Moreover, the equivalence $\ho \Ch^{\unip,\gr,\weightheart}_G \simeq \Shv_{\gr,c}(\mathcal{N}/G)^{\weightheart}$ formulates precisely the sense in which $\Shv_{\gr,c}(\mathcal{N}/G)$ as the truncated trace of $\Hecke^\gr_n$, as was speculated in~\cite{trinh_hecke_2021}.

\subsubsection{Identifying the functor $\lbar{\Gamma}$}
\label{subsubsec:identifying_lbar_Gamma}
We will now compute $\lbar{\Gamma}$ appearing in \cref{eq:wHH_via_tr} more explicitly: we decompose it into a direct sum (according to the $a$-degrees) of co-representable functors. In other words, the $\alpha$-th piece captures the part of $\Gamma$ that is $\alpha$ away from being pure.

For each integer $\alpha$, consider
\[
  \tilde{\varphi}_\alpha: \mathbb{Z} \to \mathbb{Z}^2
\]
such that
\[
  \tilde{\varphi}_\alpha(X) \defeq \tau_\alpha + \varphi(X) \defeq (2\alpha, \alpha) + (X, X), \qquad X \in \mathbb{Z}.
\]
We define the associated functors $\tilde{\iota}_\alpha, \iota_\alpha: \Vect^\gr \to \Vect^{\gr_{Q'}, \gr_{A'}}$ given by $\tilde{\varphi}_{\alpha,*}$ and $\varphi_*$, both followed by a cohomological shear to the left whose amount is given by the $X$-degree. More precisely, for $(V_X)_{X\in \mathbb{Z}} \in \Vect^\gr$ (note that $T'$ is the cohomological degree of $\Vect^{\gr_{Q'},\gr_{A'}}$, see \cref{subsubsec:AQT_aqt_gradings}),
\begin{align*}
  \tilde{\iota}_\alpha((V_X)_{X\in \mathbb{Z}})_{Q',A'}
   & = \begin{cases}
         V_X\{X\}_{T'}, & (Q', A') = \tilde{\varphi}_\alpha(X), \\
         0,             & \text{otherwise},
       \end{cases} \\
  \iota_\alpha((V_X)_{X\in \mathbb{Z}})_{Q',A'}
   & = \begin{cases}
         V_X\{X\}_{T'}, & (Q', A') = \varphi(X), \\
         0,             & \text{otherwise}.
       \end{cases}
\end{align*}
Note that the shear to the left by $X$ is there to account for the fact that something in graded degree $X$ and cohomological degree $X$ in $\Vect^\gr$ is sent to something of cohomological degree $0$ when we apply the construction in \cref{subsubsec:turning_coh_gr_to_formal_gr} to turn a cohomological grading to a formal grading.

\begin{prop}
  \label{prop:decomposition_lbar_Gamma_a_grading}
  The functor $\lbar{\Gamma}: \lbar{\mathcal{A}}_n\hphMod^{\gr,\perf} \to \Vect^{\gr_{Q'}, \gr_{A'}}$ breaks up into a finite direct sum
  \begin{align*}
    \lbar{\Gamma}
    \simeq \bigoplus_\alpha \tilde{\iota}_\alpha\left( \cHom^\gr_{\lbar{\mathcal{A}}_n\hphMod^{\gr,\perf}}(\wedge^\alpha P \otimes \Qlbar[\uline{x}], -) \right),
  \end{align*}
  where $P \in \Rep \SymGrp_n$ is the permutation representation.
\end{prop}
\begin{proof}
  By \cref{lem:restr_to_wheart_eq_cats}, $\lbar{\Gamma}$ is determined by its restriction to the source's weight heart,
  \[
    \lbar{\Gamma}^{\weightheart}: \lbar{\mathcal{A}}_n\hphMod^{\gr,\perf,\weightheart} \to \Vect^{\theart,\gr_{Q'},\gr_{A'}} \to \Vect^{\gr_{Q'}, \gr_{A'}}.
  \]
  This functor, in turn, is determined by its action on $\lbar{\mathcal{A}}_n$. Under the equivalence of categories in \cref{thm:explicit_Char_type_A}, the constant sheaf corresponds to $\Qlbar[\uline{x},\uline{\theta}]$, from which we obtain the following natural equivalence
  \begin{align*}
     & \alignsep\lbar{\Gamma}^{\weightheart}(\lbar{\mathcal{A}}_n)
    \simeq \Gamma(\mathcal{A}_n)
    \simeq \bigoplus_\alpha \tilde{\iota}_\alpha (\wedge^\alpha P \otimes \Qlbar[\uline{x}])                                                                                                                                \\
     & \simeq \bigoplus_\alpha \tilde{\iota}_\alpha \left( \cHom^\gr_{\lbar{\mathcal{A}}_n\hphMod^{\gr,\perf,\weightheart}}(\wedge^\alpha P \otimes \Qlbar[\uline{x}], \Qlbar[\uline{x}]\rotimes \Qlbar[\SymGrp_n]) \right) \\
     & \simeq \bigoplus_\alpha \tilde{\iota}_\alpha \left( \cHom^\gr_{\lbar{\mathcal{A}}_n\hphMod^{\gr,\perf,\weightheart}}(\wedge^\alpha P \otimes \Qlbar[\uline{x}], \lbar{\mathcal{A}}_n) \right),
  \end{align*}
  where we have implicitly picked an isomorphism of $\SymGrp_n$-representations $P \simeq P^\vee$. Moreover, the factor $(2\alpha, \alpha)$ in $\tilde{\varphi}_\alpha$ is to account for the degrees of $\uline{\theta}$ in $\mathcal{A}_n$ which are not there anymore in $\wedge^\alpha P$: these are precisely the $Q',A'$-degrees of homogeneous wedges of $\uline{\theta}$ of order $\alpha$. Thus, we have an equivalence of functors
  \[
    \lbar{\Gamma}^{\weightheart} \simeq \bigoplus_\alpha \tilde{\iota}_\alpha \left( \cHom^\gr_{\lbar{\mathcal{A}}_n\hphMod^{\gr,\perf,\weightheart}}(\wedge^\alpha P \otimes \Qlbar[x], -) \right) : \lbar{\mathcal{A}}_n\hphMod^{\gr,\perf,\weightheart} \to \Vect^{\theart,\gr_{Q'},\gr_{A'}}.
  \]

  But now, the second functor has an extension to the whole of $\lbar{\mathcal{A}}_n\hphMod^{\gr,\perf}$, which is given by
  \[
    \bigoplus_\alpha \tilde{\iota}_\alpha \left( \cHom^\gr_{\lbar{\mathcal{A}}_n\hphMod^{\gr,\perf}}(\wedge^\alpha P \otimes \Qlbar[x], -) \right) : \lbar{\mathcal{A}}_n\hphMod^{\gr,\perf} \to \Vect^{\gr_{Q'},\gr_{A'}}.
  \]
  The proof thus concludes by \cref{lem:restr_to_wheart_eq_cats}.
\end{proof}

\begin{defn} \label{defn:lbar_Gamma_a_wHH_a_Gamma_a}
  We define
  \begin{align*}
    \lbar{\Gamma}_\alpha \defeq \cHom^\gr_{\lbar{\mathcal{A}}_n\hphMod^{\gr,\perf}}(\Qlbar[\uline{x}] \otimes \wedge^\alpha P, -)
     & : \lbar{\mathcal{A}}_n\hphMod^{\gr,\perf} \to \Vect^{\gr}, \\
    \Gamma_\alpha \defeq \lbar{\Gamma}_\alpha \circ \wt \simeq \lbar{\Gamma}_\alpha \circ \wtilde{\inv}_\theta
     & : \mathcal{A}_n\hphMod^{\gr,\perf} \to \Vect^{\gr},        \\
    \wHH_\alpha \defeq \Gamma_\alpha \circ \tr
     & : \Hecke^\gr_n \to \Vect^{\gr},
  \end{align*}
  Here, $\wt$ is the weight complex functor of $\Ch^{\unip,\gr}_G \simeq \mathcal{A}_n\hphMod^{\gr,\perf}$, which is identified with $\wtilde{\inv}_\theta$ by \cref{prop:wt_cmplx_functor_Ch_u_gr_G_vs_wtilde(inv)_theta}.
\end{defn}

As a direct consequence of the direct sum decomposition of $\lbar{\Gamma}$ in \cref{prop:decomposition_lbar_Gamma_a_grading}, we get the following direct sum decompositions of $\Gamma$ and $\wHH$.

\begin{cor} \label{cor:decomposing_lbar_Gamma_Gamma_wHH}
  We have the following equivalences of functors
  \[
    \lbar{\Gamma} \simeq \bigoplus_\alpha \tilde{\iota}_\alpha \circ \lbar{\Gamma}_\alpha, \quad
    \Gamma \simeq \bigoplus_\alpha \tilde{\iota}_\alpha \circ \Gamma_\alpha, \quad\text{\and}\quad
    \wHH \simeq \bigoplus_\alpha \tilde{\iota}_\alpha \circ \wHH_\alpha.
  \]
\end{cor}

\subsubsection{Gradings}
By construction, for any integer $\alpha$, $\tilde{\iota}_\alpha$ sends a graded vector space of graded degree $X$ and cohomological degree $C$ to an object of $Q',A',T'$-degree $(X+2\alpha, X+\alpha, C-X)$. Written multiplicatively (see also \cref{subsubsec:AQT_aqt_gradings}), this object has degree
\[
  Q'^{X+2\alpha} A'^{X+\alpha} T'^{C-X} = Q'^{2\alpha} A'^\alpha Q'^X A'^X T'^{C-X} = a^\alpha Q^X T^{C-X}.
\]
In other words, the target of $\tilde{\iota}_\alpha$ has $a$-degree $\alpha$. Consequently, $\tilde{\iota}_\alpha \circ \lbar{\Gamma}_\alpha$, $\tilde{\iota}_\alpha \circ \Gamma_\alpha$, and $\tilde{\iota}_\alpha \circ \wHH_\alpha$ all land in the $a$-degree $\alpha$ part. This justifies the notation given in \cref{defn:lbar_Gamma_a_wHH_a_Gamma_a}: $\lbar{\Gamma}_\alpha$, $\Gamma_\alpha$, and $\wHH_\alpha$ are the $a$-degree $\alpha$ parts of $\lbar{\Gamma}$, $\Gamma$, and $\wHH$, respectively.

Observe that the factor $(2\alpha, \alpha)$ in $\tilde{\varphi}_\alpha$ is responsible precisely for the $a$-degree. In what follows, we will treat each $a$-degree separately, and for each $a$-degree, we will ignore the $a$-factor and consider only the $Q,T$-degrees. \cref{prop:decomposition_lbar_Gamma_a_grading} thus has the following reformulation.

\begin{cor} \label{cor:a_term_QT_grading}
  The $a$-degree $\alpha$ part $\lbar{\Gamma}_\alpha$ of $\lbar{\Gamma}$ is given by $\cHom^\gr_{\lbar{\mathcal{A}}_n\hphMod^{\gr,\perf}}(\wedge^\alpha P \otimes \Qlbar[\uline{x}], -) \in \Vect^{\gr}$ where the degree $(X, C)$ corresponds to $(X, C-X)$ in the $Q,T$-grading. Here, $X$ is the graded degree and $C$ is the cohomological degree in $\Vect^{\gr}$.
\end{cor}

\subsubsection{Identifying the functor $\Gamma$}
\label{subsubsec:identifying_Gamma}
Our goal is now to show an analog of \cref{cor:a_term_QT_grading} for the functor $\Gamma$ appearing in \cref{eq:wHH_via_tr}. We already saw in \cref{cor:decomposing_lbar_Gamma_Gamma_wHH} above that $\Gamma$ is decomposed into a direct sum of $\iota_\alpha \circ \Gamma_\alpha$. We will now show that $\Gamma_\alpha$'s are co-representable, except that now, the co-representing objects live in $\mathcal{A}_n\hphMod^{\gr,\coh} \subseteq \mathcal{A}_n\hphMod^\gr \simeq \Ind(\mathcal{A}_n\hphMod^{\gr,\perf})$. In other words, they are co-represented by ind-objects (which happen to be coherent)!

Recall the functor
\[
  \wtilde{\triv}_\theta \defeq \Sym^n V[-1] \otimes \triv_\theta: \lbar{\mathcal{A}}_n\hphMod^{\gr,\perf} \to \mathcal{A}_n\hphMod^{\gr,\coh},
\]
where $\triv_\theta$ sends a perfect $\lbar{\mathcal{A}}_n$ complex to a coherent object where $\uline{\theta}$ act by $0$. This functor has a partially defined right adjoint
\[
  \wtilde{\inv}_\theta \defeq \Sym^n V^\vee[1] \otimes \inv_\theta: \mathcal{A}_n\hphMod^{\gr,\perf} \to \lbar{\mathcal{A}}_n\hphMod^{\gr,\perf},
\]
where $\mathcal{A}_n\hphMod^{\gr,\perf}$ is naturally a full subcategory of $\mathcal{A}_n\hphMod^{\gr,\coh}$. This pair of partial adjoints comes from the standard adjunction pair $\triv_\theta \dashv \inv_\theta$ between the $\Ind$-completed categories $\lbar{\mathcal{A}}_n\hphMod^\gr = \Ind(\lbar{\mathcal{A}}_n\hphMod^{\gr,\perf})$ and $\Ind(\mathcal{A}_n\hphMod^{\gr,\coh})$.

\begin{rmk}[Abuse of notation]
  In what follows, we will abuse notation and implicitly view objects of $\mathcal{A}_n\hphMod^{\gr,\perf}$ as living in $\mathcal{A}_n\hphMod^{\gr,\coh}$ whenever necessary. For example, when $c \in \mathcal{A}_n\hphMod^{\gr,\coh}$ and $p \in \mathcal{A}_n\hphMod^{\gr,\perf}$, we write
  \[
    \cHom^\gr_{\mathcal{A}_n\hphMod^{\gr,\coh}}(c, p) \defeq \cHom^\gr_{\mathcal{A}_n\hphMod^{\gr,\coh}}(c, \iota_{\perf\hookrightarrow \coh} p),
  \]
  where $\iota_{\perf\hookrightarrow \coh}: \mathcal{A}_n\hphMod^{\gr,\perf} \hookrightarrow \mathcal{A}_n\hphMod^{\gr,\coh}$ is the natural inclusion.
\end{rmk}

In view of \cref{cor:decomposing_lbar_Gamma_Gamma_wHH,cor:a_term_QT_grading}, we have the following analog of \cref{prop:decomposition_lbar_Gamma_a_grading} for $\Gamma$.

\begin{prop} \label{prop:compute_Gamma_a}
  The $a$-degree $\alpha$ part $\Gamma_\alpha$ of $\Gamma$ is given by
  \[
    \Gamma_\alpha \simeq \cHom^\gr_{\mathcal{A}_n\hphMod^{\gr,\coh}}(\wedge^\alpha P \otimes \Sym^n V[-1] \otimes \Qlbar[\uline{x}], -): \mathcal{A}_n\hphMod^{\gr,\perf} \to \Vect^{\gr},
  \]
  where the degree $(X, C)$ in $\Vect^\gr$ corresponds to $(X, C-X)$ in the $Q,T$-grading. Here, $X$ is the graded degree and $C$ is the cohomological degree in $\Vect^{\gr}$.
\end{prop}
\begin{proof}
  Using the adjunction $\wtilde{\triv}_\theta \dashv \wtilde{\inv}_\theta$ discussed above, for any $M \in \mathcal{A}_n\hphMod^{\gr,\perf}$, we have
  \begin{align*}
    \Gamma_\alpha(M)
     & \simeq \cHom^\gr_{\lbar{\mathcal{A}}_n\hphMod^{\gr,\perf}}(\wedge^\alpha P \otimes \Qlbar[\uline{x}], \wtilde{\inv}_\theta M) \\
     & \simeq \cHom^\gr_{\mathcal{A}_n\hphMod^{\gr,\coh}}(\wedge^\alpha P \otimes \wtilde{\triv}_\theta(\Qlbar[\uline{x}]), M)       \\
     & \simeq \cHom^\gr_{\mathcal{A}_n\hphMod^{\gr,\coh}}(\wedge^\alpha P \otimes \Sym^n V[-1] \otimes \Qlbar[\uline{x}], M),
  \end{align*}
  and the proof concludes.
\end{proof}

We are now ready to state the main theorem of this subsection, which is a direct consequence of the results proved so far.

\begin{thm} \label{thm:decomposing_wHH}
  The $a$-degree $\alpha$ part $\wHH_\alpha$ of $\wHH$ is given by
  \begin{align*}
    \wHH_\alpha
    &\simeq \cHom^\gr_{\lbar{\mathcal{A}}_n\hphMod^{\gr,\perf}}(\wedge^\alpha P \otimes \Qlbar[\uline{x}], \wtilde{\inv}_\theta\tr(-)) \\
    &\simeq \cHom^\gr_{\mathcal{A}_n\hphMod^{\gr,\coh}} (\wedge^\alpha P \otimes \Sym^n V[-1] \otimes \Qlbar[\uline{x}], \tr(-))
  \end{align*}
  as functors $\Hecke^\gr_n \to \Vect^\gr$, where the degree $(X, C)$ corresponds to $(X, C-X)$ in the $Q,T$-grading, where $X$ is the graded degree and $C$ is the cohomological degree in $\Vect^{\gr}$.
\end{thm}
\begin{proof}
  This follows directly from the computation of $\Gamma_\alpha$ in \cref{prop:compute_Gamma_a} and the fact that $\wHH_\alpha = \Gamma_\alpha \circ \tr$ (see \cref{defn:lbar_Gamma_a_wHH_a_Gamma_a,cor:decomposing_lbar_Gamma_Gamma_wHH}).
\end{proof}

\begin{expl} \label{expl:unlink}
  The $a$-degree $\alpha$ part of the \homflypt{} link homology of the unlink of $n$ components is given by
  \[
    \cHom^\gr_{\lbar{\mathcal{A}}_n\hphMod^{\gr,\perf}}(\wedge^\alpha P \otimes \Qlbar[\uline{x}], \wtilde{\inv}_\theta\tr(1))
    \simeq \cHom^\gr_{\mathcal{A}_n\hphMod^{\gr,\coh}} (\wedge^\alpha P \otimes \Sym^n V[-1] \otimes \Qlbar[\uline{x}], \tr(1))
  \]
  where $1$ is the monoidal unit of $\Hecke^\gr_n$ and where the $Q$ and $T$ gradings are given as in \cref{thm:decomposing_wHH}. Now, recall that $\wtilde{\inv}_\theta\tr(1) \simeq \lbar{\mathcal{A}}_n$. The computation in the proof of \cref{prop:decomposition_lbar_Gamma_a_grading} then shows that the left-hand side of the equivalence above is simply $\wedge^\alpha P \otimes \Qlbar[\uline{x}]$. The associated three-variable polynomial is thus
  \begin{align*}
    \left( \sum_{\alpha=0}^n \binom{n}{\alpha} a^\alpha \right) \left( \sum_{l=0}^\infty Q^{2l} \right)^n = \left( \frac{1+a}{1-Q^2} \right)^n = \left( \frac{1+a}{1-q} \right)^n.
  \end{align*}
\end{expl}

\subsection{\texorpdfstring{\homflypt{}}{HOMFLY-PT} homology via \texorpdfstring{$\Hilb(\CC^2)$}{Hilb(ℂ²)}}
\label{subsec:homflypt_via_Hilb}
We will now finally establish the realization of \homflypt{} homology via Hilbert schemes of points on $\CC^2$. The main point is to transport \cref{thm:decomposing_wHH} to the Hilbert scheme side via Koszul duality and the result of Krug as encapsulated in \cref{sec:char_sheaves_vs_Hilb}.

\subsubsection{Koszul duality}
We will now reformulate \cref{thm:decomposing_wHH} using Koszul duality, which is an equivalence of categories \cref{eq:twisted_KD_equivalence}. We start with the following matching of objects.

\begin{lem} \label{lem:identify_corepresenting_obj_under_KD}
  Under the equivalence of categories \cref{eq:twisted_KD_equivalence}, we have the following matching of objects
  \begin{align*}
    \wtilde{\inv}_\theta^\enh(\wedge^\alpha P \otimes \Sym^n V[-1] \otimes \Qlbar[\uline{x}])
     & \simeq \wedge^\alpha P \otimes \Qlbar[\uline{x},\uline{y}]                                         \\
     & \simeq \wedge^\alpha P \otimes \Sym(V_x[-2] \oplus V_y^\vee) \in \mathcal{B}_n\hphMod^{\gr,\perf}.
  \end{align*}
\end{lem}
\begin{proof}
  Directly from the construction (see also \cref{subsubsec:KD_sends_trivial_to_free}), we have
  \[
    \inv_\theta^\enh(\Qlbar[\uline{x}]) \simeq \Qlbar[\uline{x},\uline{y}].
  \]
  Thus,
  \[
    \wtilde{\inv}_\theta^\enh(\Qlbar[\uline{x}]) \defeq \Sym^n V^\vee[1] \otimes \inv_\theta^\enh(\Qlbar[\uline{x}]) \simeq \Sym^n V^\vee[1] \otimes \Qlbar[\uline{x}, \uline{y}].
  \]
  Since all functors involved are linear over $\Qlbar[\SymGrp_n]\hphMod^{\gr,\perf}$, we get
  \[
    \wtilde{\inv}_\theta^\enh(\wedge^\alpha P \otimes \Sym^n V[-1] \otimes \Qlbar[\uline{x}]) \simeq \wedge^\alpha P \otimes \Qlbar[\uline{x},\uline{y}] \in \mathcal{B}_n\hphMod^{\gr,\perf},
  \]
  and the proof concludes.
\end{proof}

\begin{cor} \label{cor:decomposition_wHH_a_grading_KD}
  The $a$-degree $\alpha$ part $\wHH_\alpha$ of $\wHH$ is given by the cohomology of
  \[
    \wHH_\alpha(R_\beta)
    \simeq \cHom^\gr_{\mathcal{B}_n\hphMod^{\gr,\perf}}(\wedge^\alpha P \otimes \Qlbar[\uline{x},\uline{y}], \wtilde{\inv}_\theta^\enh(\tr(R_\beta))) \in \Vect^{\gr}
  \]
  where the degree $(X, C)$ corresponds to $(X, C-X)$ in the $Q,T$-grading. Here, $X$ is the graded degree and $C$ is the cohomological degree in $\Vect^{\gr}$.
\end{cor}
\begin{proof}
  This is simply \cref{thm:decomposing_wHH} transported to $\mathcal{B}_n\hphMod^{\gr,\perf}$ using the equivalence of categories \cref{eq:twisted_KD_equivalence} and the matching of objects done in \cref{lem:identify_corepresenting_obj_under_KD}.
\end{proof}

\subsubsection{\homflypt{} homology via $2$-periodized Hilbert schemes}
We will now relate \homflypt{} homology and $2$-periodized Hilbert schemes of points. Since there are now more than one gradings, we will include them in the notation. Recall that on the $\mathcal{B}_n\hphMod^{\gr,\perf}$ side, there are two sets of gradings: $X,Y$ and $\wtilde{X},\wtilde{Y}$. Here, $X$ is the original grading coming from graded sheaves whereas $Y$ is the extra grading (to be ``canceled out'' by $2$-periodization). Moreover, the $\wtilde{X},\wtilde{Y}$-grading are related to the $X,Y$-grading via a simple change of coordinates. See \cref{subsubsec:introduce_extra_grading,subsubsec:X_Y_gradings_to_X'_Y'}.

\begin{defn} \label{defn:2-periodic_wHH}
  We define
  \[
    \wHH_\alpha^{\per}: \Hecke^\gr_n \to \Qlbar[\beta,\beta^{-1}]\hphMod^{\gr_X,\gr_Y} \eqdef \Vect^{\gr_X,\gr_Y,\per_-} \simeq \Vect^{\gr_{\wtilde{X}}, \gr_{\wtilde{Y}},\per}
  \]
  by $\Qlbar[\beta,\beta^{-1}] \otimes \wHH_\alpha$, where $\beta$ lives in cohomological degree $-2$ and $X,Y$-degree $(0, 1)$, and where $\wHH_\alpha$ is viewed as a complex with two formal gradings by setting the $Y$-degree to $0$ while keeping the $X$-degree the same. See also~\cref{subsubsec:per_-_in_terms_of_wtilde(X)_wtilde(Y)} for the various definitions and grading conventions for $2$-periodization.
\end{defn}

We note that $\wHH_\alpha^{\per}$ and $\wHH_\alpha$ contain the same amount of information. The introduction of $2$-periodization allows to introduce a formal grading $Y$ that is twice the cohomological degree. For example, the $Y=l$ part of $\Ho^0(\wHH_\alpha^{\per}(-))$ is precisely $\Ho^{2l}(\wHH_\alpha(-))$. In this precise sense, we have $Y = C^2$ (in multiplicative notation) if we use $C$ to denote the cohomological degree. All the cohomology groups of $\wHH_\alpha$ can then be recovered by looking only at $\Ho^0$ and $\Ho^1$ of $\wHH^{\per}_\alpha$.

\begin{thm}[{\cite[Conjecture 7.2.(a)]{gorsky_algebra_2022}}] \label{thm:GNR_conj_a}
  We have the following equivalence of functors from $\Hecke^\gr_n \to \Vect^{\gr_{\wtilde{X}}, \gr_{\wtilde{Y}},\per}$
  \[
    \wHH^{\per}_\alpha \simeq \cHom^{\gr_{\wtilde{X}}, \gr_{\wtilde{Y}}}_{\Perf(\Hilb_n/\Gm^2)^{\per}}((\wedge^\alpha \mathcal{T})^{\per}, \Psi^{\per}(\wtilde{\inv}_\theta^\enh(\tr(-))^{\per})),
  \]
  where $\mathcal{T}$ is the tautological bundle on $\Hilb_n$. The $\wtilde{X},\wtilde{Y}$-grading matches with the $q,t$-grading as follows: $\wtilde{X} = q$, $\wtilde{Y} = \sqrt{t}$. Moreover, the cohomological degree $C$ on the RHS corresponds to $\sqrt{qt}$. Note also that for any $R \in \Hecke^\gr_n$, $\Psi^{\per}(\wtilde{\inv}_\theta^\enh(\tr(R))^{\per}) \in \Perf(\Hilb_n/\Gm^2)^{\per}_{\Hilb_{n,\uline{x}}}$, i.e., it is supported along the $x$-axis.

  In particular, for any $R_\beta \in \Hecke^\gr_n$ associated to a braid $\beta$, there exists a natural
  \[
    \mathcal{F}_\beta \defeq \Psi^{\per}(\wtilde{\inv}_\theta^\enh(\tr(R_\beta))^{\per}) \in \Perf(\Hilb_n/\Gm^2)^{\per}_{\Hilb_{n,\tilde{\uline{x}}}}
  \]
  such that the $a$-degree $\alpha$ component of the \homflypt{} homology of $\beta$ is given by
  \[
    \cHom^{\gr_{\wtilde{X}},\gr_{\wtilde{Y}}}_{\Perf(\Hilb_n/\Gm^2)^{\per}} (\wedge^\alpha \mathcal{T}^{\per}, \mathcal{F}_\beta).
  \]
\end{thm}
\begin{proof}
  For any $R \in \Hecke^\gr_n$, we have the following natural equivalences
  \begin{align*}
    \wHH^{\per}_\alpha(R)
     & \defeq \Qlbar[\beta,\beta^{-1}] \otimes \wHH_\alpha(R)                                                                                                                                                                                                               \\
     & \simeq \Qlbar[\beta,\beta^{-1}] \otimes \cHom^{\gr_X}_{\mathcal{B}_n\hphMod^{\gr,\perf}}(\wedge^\alpha P \otimes \Qlbar[\uline{x}, \uline{y}], \wtilde{\inv}_\theta^\enh(\tr(R))) \tag{\cref{cor:decomposition_wHH_a_grading_KD}}                                    \\
     & \simeq \cHom^{\gr_X,\gr_Y}_{{}^\Rightarrow(\mathcal{B}_n\hphMod^{\gr,\perf})}(\wedge^\alpha P \otimes \Qlbar[\uline{x},\uline{y}], \wtilde{\inv}_\theta^\enh(\tr(R))) \tag{\cref{lem:sheared_mod_cat_str_vs_enriched_Hom}}                                           \\
     & \simeq \cHom^{\gr_{\wtilde{X}}, \gr_{\wtilde{Y}}}_{\wtilde{\mathcal{B}}_n^\gr\hphMod^{\gr_{\wtilde{X}},\gr_{\wtilde{Y}}, \per, \perf}}((\wedge^\alpha T)^{\per}, (\wtilde{\inv}_\theta^\enh(\tr(R)))^{\per}) \tag{\cref{cor:sheared_B_n-Mod_vs_2-periodic}}              \\
     & \simeq \cHom^{\gr_{\wtilde{X}}, \gr_{\wtilde{Y}}}_{\Perf(\Hilb_n/\Gm^2)^{\per}}((\wedge^\alpha \mathcal{T})^{\per}, \Psi^{\per}(\wtilde{\inv}_\theta^\enh(\tr(R))^{\per})) \tag{\cref{thm:2-periodized_hilb_vs_ch_sheaves,cor:Krug_match_wedges_taut_2-periodized}}.
  \end{align*}
  Note that as before, we implicitly passed from $X,Y$-grading to $\wtilde{X},\wtilde{Y}$-grading.

  For the matching of gradings, observe that $\wtilde{X}^k \wtilde{Y}^l = X^{2k}Y^k X^{-l} = X^{2k-l}C^{2k}$ corresponds to
  \[
    Q^{2k-l} T^{2k-2k+l} = Q^{2k-l}T^l = Q^{2k} Q^{-l}T^l = q^k (\sqrt{t})^l
  \]
  in the $q,t$-degree. In other words, $\wtilde{X},\wtilde{Y}$-grading corresponds to $q,\sqrt{t}$-grading precisely. Moreover, cohomological degree $1$, which is $C$, corresponds to
  \[
    T = Q Q^{-1} T = \sqrt{qt}.
  \]
\end{proof}

\begin{rmk}
  When decategorified, a term in cohomological degree $1$ of $\wHH_\alpha^{\per}$ contributes a factor of $-\sqrt{qt}$. The sign is there because of the odd cohomological degree.
\end{rmk}

\begin{expl} \label{expl:unlink_periodic}
  We will now rephrase \cref{expl:unlink} in terms of $2$-periodic complexes. By \cref{defn:2-periodic_wHH}, the $a$-degree $\alpha$ part of the $2$-periodized \homflypt{} homology is given by $\wHH_\alpha^{\per} = \Qlbar[\beta, \beta^{-1}]\otimes \wedge^\alpha P \otimes \Qlbar[\uline{x}]$. Since there's no odd cohomological degree, the associated polynomial is extracted from
  \[
    \Ho^0(\wHH_\alpha^{\per}) \simeq \wedge^\alpha P \otimes \Qlbar[\beta \uline{x}],
  \]
  where $\beta \uline{x}$ has cohomological degree $0$ and $X,Y$-degree $(2, 1)$, or equivalently $\wtilde{X},\wtilde{Y}$-degree $(1, 0)$, which is the same as $q,t$-degree $(1, 0)$. The associated polynomial is thus,
  \[
    \left( \sum_{\alpha=0}^n \binom{n}{\alpha} a^\alpha \right) \left( \sum_{l=0}^\infty q^l \right)^n = \left( \frac{1+a}{1-q} \right)^n.
  \]
\end{expl}

\subsection{Matching the actions \texorpdfstring{of $\Qlbar[\uline{x}]^{\SymGrp_n}$}{}} \label{subsec:match_actions_Qlbar[x]^{S_n}}
In \cite[\S5.1]{gorsky_algebra_2022}, an action of the variables $\uline{x}$ (or equivalently, $\tilde{\uline{x}}$, depending on whether we are working with the sheared or the $2$-periodized version) on $\HHH$ is constructed. In particular, symmetric functions on $\uline{x}$ (or equivalently, $\tilde{\uline{x}}$) also act. On the other hand, by construction, the geometric realization of $\wHH^{\per}$ via Hilbert schemes in \cref{thm:GNR_conj_a} automatically equips it with an action of $\Qlbar[\tilde{\uline{x}}, \tilde{\uline{y}}]^{\SymGrp_n} = \Ho^0(\Hilb_n, \mathcal{O})$. In particular, it admits an action of $ \Qlbar[\tilde{\uline{x}}]^{\SymGrp_n} \subset \Qlbar[\tilde{\uline{x}}, \tilde{\uline{y}}]^{\SymGrp_n}$.

This subsection is dedicated to showing the following result, whose proof will conclude in \cref{subsubsec:complete_proof:thm:GNR_conj_b}.

\begin{thm}[{part of~\cite[Conjecture 7.2.(b)]{gorsky_algebra_2022}}] \label{thm:GNR_conj_b}
  The two actions above coincide.
\end{thm}

The proof is a simple matter of chasing through the construction. To keep the main ideas evident, we will elide the difference between sheared and $2$-periodic versions. For example, we will pass seamlessly between $\uline{x},\uline{y}$ and $\tilde{\uline{x}},\tilde{\uline{y}}$.

The argument can be most conceptually and conveniently phrased in terms of module categories. For example, instead of saying that a (graded) commutative ring $R$ acts on objects of a ($\Vect^\gr$-module) category $\mathcal{C}$, we formulate everything in terms of the symmetric monoidal category $\mathcal{O} = R\hphMod$ (or $\mathcal{O} = R\hphMod^\gr$) acting on $\mathcal{C}$ itself. Then, the core of the argument revolves around observing that all functors involved are linear over $\mathcal{O}$.

\subsubsection{Module categories and actions of a commutative ring}
We start with some generalities regarding module categories and actions of commutative rings.

Let $\mathcal{O}$ and $\mathcal{O}'$ be compactly generated rigid symmetric monoidal categories equipped with a symmetric monoidal functor $F: \mathcal{O} \to \mathcal{O}'$ and a right adjoint $G$, which is necessarily right-lax symmetric monoidal. The symmetric monoidal functor $F$ equips $\mathcal{O}'$ with the structure of an $\mathcal{O}$-module category. $F$ can thus be written as $F = -\otimes 1_{\mathcal{O}'}$. The functor $G$ can thus be written as $G \simeq \cHom^{\mathcal{O}}_{\mathcal{O}'}(1_{\mathcal{O}'}, -)$, where the superscript $\mathcal{O}$ in $\cHom$ denotes the $\mathcal{O}$-enriched $\Hom$. See~\cite[\crefnolink{mg:subsubsec:enriched_Hom}]{ho_revisiting_2022} for a quick review on enriched $\Hom$-spaces.

\begin{lem} \label{lem:O-ModCat_to_mod_str_on_Hom}
  Let $\mathcal{O}$ and $\mathcal{O}'$ be as above, and $\mathcal{C}$ an $\mathcal{O}'$-module category. Then, for any $c_1, c_2 \in \mathcal{C}$, $\cHom^{\mathcal{O}}_\mathcal{C}(c_1,c_2)$ has a natural $\cEnd^\mathcal{O}_{\mathcal{O}'}(1_{\mathcal{O}'})$-module structure. Moreover, this module structure is compatible with
  \begin{enumerate}
    \item \label{item:lem_O-ModCat_to_mod_str_on_Hom:func_C} $\mathcal{O}'$-linear functors: if $F: \mathcal{C} \to \mathcal{D}$ is a morphism of $\mathcal{O}'$-module categories, then $\cHom^{\mathcal{O}}_\mathcal{C}(c_1, c_2) \to \cHom^{\mathcal{O}}_{\mathcal{D}}(F(c_1), F(c_2))$ is a morphism of modules over $\cEnd^\mathcal{O}_{\mathcal{O}'}(1_{\mathcal{O}'})$;
    \item \label{item:lem_O-ModCat_to_mod_str_on_Hom:func_O}  the restriction of structure via a symmetric monoidal functor $F': \mathcal{O}' \to \mathcal{O}''$: if the $\mathcal{O}'$-module structure on $\mathcal{C}$ comes restricting an $\mathcal{O}''$-module structure, then the $\cEnd^{\mathcal{O}}_{\mathcal{O}'}(1_{\mathcal{O}'})$-module structure on $\cHom^\mathcal{O}_\mathcal{C}(c_1, c_2)$ is obtained by restriction of scalars along the commutative algebra map $\cEnd^{\mathcal{O}}_{\mathcal{O}'}(1_{\mathcal{O}'}) \to \cEnd^{\mathcal{O}}_{\mathcal{O}''}(1_{\mathcal{O}''})$.
  \end{enumerate}
\end{lem}
\begin{proof}
  We have $\cHom^{\mathcal{O}'}_\mathcal{C}(c_1,c_2) \in \mathcal{O}'$ and thus, it naturally has a module structure over $1_{\mathcal{O}'}$. Since $G$ is right-lax symmetric monoidal, $G\cHom^{\mathcal{O}'}_\mathcal{C}(c_1,c_2)$ has a natural module structure over $G(1_{\mathcal{O}'}) \simeq \cEnd^{\mathcal{O}}_{\mathcal{O}'}(1_{\mathcal{O}'})$. But, by adjunction, $G\cHom^{\mathcal{O}'}_\mathcal{C}(c_1,c_2) \simeq \cHom^{\mathcal{O}}_\mathcal{C}(c_1, c_2)$ and we are done.

  The second part regarding various compatibilities is an easy diagram chase.
\end{proof}

\begin{lem} \label{lem:O-ModCat_to_alg_hom}
  Let $\mathcal{O}$, $\mathcal{O}'$, and $\mathcal{C}$ be as above. Then, for any $c$, there is a natural map of algebra objects in $\mathcal{O}$
  \[
    \varphi_c: \cEnd^{\mathcal{O}}_{\mathcal{O}'}(1_{\mathcal{O}'}) \to \cEnd^{\mathcal{O}}_\mathcal{C}(c).
  \]
  Moreover, this map is compatible with $\mathcal{O}'$-linear functors and with restriction of structure via a symmetric monoidal functor $\mathcal{O}' \to \mathcal{O}''$.
\end{lem}
\begin{proof}
  $\cEnd^{\mathcal{O}'}_\mathcal{C}(c)$ is an algebra object in $\mathcal{O}'$ and hence, it receives a natural algebra map from $1_{\mathcal{O}'}$ since $1_{\mathcal{O}'}$ is the initial algebra object. Applying $G$ to this map, we obtained the desire map of algebras in $\mathcal{O}$.

  The second part follows from \cref{lem:O-ModCat_to_mod_str_on_Hom}.
\end{proof}

\begin{rmk}
  In the case where $\mathcal{O} = \Vect$ (resp. $\mathcal{O} = \Vect^\gr$), \cref{lem:O-ModCat_to_alg_hom} states that any object $c \in \mathcal{C}$ admits a natural action (resp. a graded action) of $\cEnd_{\mathcal{O}'}(1_{\mathcal{O}'})$ (resp. $\cEnd^\gr_{\mathcal{O}'}(1_{\mathcal{O}'})$).
\end{rmk}

\begin{cor} \label{cor:O-ModCat_alg_hom_and_action}
  Let $\mathcal{O}$, $\mathcal{O}'$, and $\mathcal{C}$ be as above. Then, for any $c_1, c_2 \in \mathcal{C}$, the action of $\cEnd^\mathcal{O}_{\mathcal{O}'}(1_{\mathcal{O}'})$ on $\cHom^\mathcal{O}_\mathcal{C}(c_1, c_2)$ factors through $\cEnd^{\mathcal{O}}_\mathcal{C}(c_1)^\rev$ and $\cEnd^{\mathcal{O}}_\mathcal{C}(c_2)$. Here, the superscript $\rev$ denotes the reverse multiplication.
\end{cor}
\begin{proof}
  Since $1_{\mathcal{O}'}$ is the initial algebra object in $\mathcal{O}'$, the action of $1_{\mathcal{O}'}$ on $\cHom^{\mathcal{O}'}_\mathcal{C}(c_1, c_2)$ factors through the natural actions of $\cEnd^{\mathcal{O}'}_\mathcal{C}(c_1)^\rev$ and $\cEnd^{\mathcal{O}'}_\mathcal{C}(c_2)$. Applying $G$, we obtain the desired conclusion.
\end{proof}

In what follows, we will apply the statements above to the case where $\mathcal{O}' = R\hphMod^\gr$ for some commutative ring $R$ and $\mathcal{O} = \Vect^\gr$. Moreover, $F$ is the functor of taking free $R$-modules and $G$ is the forgetful functor.

\begin{rmk}
  In the situations that we are interested in, the functor $F$ and all the actions (in various module category structures) are compact-preserving. Note, however, that $G$ does not preserve compactness. Thus, when working with the small full subcategory spanned by compact objects, \cref{lem:O-ModCat_to_mod_str_on_Hom,lem:O-ModCat_to_alg_hom,cor:O-ModCat_alg_hom_and_action} still apply, except that the enriched $\Hom$'s are not necessarily compact.
\end{rmk}

\subsubsection{The action of $\Qlbar[\uline{x}]$ on $\HHH$}
We will now recast the action of $\Qlbar[\uline{x}]$ on $\HHH$ as described in \cite[\S5.1]{gorsky_algebra_2022} in the categorical language above. To start, observe that $\Hecke^\gr_n \simeq \Shv_{\gr,c}(BB\times_{BG} BB)$ has an action of $\Shv_{\gr, c}(BB\times BB) \simeq \Qlbar[\uline{x}]^{\otimes 2}\hphMod^{\gr,\perf}$ where the two factors act on the left and the right, respectively. The usual action of $\Qlbar[\uline{x}]^{\otimes 2}$ on any $R \in \Hecke^\gr_n$ (i.e., any Soergel bimodule) can thus be recovered from \cref{lem:O-ModCat_to_mod_str_on_Hom}. Note also that the commutative diagram
\[
  \begin{tikzcd}
    BB\times_{BG} BB \ar{d} \ar{r} & BB \ar{d} \\
    BB \ar{r} & BG
  \end{tikzcd}
\]
implies that the induced actions of $\Shv_{\gr,c}(BG) \simeq \Co^*_\gr(BG)\hphMod^{\gr,\perf}\simeq \Qlbar[\uline{x}]^{\SymGrp_n}\hphMod^{\gr,\perf}$ via these two actions (on the left and on the right) agree. In particular, the two actions of $\Co^*_\gr(BG)$ on any $R \in \Hecke^{\gr}_n$ obtained from restricting the two actions of $\Co^*_\gr(BB)$ along $\Co^*_\gr(BG) \to \Co^*_\gr(BB)$ agree (see also~\cite[Lemma 5.1]{gorsky_algebra_2022}).

\subsubsection{Reducing to the weight heart} \label{subsubsec:action_x_reduce_weightheart}
By construction, $\HHH$ is obtained by applying Hochschild homology \emph{term-wise} (via the weight complex functor) to an element $R \in \Hecke^\gr_n$. This was formulated geometrically in \cref{subsubsec:geometric_interpretation_homflypt}. Since all the categorical actions are weight exact, for the purpose of describing the action of $\Qlbar[\uline{x}]$ on $\HHH$, we can assume that $R \in \Hecke^{\gr,\weightheart}_n$ and forget about the weight complex functor. The action of $\Qlbar[\uline{x}] \simeq \Co^*_\gr(BB)$ on $\HHH(R)$ can thus be obtained by applying \cref{lem:O-ModCat_to_mod_str_on_Hom} where $\mathcal{O} = \Vect^\gr$, $\mathcal{O}' = \Shv_\gr(BB)^\ren$, $\mathcal{C} = \Shv_\gr(\frac{G}{B})^\ren$, $c_1=\QlbarA_{\frac{G}{B}}$, and $c_2 = p^* R$. Indeed, this is because $\HHH(R)$ is given by (the cohomology of)
\[
  \Co^*_\gr(BB, \pi_* p^* R) \simeq \cHom^\gr_{\Shv_{\gr,c}(\frac{G}{B})}(\QlbarA_{\frac{G}{B}}, p^* R).
\]

\subsubsection{The action of $\Qlbar[\uline{x}]^{\SymGrp_n}$ on $\HHH$}
We can restrict the action of $\Qlbar[\uline{x}]$ on $\HHH$ to $\Qlbar[\uline{x}]^{\SymGrp_n} \subseteq \Qlbar[\uline{x}]$ and obtain an action of $\Qlbar[\uline{x}]^{\SymGrp_n}$, which is the action described in~\cite{gorsky_algebra_2022}. We will now realize this action more geometrically. We start with the following result.

\begin{lem} \label{lem:C^*(BG)-action_via_G/B}
  For $R\in \Hecke^{\gr,\weightheart}_n$, the action of $\Qlbar[\uline{x}]^{\SymGrp_n} = \Co^*_\gr(BG)$ on $\HHH(R) = \Ho^*(\cHom^\gr_{\Shv_{\gr,c}(\frac{G}{B})}(\Qlbar, p^* R))$ is given by the $\Shv_{\gr,c}(BG)$-module structure on $\Shv_{\gr,c}(\frac{G}{B})$ via \cref{lem:O-ModCat_to_mod_str_on_Hom}, where the former acts via pulling back along $\frac{G}{B} \to BB \to BG$.
\end{lem}
\begin{proof}
  By definition, the $\Co^*_\gr(BG)$-action is given by restricting along $\Co^*_\gr(BG) \to \Co^*_\gr(BB)$. The result thus follows by applying \cref{lem:O-ModCat_to_mod_str_on_Hom}.\ref{item:lem_O-ModCat_to_mod_str_on_Hom:func_O} to the case where $\mathcal{O} = \Vect^\gr$, $\mathcal{O}' = \Shv_\gr(BG)^\ren$, $\mathcal{O}'' = \Shv_\gr(BB)^\ren$, $\mathcal{C} = \Shv_{\gr}(\frac{G}{B})^\ren$, $c_1 = \QlbarA_{\frac{G}{B}}$, and $c_2$ is $p^* R$.
\end{proof}

\begin{lem} \label{lem:sym_pol_x_action_via_tr}
  For $R \in \Hecke^{\gr,\weightheart}_n$, the action of $\Co^*_\gr(BG)$ on $\HHH(R)$ described above agrees with the one coming from the $\Shv_{\gr,c}(BG)$-module structure on $\Shv_{\gr,c}(G/G)$ via \cref{lem:O-ModCat_to_mod_str_on_Hom} and the equivalence
  \[
    \HHH(R) = \Ho^*(\cHom^\gr_{\Shv_{\gr,c}(\frac{G}{G})}(\QlbarA_{\frac{G}{G}}, q_* p^* R)) \simeq \Ho^*(\cHom^\gr_{\Ch^{\unip,\gr}_G}(\QlbarA_{\frac{G}{G}}, \tr(R))).
  \]
\end{lem}
\begin{proof}
  Consider the following (non-Cartesian) commutative square
  \[
    \begin{tikzcd}[sep=small]
      & \ar{dl}[swap]{p} \frac{G}{B} \ar{dr}{q} \\
      BB \times_{BG} BB \ar{ddr} && \ar{ddl} \frac{G}{G} \\
      \\
      & BG
    \end{tikzcd}
  \]
  where $p$ and $q$ form the horocycle correspondence. Note that $p^*$ and $q^*$ are $\Shv_{\gr,c}(BG)$-linear. By rigidity of $\Shv_{\gr,c}(BG)$, $p_*$ and $q_*$ are also $\Shv_{\gr,c}(BG)$-linear. But now, we have a natural equivalence
  \begin{align*}
    \cHom^{\Shv_\gr(BG)^\ren}_{\Shv_\gr(\frac{G}{B})^\ren}(\Qlbar, p^* R)
     & \simeq \cHom^{\Shv_\gr(BG)^\ren}_{\Shv_\gr(\frac{G}{B})^\ren}(q^* \Qlbar, p^* R)                                                                                              \\
     & \simeq \cHom^{\Shv_\gr(BG)^\ren}_{\Shv_\gr(\frac{G}{G})^\ren}(\Qlbar, q_* p^* R) \in \Shv_{\gr,c} \left(\frac{G}{G}\right) \subseteq \Shv_\gr\left( \frac{G}{G} \right)^\ren.
  \end{align*}
  The proof concludes by applying $\cHom^\gr_{\Shv_{\gr,c}(\frac{G}{G})}(\Qlbar, -)$ as in \cref{lem:O-ModCat_to_mod_str_on_Hom} and using \cref{lem:C^*(BG)-action_via_G/B}.
\end{proof}

The action above has a more explicit description in terms of \cref{thm:explicit_Char_type_A}. We start with the following result.

\begin{lem} \label{lem:explicit_Char_type_A_refined}
  We have an equivalence of $\Qlbar[\uline{x}]^{\SymGrp_n}\hphMod^{\gr,\perf}$-module categories (or equivalently, of $\Shv_{\gr,c}(BG)$-module categories)
  \[
    \Hom^\gr_{\Ch^{\unip,\gr}_G}(\Spr^\gr_G, -): \Ch^{\unip,\gr}_G \xrightarrow{\simeq} \Qlbar[\uline{x},\uline{\theta}]\rotimes \Qlbar[\SymGrp_n]\hphMod^{\gr,\perf}
  \]
  where the $\Shv_{\gr,c}(BG)$-module structure on the left is inherited from the embedding $\Ch^{\unip,\gr}_G \hookrightarrow \Shv_{\gr,c}\left( \frac{G}{G} \right)$ and where the $\Co^*_\gr(BG)\hphMod^{\gr,\perf}$-module structure on the right comes from the natural morphism of algebras
  \[
    \Qlbar[\uline{x}]^{\SymGrp_n} \simeq \Co^*_\gr(BG) \simeq \cEnd^\gr_{\Shv_{\gr,c}(BG)}(\Qlbar) \to \cEnd^\gr_{\Ch^{\unip,\gr}_G}(\Spr^\gr_G) \simeq \Qlbar[\uline{x},\uline{\theta}] \rotimes \Qlbar[\SymGrp_n], \teq \label{eq:morph_alg_symmetric_to_End_Spr}
  \]
  given by \cref{lem:O-ModCat_to_alg_hom}. Moreover, this morphism of algebra is the obvious one, given by the embedding of symmetric polynomials on $\uline{x}$ to all polynomials.
\end{lem}
\begin{proof}
  Instead of running the proof of \cref{thm:explicit_Char_type_A} relatively over $\Vect^\gr$, we could have run it relatively over $\Shv_\gr(BG)^\ren$. Then, we obtain an equivalence of $\Shv_\gr(BG)^\ren$-module categories (or equivalently, of $\Qlbar[\uline{x}]^{\SymGrp_n}$-module categories)
  \[
    \Ch^{\unip,\gr,\ren}_G \xrightarrow{\simeq} \Qlbar[\uline{x},\uline{\theta}]\rotimes \Qlbar[\SymGrp_n]\hphMod(\Shv_\gr(BG)^\ren) \simeq \Qlbar[\uline{x},\uline{\theta}]\rotimes \Qlbar[\SymGrp_n]\hphMod(\Qlbar[\uline{x}]^{\SymGrp_n}\hphMod),
  \]
  where the category on the right-hand side is formed using the natural morphism of algebras \cref{eq:morph_alg_symmetric_to_End_Spr}, by \cref{cor:O-ModCat_alg_hom_and_action}. Note also that the equivalence of categories in \cref{thm:explicit_Char_type_A} is obtained by further composing with the forgetful functor
  \begin{align*}
     & \alignsep \Qlbar[\uline{x},\uline{\theta}]\rotimes\Qlbar[\SymGrp_n]\hphMod(\Shv_\gr(BG)^\ren) \simeq \Qlbar[\uline{x},\uline{\theta}]\rotimes \Qlbar[\SymGrp_n]\hphMod(\Qlbar[\uline{x}]^{\SymGrp_n}\hphMod) \\
     & \xrightarrow{\text{forgetful}} \Qlbar[\uline{x},\uline{\theta}] \rotimes \Qlbar[\SymGrp_n]\hphMod(\Vect^\gr) \simeq \Qlbar[\uline{x},\uline{\theta}] \rotimes \Qlbar[\SymGrp_n]\hphMod^\gr.
  \end{align*}
  Now, observe that the forgetful functor is in fact an equivalence of categories. The proof of the first part thus concludes.

  We will compute the morphism \cref{eq:morph_alg_symmetric_to_End_Spr} explicitly, as stated in the last sentence of the lemma. Consider the following (non-Cartesian) commutative diagram
  \[
    \begin{tikzcd}
      \frac{B}{B} \ar{d}{h} \ar{r}{q} & \frac{G}{G} \ar{d}{h} \\
      BB \ar{r}{q} & BG.
    \end{tikzcd}
  \]
  Now, applying \cref{lem:O-ModCat_to_alg_hom} and the fact that $\Spr^\gr_G \simeq q_* \Qlbar$, we see that the natural algebra morphism
  \[
    \Co^*_\gr(BG) \simeq \cEnd^\gr_{\Shv_{\gr,c}(BG)}(\Qlbar) \to \cEnd^\gr_{\Shv_{\gr,c}(G/G)}(\Spr^\gr_G)
  \]
  factors as in the following diagram
  \[
    \begin{tikzcd}
      \Co^*_\gr(BG) \ar{r} \ar{d} & \Co^*_{\gr}(BB) \ar{r} \ar{d} & \cEnd^\gr_{\Shv_{\gr,c}(B/B)}(\Qlbar) \ar{r} \ar{d} & \cEnd^\gr_{\Shv_{\gr,c}(G/G)}(\Spr^\gr_G) \ar{d} \\
      \Qlbar[\uline{x}]^{\SymGrp_n} \ar{r} & \Qlbar[\uline{x}] \ar{r} & \Qlbar[\uline{x},\uline{\theta}] \ar{r} & \Qlbar[\uline{x},\uline{\theta}] \rotimes \Qlbar[\SymGrp_n],
    \end{tikzcd}
  \]
  where the maps in the bottom row are the obvious ones. The commutativity of the third square is due to \cref{thm:formality_Spr^gr}.\ref{item:compute_ind_morphism_thm:formality_Spr^gr}. The commutativity of the other two squares is obvious. The desired statement then follows from the commutativity of the outer square in the diagram above.
\end{proof}

We obtain the following refinement of \cref{thm:2-periodized_hilb_vs_ch_sheaves}.

\begin{cor} \label{cor:2-periodized_hilb_vs_ch_sheaves_refined}
  We have an equivalence of $\Qlbar[\tilde{\uline{x}}]^{\SymGrp_n}\hphMod^{\gr_{\wtilde{X}},\gr_{\wtilde{Y}},\per,\perf}$-module categories
  \[
    {}^\Rightarrow\Ch^{\unip,\gr}_G \simeq
    \wtilde{\mathcal{B}}_n^\gr\hphMod^{\gr_{\wtilde{X}},\gr_{\wtilde{Y}},\per}_{\nilp_{\tilde{\uline{y}}}} \xrightarrow[\simeq]{\Psi^{\per}} \QCoh(\Hilb_n/\Gm^2)^{\per}_{\Hilb_{n,\tilde{\uline{x}}}}
  \]
  where the $\Qlbar[\uline{x}]^{\SymGrp_n}\hphMod^{\gr_{\wtilde{X}},\gr_{\wtilde{Y}},\per,\perf}$-module category structure on the first two (resp. last) categories (resp. category) is given by \cref{lem:explicit_Char_type_A_refined} (resp. by the Hilbert--Chow map and forgetting the $\tilde{\uline{y}}$-variables).

  Similarly, we have the corresponding statement for the small category variant by passing to the full subcategories of compact objects.
\end{cor}
\begin{proof}
  \cref{lem:explicit_Char_type_A_refined} above, combined with Koszul duality, implies that we have an equivalence of $\Qlbar[\tilde{\uline{x}}]^{\SymGrp_n}\hphMod^{\gr_{\wtilde{X}},\gr_{\wtilde{Y}},\per,\perf}$-module categories
  \[
    {}^\Rightarrow \Ch^{\unip,\gr, \ren}_G \simeq
    \wtilde{\mathcal{B}}_n^\gr\hphMod^{\gr_{\wtilde{X}},\gr_{\wtilde{Y}},\per}_{\nilp_{\tilde{\uline{y}}}},
  \]
  where on the right, the module structure is induced by the algebra map
  \[
    \Qlbar[\tilde{\uline{x}}]^{\SymGrp_n} \to \Qlbar[\tilde{\uline{x}}, \tilde{\uline{y}}]^{\SymGrp_n} \to \Qlbar[\tilde{\uline{x}},\tilde{\uline{y}}]\rotimes \Qlbar[\SymGrp_n] \eqdef \wtilde{\mathcal{B}}_n^\gr,
  \]
  which is geometrically realized as coming from the natural morphism
  \[
    \mathbb{A}^{2n}/\SymGrp_n \to \mathbb{A}^{2n}\sslash \SymGrp_n \to \mathbb{A}^n \sslash \SymGrp_n.
  \]

  The commutative diagram \cref{eq:iso_spectral_Hilb_stack} implies that $\Psi^{\per}$ is linear over $\Qlbar[\tilde{\uline{x}},\tilde{\uline{y}}]^{\SymGrp_n}\hphMod^{\gr_{\wtilde{X}},\gr_{\wtilde{Y}},\per}$ (see also~\cite[Remark 3.11]{krug_remarks_2018}), where the $\Qlbar[\tilde{\uline{x}},\tilde{\uline{y}}]^{\SymGrp_n}\hphMod^{\gr_{\wtilde{X}},\gr_{\wtilde{Y}},\per}$-module category structure is given by the Hilbert--Chow morphism $\Hilb_n \to \mathbb{A}^{2n}\sslash \SymGrp_n$. In particular, $\Psi^{\per}$ is linear over $\Qlbar[\tilde{\uline{x}}]^{\SymGrp_n}\hphMod^{\gr_{\wtilde{X}},\gr_{\wtilde{Y}},\per}$. Thus, we obtain an equivalence of $\Qlbar[\tilde{\uline{x}}]^{\SymGrp_n}\hphMod^{\gr_{\wtilde{X}},\gr_{\wtilde{Y}},\per}$-module categories
  \[
    \wtilde{\mathcal{B}}_n^\gr\hphMod^{\gr_{\wtilde{X}},\gr_{\wtilde{Y}},\per}_{\nilp_{\tilde{\uline{y}}}} \xrightarrow[\simeq]{\Psi^{\per}} \QCoh(\Hilb_n/\Gm^2)^{\per}_{\Hilb_{n,\tilde{\uline{x}}}}
  \]
  and the proof concludes.
\end{proof}

\subsubsection{Completing the proof of \cref{thm:GNR_conj_b}} \label{subsubsec:complete_proof:thm:GNR_conj_b}
In \cref{lem:sym_pol_x_action_via_tr}, we show that the action of $\Qlbar[\uline{x}]^{\SymGrp_n}$ on $\HHH$ as given in~\cite{gorsky_algebra_2022} comes from the $\Shv_{\gr,c}(BG)$-module category structure on $\Ch^{\unip,\gr}_G$. In \cref{lem:explicit_Char_type_A_refined}, we show that this module category structure is induced by an explicit map of algebras, which is then used in \cref{cor:2-periodized_hilb_vs_ch_sheaves_refined} to show that this module structure is compatible with the one on the Hilbert scheme side, where the structure is induced by the Hilbert--Chow morphism. But now, this structure is responsible for the geometric construction of the action of $\Qlbar[\uline{x}]^{\SymGrp_n}$ on $\HHH$ and the proof concludes.
\qed

\subsection{Support of \texorpdfstring{$\mathcal{F}_\beta$}{ℱᵦ} on \texorpdfstring{$\Hilb_n$}{Hilbₙ}}
\label{subsec:support_F_beta}
We will now show that for a braid $\beta$, the support of $\mathcal{F}_\beta$ on $\Hilb_n$, where $\mathcal{F}_\beta$ is as in \cref{thm:GNR_conj_a}, can be bounded above using the number of components of the link associated to $\beta$. This is the content of~\cite[Conjecture 7.2.(b) and (c)]{gorsky_algebra_2022}. Unlike what is implied over there, however, we do not deduce this as a consequence of \cref{thm:GNR_conj_b}. The obstacle is that \cref{thm:GNR_conj_b} is a statement about global sections of a sheaf whereas the statement we are after is one about the support \emph{before} taking global sections, which is more local in nature. The route we take is thus via the theory of supports as developed in~\cite{benson_local_2008,arinkin_singular_2015}, which is inherently local.

As in the previous subsection, to keep the exposition not heavy, we will elide the difference between the $2$-periodic and sheared versions of the various categories involved, passing seamlessly between, for instance, $\uline{x},\uline{y}$ and $\tilde{\uline{x}},\tilde{\uline{y}}$.

\subsubsection{Subspaces of \texorpdfstring{$\Hilb_n$}{Hilbₙ} and supports of \texorpdfstring{$\mathcal{F}_\beta$}{ℱᵦ}}
We start with some notation regarding various subspaces of $\Hilb_n$ which will appear as supports of $\mathcal{F}_\beta$.

\begin{defn} \label{defn:closed_subset_A^n}
  \begin{enumerate}
    \item For $\beta \in \Br_n$ a braid on $n$ strands, we define $w_\beta \in \SymGrp_n$ to be the corresponding permutation.
    \item For each $w \in \SymGrp_n$, we let $\mathbb{A}^{2n}_{w} \subseteq \mathbb{A}^{2n}$ denote the closed subscheme of $\mathbb{A}^{2n}$ defined by the relations $x_i = x'_{w(i)}$, where $x_1, \dots, x_n, x'_1, \dots, x'_n$ are the coordinates of the $\mathbb{A}^{2n}$.\footnote{Note that $\mathbb{A}^{2n}_w$ has dimension $n$; the superscript is simply to indicate that it is a closed subscheme of $\mathbb{A}^{2n}$.}
    \item We let $\mathbb{A}^n_{w} \defeq \mathbb{A}^{2n}_{w} \times_{\mathbb{A}^{2n}} \mathbb{A}^n$ be the closed subscheme of the diagonal $\mathbb{A}^n$, where $\mathbb{A}^n \to \mathbb{A}^{2n}$ is the diagonal map. Alternatively, if we let $x_i$'s denote the coordinates of the diagonal $\mathbb{A}^n$, then $\mathbb{A}^n_{w}$ is defined by the relations $x_i = x_{w(i)}$.
    \item For each $w\in \SymGrp_n$, we let $\lbar{w}$ denote its conjugacy class, $(\mathbb{A}^n_{\lbar{w}})\sslash \SymGrp_n \subseteq \mathbb{A}^n\sslash \SymGrp_n$ the image of $\mathbb{A}^n_{w}$ in $\mathbb{A}^n\sslash \SymGrp_n$, and $\mathbb{A}^n_{\lbar{w}} \subseteq \mathbb{A}^n$ the preimage of $(\mathbb{A}^n_{\lbar{w}})\sslash \SymGrp_n$ in $\mathbb{A}^n$. In other words, $\mathbb{A}^n_{\lbar{w}}$ is the orbit of $\mathbb{A}^n_w$ under the $\SymGrp_n$ action and $\mathbb{A}^n_{\lbar{w}}\sslash \SymGrp_n$ its GIT quotient.

    \item \label{item:defn:closed_subset_A^n:Hilb_n_w} For each $w \in \SymGrp_n$, we let $\Hilb_{n,\lbar{w}} \subseteq \Hilb_n$ denote the preimage of $\mathbb{A}^n_{\lbar{w}}\sslash \SymGrp_n$ under $\Hilb_n \to \mathbb{A}^{2n}\sslash \SymGrp_n \to \mathbb{A}^{n}\sslash \SymGrp_n$ where the first map is the Hilbert--Chow map and the second map is the projection onto the first factor.

    \item For each $w\in \SymGrp_n$, we let $\Hilb_{n,\tilde{\uline{x}}, \lbar{w}}$ be the closed subscheme of $\Hilb_n$ fitting in the following Cartesian square
          \[
            \begin{tikzcd}
              \Hilb_{n,\tilde{\uline{x}},\lbar{w}} \ar{d} \ar{r} & \Hilb_{n,\tilde{\uline{x}}} \ar{d} \ar{r} & \Hilb_n \ar{d} \\
              \mathbb{A}^n_{\lbar{w}} \sslash \SymGrp_n \ar{r} & \mathbb{A}^n \sslash \SymGrp_n \ar{r} & \mathbb{A}^{2n}\sslash \SymGrp_n,
            \end{tikzcd}
          \]
          where the bottom right horizontal map is induced by $\mathbb{A}^n \simeq \mathbb{A}^n \times \{0\} \to \mathbb{A}^{2n}$. Equivalently, $\Hilb_{n,\tilde{\uline{x}},\lbar{w}} = \Hilb_{n,\tilde{\uline{x}}} \cap \Hilb_{n, \lbar{w}}$.
  \end{enumerate}
  Note that all of these subschemes are stable under the scaling action of $\Gm$. It thus makes sense to talk about $\Gm$-equivariant (quasi-)coherent sheaves on these spaces.
\end{defn}

\begin{rmk}
  Since these schemes are used for the sole purpose of making statements about set-theoretical supports of (quasi-)coherent sheaves, any possible non-reduced or derived structure is not relevant to us. For definiteness, we take the underlying reduced classical schemes if any non-reduced or derived structure is present.
\end{rmk}

With the definitions in place, we are now ready to state the main result of this subsection, whose proof will conclude in \cref{subsubsec:complete_proof:thm:GNR_conj_b_c} below.

\begin{thm}[{\cite[Conjecture 7.2.(b) and (c)]{gorsky_algebra_2022}}] \label{thm:GNR_conj_b_c}
  Let $\beta$ be a braid on $n$ strands and $w_\beta \in \SymGrp_n$ the associated permutation. Then, $\mathcal{F}_\beta \in \Perf(\Hilb_n/\Gm^2)^{\per}_{\Hilb_{n,\tilde{\uline{x}}}}$ (see \cref{thm:GNR_conj_a}) is set-theoretically supported on $\Hilb_{n,\tilde{\uline{x}},\lbar{w}_\beta}$, i.e., $\mathcal{F}_\beta \in \Perf(\Hilb_n/\Gm^2)^{\per}_{\Hilb_{n, \tilde{\uline{x}}, \lbar{w}_\beta}}$.
\end{thm}

\subsubsection{Support in $\DG$-categories}
As in \cref{subsec:match_actions_Qlbar[x]^{S_n}}, the proof is most conceptually explained in terms of module category structures. These structures are naturally in contact with support conditions in the sense of Arinkin--Gaitsgory in~\cite{arinkin_singular_2015}.

More precisely, let $\mathcal{O} = \mathcal{A}\hphMod^{\gr,\perf}$ be a symmetric monoidal category where $\mathcal{A}$ is a graded commutative ($\DG$) algebra and $\mathcal{C}$ an $\mathcal{O}$-module category. Consider the commutative algebra $A = \bigoplus_{n}\Ho^{2n}(\mathcal{A})$ which is doubly graded. Equivalently, $\Spec A$ is equipped with an action of $\Gm^2$. Then, for any $\Gm^2$-invariant closed subset $Y$ of $\Spec A$,~\cite{arinkin_singular_2015} defines the full subcategory $\mathcal{C}_Y$ of $\mathcal{C}$ consisting of objects supported along $Y$.

The algebra $\mathcal{A}$ in our cases are of the form $\Co^*_\gr(BB\times BB) \simeq \Qlbar[\uline{x},\uline{x}']$, $\Co^*_\gr(BB)\simeq \Qlbar[\uline{x}]$, and $\Co^*_\gr(BG) \simeq \Qlbar[\uline{x}]^{\SymGrp_n}$ which are pure and concentrate in only even degrees. The two $\Gm$ actions thus coincide and moreover, the algebra $A$ is simply a shear of $\mathcal{A}$. We can therefore ignore one of the gradings and the closed subsets $Y$ we will consider are simply $\Gm$-invariant closed subsets of $\mathbb{A}^{2n}$, $\mathbb{A}^n$, and $\mathbb{A}^n \sslash \SymGrp_n$. These are precisely the ones that appear in \cref{defn:closed_subset_A^n}.

\begin{rmk}
  Since these are the only cases that we are interested in, we will assume that all of our commutative graded $\DG$-algebras $\mathcal{A}$ used to study supports are pure and concentrated in even degrees.
\end{rmk}

\begin{rmk}
  Strictly speaking, \cite{arinkin_singular_2015} works with big categories whereas we formulate everything using small categories. This is not a problem, however, since all of our categories are compactly generated and all actions are compact preserving. Alternatively, we could also just formulate everything using big categories by working with, for example, $\Hecke^{\gr,\ren}_n$ and $\Shv_{\gr}(BG)^\ren$, etc. We chose not to do so since, for example, $\Hecke^\gr_n$ is a much more familiar object than $\Hecke^{\gr,\ren}_n$.
\end{rmk}

\subsubsection{Supports of Rouquier complexes}
We will now formulate the support conditions satisfied by Rouquier complexes. Recall that for a braid $\beta$ on $n$ strands, we use $R_\beta \in \Hecke^\gr_n$ to denote the associated Rouquier complex. As seen above, $\Hecke^\gr_n \defeq \Shv_{\gr,c}(BB\times_{BG} BB)$ has a natural module category structure over $\Shv_{\gr,c}(BB\times BB) \simeq \Qlbar[\uline{x},\uline{x}']\hphMod^{\gr,\perf}$. By~\cite[\S3.5 and E.3.2]{arinkin_singular_2015}, it makes sense to talk about the support of any element in $\Hecke^\gr_n$ relative to $\Qlbar[\uline{x},\uline{x}']\hphMod^{\gr,\perf}$. More precisely, given any conical closed subset $Y \subset \mathbb{A}^{2n}$ (with respect to the usual scaling action of $\Gm$), we can talk about the full subcategory
\[
  (\Hecke^{\gr}_n)_Y \simeq \Hecke^\gr \otimes_{\Shv_{\gr,c}(BB\times BB)} \Shv_{\gr,c}(BB\times BB)_Y \simeq \Hecke^\gr \otimes_{\Qlbar[\uline{x},\uline{x}']\hphMod^{\gr,\perf}} \Qlbar[\uline{x},\uline{x}']\hphMod^{\gr,\perf}_Y \hookrightarrow \Hecke^\gr_n
\]
of objects supported on $Y$.

The main computational input is the following result.
\begin{prop} \label{prop:support_R_beta}
  Let $\beta \in \Br_n$ be a braid on $n$ strands. Then, $R_\beta \in \Hecke^\gr_n$ is supported on $\mathbb{A}^{2n}_{w_\beta}$ with respect to the action of $\Shv_{\gr,c}(BB\times BB)$ on $\Hecke^\gr_n$.
\end{prop}
\begin{proof}
  This is essentially~\cite[Theorem 5.2]{gorsky_algebra_2022} but formulated in a more ``local'' way. For the sake of completeness, let us also sketch a more geometric proof, which follows the usual strategy. Namely, we decompose $\beta$ into a product of simple crossings.

  We thus start with the case of a single simple crossing $\beta_i$. Then $w_{\beta_i}$ is a simple permutation. Consider the stack $B\backslash B w_{\beta_i} B/B$. Observe that
  \[
    \Co^*_\gr(B\backslash Bw_{\beta_i} B/B) \simeq \Qlbar[\uline{x}, \uline{x}']/(x_k - x_{w_{\beta_i}(k)}),
  \]
  and hence,
  \[
    \Shv_{\gr,c}(B\backslash B w_{\beta_i} B/B) \simeq \Qlbar[\uline{x}, \uline{x}']/(x_k - x_{w_{\beta_i}(k)})\hphMod^{\gr,\perf}.
  \]
  Consequently, any object in $\Shv_{\gr, c}(B\backslash B w_{\beta_i} B/B)$ is supported on $\mathbb{A}^{2n}_{w_{\beta_i}}$, relative to
  \[
    \Co^*_\gr(BB)\otimes \Co^*_\gr(BB)\hphMod^{\gr,\perf} \simeq \Shv_{\gr,c}(BB\times BB).
  \]
  Now, let $j_{w_{\beta_i}}: B\backslash B w_{\beta_i} B/B \to B\backslash G/B$. Then by rigidity of $\Shv_{\gr,c}(BB\times BB)$,
  \[
    j_{w_{\beta_i},*}, j_{w_{\beta_i},!}: \Shv_{\gr,c}(B\backslash B w_{\beta_i} B/B) \to \Shv_{\gr,c}(B\backslash G/B)
  \]
  are linear over $\Shv_{\gr,c}(BB\times BB)$ and hence, they both preserve supports. In particular, the Rouquier complexes $R_{\beta_i}$ and $R_{\beta_i^{-1}}$ are both supported on $\mathbb{A}^{2n}_{w_{\beta_i}}$.

  For a general braid $\beta$, we write $\beta = \beta_1 \dots \beta_m$. Then $R_\beta$ is obtained using the usual convolution diagram for Hecke categories. The proof concludes by applying~\cite[Proposition 3.5.9]{arinkin_singular_2015} to the convolution diagram realizing $R_\beta$ as a product of $R_{\beta_i}$'s.
\end{proof}

\subsubsection{Support of \texorpdfstring{$\tr(R_\beta)$}{tr(Rᵦ)}}
The rest of the proof of \cref{thm:GNR_conj_b_c} is straightforward and follows essentially the same strategy as the one used in \cref{subsec:match_actions_Qlbar[x]^{S_n}}. Namely, it amounts to transporting the support condition on $R_\beta$ established in \cref{prop:support_R_beta} around to $\Ch^{\unip,\gr}_G$ and then to $\Hilb_n$. We will now consider the first part.

\begin{cor} \label{cor:support_tr(R_beta)}
  Let $R_\beta$ be the Rouquier complex associated to a braid $\beta$. Then, the support of $\tr(R_\beta) \in \Ch^{\unip,\gr}_G$, relative to $\Shv_{\gr,c}(BG) \simeq \Qlbar[\uline{x}]^{\SymGrp_n}\hphMod^{\gr,\perf}$, lies inside $\mathbb{A}^n_{\lbar{w}}\sslash \SymGrp_n$.
\end{cor}
\begin{proof}
  Recall that $\tr=q_!p^*$, where $p$ and $q$ are given in the following diagram
  \[
    \begin{tikzcd}
      BB\times_{BG} BB \ar{d} & \ar{l}[swap]{p} \frac{G}{B} \ar{d} \ar{r}{q} & \frac{G}{G} \ar{d} \\
      BB\times BB & \ar{l} BB \ar{r} & BG.
    \end{tikzcd}
  \]
  Applying~\cite[Proposition 3.5.9]{arinkin_singular_2015} and the fact that supports are preserved under functors compatible with module category structures, we see that $p^* R_\beta$ has support in $\mathbb{A}^n_{w_\beta} \defeq \mathbb{A}^{2n}_{w_\beta} \times_{\mathbb{A}^{2n}} \mathbb{A}^n$ relative to $\Shv_{\gr,c}(BB)$, as a consequence of \cref{prop:support_R_beta}. Applying~\cite[Proposition 3.5.9]{arinkin_singular_2015} again, we see that $p^* R_\beta$ has support $\mathbb{A}^n_{\lbar{w}_\beta}\sslash \SymGrp_n$ relative to $\Shv_{\gr,c}(BG)$. Now, the proof concludes by using the fact that $q_!$ is $\Shv_{\gr,c}(BG)$-linear.
\end{proof}

\subsubsection{Completing the proof of \cref{thm:GNR_conj_b_c}} \label{subsubsec:complete_proof:thm:GNR_conj_b_c}
We are now ready to complete the proof of \cref{thm:GNR_conj_b_c}. First, recall that in geometric settings, i.e. when all categories involved are $\QCoh$ of Artin stacks of finite type and module category structures are given by pulling back along morphisms of stacks, categorical support in the sense above agrees with the usual notion of set-theoretic support. Thus, it suffices to show that the sheaf $\mathcal{F}_\beta$ is supported on $\mathbb{A}^n_{\lbar{w}_\beta}\sslash \SymGrp_n$ relative to $\mathbb{A}^n \sslash \SymGrp_n$. Indeed, this is equivalent to stating that $\mathcal{F}_\beta$ is set-theoretically supported on $\Hilb_{n,\lbar{w}} \subseteq \Hilb_n$, and hence, also on $\Hilb_{n,\tilde{\uline{x}},\lbar{w}} = \Hilb_{n,\tilde{\uline{x}}} \cap \Hilb_{n, \lbar{w}}$ since $\mathcal{F}_\beta$ is known to have support on $\Hilb_{n,\tilde{\uline{x}}}$ by \cref{thm:GNR_conj_a} already.

But now, by \cref{cor:2-periodized_hilb_vs_ch_sheaves_refined} and, again, the fact that supports are preserved by functors that are compatible with the module category structures, this support condition is equivalent to the one proved in \cref{cor:support_tr(R_beta)} above. The proof thus concludes.
\qed

\section*{Acknowledgements}
Q. Ho would like to thank the organizers of the workshop on \emph{Categorification in Quantum Topology and Beyond} at the Erwin Schrödinger International Institute for Mathematics and Physics (ESI) in 2019, where he learned about the conjecture in~\cite{gorsky_flag_2021}. He would also like to thank the organizers of the American Institute of Mathematics's \emph{Link Homology Research Community} for providing a stimulating (online) research environment.

Q. Ho would like to thank Minh-Tam Quang Trinh for many conversations over the years on \homflypt{} homology, Paul Wedrich for helpful email exchanges at the initial stage of this project, and Oscar Kivinen for many helpful email exchanges (especially regarding the various grading conventions for $\HHH$) and for pointing out the reference~\cite{krug_remarks_2018}. P. Li would like to thank Peng Shan and Yin Tian for many helpful conversations on Soergel bimodules and Hilbert schemes. We thank the anonymous referee for carefully reading the manuscript and providing extensive feedback which greatly helps improve the clarity of the exposition.

Q. Ho is partially supported by Hong Kong RGC ECS grant 26305322 and RGC GRF grant 16304923. P. Li is partially supported by the National Natural Science Foundation of China (Grant No. 12101348).

\printbibliography
\end{document}

